\documentclass[toc=bibnumbered,ngerman,english,titlepage=false]{scrartcl}

\usepackage[utf8]{inputenc}
\usepackage{babel}
\usepackage[T1]{fontenc}
\usepackage{lmodern}
\usepackage{microtype}

\usepackage{amsmath}
\usepackage{amsfonts}
\usepackage{amssymb}
\usepackage{bm}
\usepackage{textcomp}
\usepackage{csquotes}
\usepackage{enumitem}
\usepackage{graphicx}

\usepackage{scrpage2,color}

\numberwithin{equation}{section} 

\newcommand{\N}{\ensuremath{\mathbb{N}}}

\newcommand{\R}{\ensuremath{\mathbb{R}}}

\newcommand{\diverg}{\textup{div}}

\newcommand{\tr}{\textup{tr}}

\newcommand{\Hc}{\ensuremath{\mathcal{H}}}

\newcommand{\Lc}{\ensuremath{\mathcal{L}}}

\newcommand{\Oc}{\ensuremath{\mathcal{O}}}

\setkomafont{disposition}{\normalcolor\bfseries}
 \setkomafont{section}{\large} 
 
\begin{document}
\thispagestyle{plain}

\oddsidemargin 6pt\evensidemargin 6pt\marginparwidth 48pt\marginparsep 10pt

\topmargin -18pt\headheight 12pt\headsep 25pt

\ifx\cs\documentclass \footheight 12pt \fi \footskip 30pt

\textheight 625pt\textwidth 431pt\columnsep 10pt\columnseprule 0pt 
                                                                      
\automark[section]{section}
\clearscrheadings
\chead{\scshape\headmark}
\cfoot{\pagemark}
\pagestyle{scrheadings}

 \renewcommand{\headfont}{\slshape}      
 \renewcommand{\pnumfont}{\upshape}      
 \setcounter{secnumdepth}{4}             
 \setcounter{tocdepth}{4}             

\newtheorem{Definition}{Definition}[section]
\newtheorem{Satz}[Definition]{Satz}
\newtheorem{Lemma}[Definition]{Lemma}
\newtheorem{Korollar}[Definition]{Korollar}
\newtheorem{Corollary}[Definition]{Corollary}
\newtheorem{Bemerkung}[Definition]{Bemerkung}
\newtheorem{Remark}[Definition]{Remark}
\newtheorem{Proposition}[Definition]{Proposition}
\newtheorem{Beispiel}[Definition]{Beispiel}
\newtheorem{Theorem}[Definition]{Theorem}

\begin{center}
\textbf{\Large Convergence of the Allen-Cahn Equation to the Mean Curvature Flow with $90$\textdegree-Contact Angle in 2D}\\
\vspace{0.5cm}
\large\textsc{Helmut Abels and Maximilian Moser}\\
\vspace{0.2cm}
\textit{Fakultät für Mathematik, Universität Regensburg, Universitätsstraße 31, D-93053 Regensburg, Germany}\\
\vspace{0.2cm}
\small E-mail: 
\textsf{helmut.abels@mathematik.uni-regensburg.de}, \textsf{Maximilian1.Moser@mathematik.uni-regensburg.de}\\
\vspace{0.2cm}
\normalsize \today
\end{center}
\begin{abstract}
\textbf{Abstract.}   We consider the sharp interface limit of the Allen-Cahn equation with homogeneous Neumann boundary condition in a two-dimensional domain $\Omega$, in the situation where an interface has developed and intersects $\partial\Omega$. Here a parameter $\varepsilon>0$ in the equation, which is related to the thickness of the diffuse interface, is sent to zero. The limit problem is given by mean curvature flow with a $90^\circ$-contact angle condition and convergence using strong norms is shown for small times. Here we assume that a smooth solution to this limit problem exists on $[0,T]$ for some $T>0$ and that it can be parametrized suitably. 
With the aid of asymptotic expansions we construct an approximate solution for the Allen-Cahn equation. In order to estimate the difference of the exact and approximate solution we use a spectral estimate for the linearized Allen-Cahn operator.
\end{abstract}

\noindent
\textbf{Mathematics Subject Classification (2010):}\\
Primary: 35K57; 
Secondary  35B25, 
35B36, 
35R37
.\\
\textbf{Key words:} Sharp interface limit, mean curvature flow, contact angle, Allen-Cahn equation
\vspace{0.1in}\\

\tableofcontents

\setcounter{secnumdepth}{4}\setcounter{tocdepth}{4} 
\section{Introduction}\label{sec_intro}\thispagestyle{plain}
Let $\Omega\subseteq\R^N$ be a bounded, smooth and connected domain and $\varepsilon>0$ be a small parameter. For $u^\varepsilon:\overline{\Omega}\times[0,T]\rightarrow\R$ we consider the Allen-Cahn equation with homogeneous Neumann boundary condition
\begin{alignat}{2}\label{eq_AC1}
\partial_tu^\varepsilon-\Delta u^\varepsilon+\frac{1}{\varepsilon^2}f'(u^\varepsilon)&=0&\qquad&\text{ in }Q_T:=\Omega\times[0,T],\\\label{eq_AC2}
\partial_{\nu}u^\varepsilon&=0&\qquad&\text{ on }\partial Q_T:=\partial\Omega\times[0,T],\\
u^\varepsilon|_{t=0}&=u_0^\varepsilon&\qquad&\text{ in }\Omega,\label{eq_AC3}
\end{alignat} 
where $f$ is a suitable smooth double well potential with wells of equal depth, e.g. $f(u)=\frac{1}{8}(1-u^2)^2$. The precise assumptions are
\[
f\in C^\infty(\R),\quad f'(\pm1)=0,\quad f''(\pm1)>0,\quad \int_{-1}^{u}f'=\int_1^u f'>0\quad \text{ for all }u\in(-1,1).
\] 
The Allen-Cahn equation was originally introduced by Allen and Cahn \cite{AC} to describe the evolution of antiphase boundaries in certain polycrystalline materials. The equation is mainly motivated by the following fact: It is the (by $\frac{1}{\varepsilon}$ accelerated) $L^2$-gradient flow of the Ginzburg-Landau energy
\[
E_\varepsilon(u^\varepsilon)=\int_\Omega\frac{\varepsilon}{2}|\nabla u^\varepsilon|^2+\frac{1}{\varepsilon}f(u^\varepsilon)\,dx.
\] 
From both the equations (fast reaction, slow diffusion) and the energy functional (solutions behave in such a way that the energy is non-increasing in time) one sees that typically after a short time $\Omega$ will be partitioned into subdomains where $u^\varepsilon$ is near $\pm1$ and transition zones develop where $|\nabla u^\varepsilon|$ is large. Heuristically (or in sufficiently smooth cases) the thickness of these zones is proportional to $\varepsilon$, so for $\varepsilon\rightarrow0$ we should obtain sharp interfaces evolving in time. Both from the theoretical and numerical point of view the understanding of such limits is of great importance, cf. e.g. Caginalp, Chen \cite{CaginalpChen} and Caginalp, Chen, Eck \cite{CaginalpChenEck}. For the evolution of the limit interface in our case, formal asymptotic analysis suggests the mean curvature flow, and, if there is boundary contact, a $90$\textdegree-contact angle, cf. Rubinstein, Sternberg, Keller \cite{RSK}. There are several rigorous results supporting this:\\
Using a maximum principle and constructing sub- and supersolutions, there are Chen \cite{ChenGenPropInt} and Evans, Soner, Souganidis \cite{ESS}. For global in time results one has to use some weak formulation, as viscosity in \cite{ESS} and varifold solutions in Ilmanen \cite{Ilmanen}. In the latter case only convergence of a subsequence can be shown. Furthermore, Saez studied the sharp interface limit of the Allen-Cahn equations in \cite{SaezCalcVar}. But there is a severe gap in central part of the proof of the main result. More precisely, on page 372 in \cite{SaezCalcVar} a subset of a sequence space is considered and it is concluded that it is the whole space since it is open and closed. But the topology is not specified such that the sequence space is connected and the conclusions hold. E.g.\ the one used in \cite{SaezIndiana2009} does not work since it is not connected. \\ 
On the other hand there is the general method by de Mottoni and Schatzman \cite{deMS}: In the $\R^N$-case they show convergence using strong norms for times when a smooth solution to mean curvature flow exists. This also works when the interface is strictly contained in $\Omega$. The strategy is as follows:
\begin{enumerate} 
\item Construct an approximate solution $(u_A^\varepsilon)_{\varepsilon\in(0,\varepsilon_0]}$ to \eqref{eq_AC1}-\eqref{eq_AC3} using asymptotic expansions based on the solution to the mean curvature flow. 
\item Estimate the difference of the approximate and exact solution using a spectral estimate for the (at $u_A^\varepsilon(.,t)$ linearized) Allen-Cahn operator $-\Delta+\frac{1}{\varepsilon^2}f''(u_A^\varepsilon(.,t))$ for all $t\in[0,T]$. Here one multiplies the difference of the equations with $u_A^\varepsilon-u^\varepsilon$, integrates and applies a Gronwall-type argument.
\end{enumerate}
This method also yields the typical profile of the solution which is not the case in the other papers mentioned above. It does not make use of a comparison principle and can be applied to other diffuse interface models as well: With general spectrum estimates in Chen \cite{ChenSpectrums} this method was successfully applied to the Cahn-Hilliard equation by Alikakos, Bates, Chen \cite{ABC}, to the phase-field equations by Caginalp, Chen \cite{CaginalpChen}, the mass-conserving Allen-Cahn equation by Chen, Hilhorst, Logak \cite{CHL}, the Cahn-Larché system by Abels, Schaubeck \cite{AbelsSchaubeck} and a Stokes/Allen-Cahn system by Abels, Liu \cite{ALiu}. See also Schaubeck \cite{Schaubeck} for a result on a convective Cahn-Hilliard equation.\\
It should be mentioned that the papers by Chen, Hilhorst, Logak \cite{CHL} and Abels, Liu \cite{ALiu} also yield results for the Allen-Cahn equation with strictly contained interface by simple modifications. The asymptotic expansion is much simpler than in de Mottoni, Schatzman \cite{deMS}: in both papers a \enquote{stretched variable}-technique going back to Caginalp, Chen, Eck \cite{CaginalpChenEck} is used. Moreover, the spectrum estimate as well as the difference estimate in \cite{ALiu} is more refined, thus lowering the number of terms needed in the asymptotic expansion and providing also estimates for first order derivatives. These results have in common that the interface is strictly contained in $\Omega$.\\
The goal of this paper is to apply this strategy to the Allen-Cahn equation in the case of boundary contact. For simplicity we restrict to two dimensions, $N=2$. Given a smooth, compact, embedded and non-closed curve $\Gamma_0\subseteq\overline{\Omega}$ intersecting $\partial\Omega$ at the two endpoints orthogonally, we assume that there is a $T>0$ such that there is a smooth solution $\Gamma:=\bigcup_{t\in[0,T]}\Gamma_t\times\{t\}$ to the mean curvature flow with $90$\textdegree-contact angle condition starting from $\Gamma_0$ and that it can be parametrized suitably with an $X_0:I\times[0,T]\rightarrow\overline{\Omega}$ where $I:=[-1,1]$. See Section \ref{sec_coord} for the precise assumptions.\\
To construct an approximate solution to \eqref{eq_AC1}-\eqref{eq_AC3} we need a boundary layer expansion near the contact points because of the Neumann boundary condition. Therefore, based on a work by Vogel \cite{Vogel}, in Section \ref{sec_coord} we construct suitable curvilinear coordinates $(r,s)$ with values in $(-2\delta,2\delta)\times I$ describing an open neighbourhood of $\Gamma$ in $\overline{\Omega}\times[0,T]$. Here $r$ has the role of a signed distance function and $s$ is like a tangential variable. Then the inner and boundary layer expansions can be compared effectively. For the inner expansion in Subsection \ref{sec_inner_exp}, we proceed similarly as in \cite{CHL} or \cite{ALiu} using the \enquote{stretched variable}-idea, ignoring the Neumann boundary condition. Here some results from ODE-theory are needed, which are summarized in Subsection \ref{sec_ODE}. In Section \ref{sec_bdry_exp}, for the boundary layer expansion at the contact points, the \enquote{stretched variable}-approach is combined with scaling the tangential variable, similarly to Alikakos, Chen, Fusco \cite{ACF}. Therefore we need to solve a model problem on $\R^2_+$ which is done in part \ref{sec_hp}. Using finitely many terms in the expansions, we obtain a candidate for an approximate solution $(u_A^\varepsilon)_{\varepsilon>0}$ in Subsection \ref{sec_approx_sol} by gluing them together.\\
Having the form of the approximate solution, for $\varepsilon\in(0,\varepsilon_0]$ and $\varepsilon_0>0$ small we prove a spectral estimate in Section \ref{sec_SE} which looks similar to that in Abels, Liu \cite{ALiu}, Section 2.7. By localization it will suffice to prove the estimate for parts of $\Omega$ near the contact points; denote them by $\Omega^{B\pm}_t$, $t\in[0,T]$ for now. For the boundary parts an idea from Alikakos, Chen, Fusco \cite{ACF} is used: For $t\in[0,T]$ we construct a suitable approximation $\phi_A^\varepsilon(.,t)$ to the first eigenfunction of $-\Delta+\frac{1}{\varepsilon^2}f''(u_A^\varepsilon(.,t))$ with homogeneous Neumann boundary condition and decompose $H^1(\Omega^{B\pm}_t)=V^\pm_{\varepsilon,t}\oplus(V^\pm_{\varepsilon,t})^\perp$ orthogonally in $L^2(\Omega^{B\pm}_t)$, where $V_{\varepsilon,t}^\pm$ consists of tangential variations $a^\pm(s(.,t))\phi_A^\varepsilon(.,t)$, where $a^\pm$ is in $H^1$ on a fixed suitable interval. Then the associated bilinear form is analyzed on $V_{\varepsilon,t}^\pm\times V_{\varepsilon,t}^\pm, V_{\varepsilon,t}^\pm\times(V_{\varepsilon,t}^\pm)^\perp$ and $(V_{\varepsilon,t}^\pm)^\perp\times(V_{\varepsilon,t}^\pm)^\perp$.\\ 
Using the spectral estimate, we are able to prove a difference estimate in Section \ref{sec_difference} similarly as in Section 5.2 in \cite{ALiu}. Altogether we obtain a convergence result in Section \ref{sec_conv_result} by showing that our approximate solution satisfies the conditions for the difference estimate.\\
\textbf{Acknowlegdement:} The second author gratefully acknowledges support through DFG, GRK 1692 \enquote{Curvature, Cycles and Cohomology} during parts of the work.
\section{Preliminaries}\label{sec_prelim}
\subsection{Notation}
Let $\N$ be the natural numbers and $\N_0:=\N\cup\{0\}$.
The Euclidean norm in $\R^n,n\in\N$ is denoted by $|.|$. The differential operators $\nabla, \diverg$ and $D^2$ always act only on spatial variables. The evaluation of a function $f$ at some point $x$ is often denoted by $f|_x$. Let $X$ be a set and $Y$ a normed space. Then the set of bounded functions $f:X\rightarrow Y$ is denoted by $B(X,Y)$.
\subsection{Curvilinear Coordinates}\label{sec_coord}
Let $\Omega\subseteq\R^N, N=2$ be a smooth, bounded and connected domain with outer unit normal $N_{\partial\Omega}$. In this section we show the existence of a curvilinear coordinate system describing a neighbourhood of an evolving hypersurface\footnote{ Cf. Depner \cite{Depner}, Definition 2.31.} in $\overline{\Omega}$ that meets the boundary $\partial\Omega$ orthogonally and can be parametrized suitably. More precisely let $I:=[-1,1]$, $X_0:I\times[0,T]\rightarrow\overline{\Omega}$ be smooth such that $X_0((-1,1),t)\subseteq\Omega$ and $X_0(.,t)$ is injective for all $t\in[0,T]$ as well as $\partial_sX_0(s,t)\neq 0$ for all $(s,t)\in I\times[0,T]$. Since continuous bijections of compact into Hausdorff topological spaces are homeomorphisms, we know that $X_0(.,t)$ is an embedding and $\Gamma_t:=X_0(I,t)$ is a manifold with boundary for all $t\in[0,T]$. Moreover,
\[
\Gamma:=\bigcup_{t\in[0,T]}\Gamma_t\times\{t\}
\]
is an evolving hypersurface and $\overline{X}_0:=(X_0,\textup{pr}_t):I\times[0,T]\rightarrow\Gamma$ is a homeomorphism.\\ 
We define smooth tangent and normal fields on $\Gamma_t$ by
\[
\vec{\tau}(s,t):=\frac{\partial_sX_0(s,t)}{|\partial_sX_0(s,t)|}\quad\text{ and }\quad\vec{n}(s,t):=
\begin{pmatrix}
0 & 1\\
-1 & 0
\end{pmatrix}\vec{\tau}(s,t)\quad\text{ for all }(s,t)\in I\times[0,T].
\]
Due to Depner, Lemma 2.40 the corresponding normal velocity is given by 
\[
V_{\Gamma_t}(s)=V(s,t)=\vec{n}(s,t)\cdot\partial_tX_0(s,t)\quad\text{ for }\quad(s,t)\in I\times[0,T].
\]
Moreover, let $H(s,t):=H_{\Gamma_t}(s)$ for $(s,t)\in I\times[0,T]$ be the curvature.\\
Additionally, let $\partial\Gamma_t\subseteq\partial\Omega$ and $\Gamma_t\perp\partial\Omega$ for all $t\in[0,T]$. We set $p^\pm(t):=X_0(\pm 1,t)$ and $\overline{p}^\pm(t):=(p^\pm(t),t)$ for $t\in[0,T]$ which are the contact points. Guided by Vogel \cite{Vogel}, Proposition 3.1 we obtain coordinates describing a neighbourhood of $\Gamma$ in $\overline{\Omega}\times[0,T]$:
\begin{Theorem}\label{th_coord}
There is a $\delta>0$ and a smooth $[-\delta,\delta]\times I\times[0,T]\ni(r,s,t)\mapsto X(r,s,t)\in\overline{\Omega}$ with the following properties: $\overline{X}:=(X,\textup{pr}_t)$ is a homeomorphism onto a neighbourhood of $\Gamma$ in $\overline{\Omega}\times[0,T]$, the inverse is smooth and
\[
\Gamma(\eta):=\overline{X}((-\eta,\eta)\times I\times[0,T])
\]
is an open neighbourhood of $\Gamma$ in $\overline{\Omega}\times[0,T]$ for $\eta\in(0,\delta]$. Moreover, $X$ coincides with the usual tubular neighbourhood coordinate system for $s\in[-1+\delta_0,1-\delta_0]$ for some $\delta_0>0$ small, $X|_{r=0}=X_0$ and $X(r,s,t)\in\partial\Omega$ if and only if $s\in\{\pm 1\}$. Furthermore
\[
|\nabla r|^2|_\Gamma=1,\quad \partial_r(|\nabla r|^2\circ\overline{X})|_{r=0}=0\quad\text{ and }\quad\nabla r\cdot\nabla s|_\Gamma=0.
\] 
Moreover, we can choose $\nabla s\circ\overline{X}_0=\partial_sX_0/|\partial_sX_0|^2$ and $\nabla r\circ\overline{X}_0=\vec{n}$. Then it holds $V=-\partial_tr\circ\overline{X}_0$ and $H=-\Delta r\circ\overline{X}_0$.
\end{Theorem}
A similar theorem should hold also in higher dimensions and for other contact angles.
\begin{Remark}\label{th_coord_bem}\upshape
\begin{enumerate}
\item There are unique $Q_T^\pm$ open and connected in $Q_T=\Omega\times[0,T]$ such that $Q_T=Q_T^-\cup Q_T^+\cup\Gamma$ and $\textup{sign}\,r=\pm 1$ on $Q_T^\pm\cap\Gamma(2\delta)$. Moreover, we set
\begin{align*}
\Gamma(\eta,\mu)&:=\overline{X}((-\eta,\eta)\times(-\mu,\mu)\times[0,T]),\\ 
\Gamma^\pm(\eta,\mu)&:=\overline{X}((-\eta,\eta)\times \pm(1-\mu,1)\times[0,T])
\end{align*}
for $\eta\in(0,\delta]$ and $\mu\in(0,1)$. For $t\in[0,T]$ fixed let $\Gamma_t(\eta), \Gamma_t(\eta,\mu)$ and $\Gamma^\pm_t(\eta,\mu)$ be the respective sets intersected with $\R^2\times\{t\}$ and then projected to $\R^2$.
\item For $\psi:\Gamma(\eta)\rightarrow\R, \eta\in(0,\delta]$ sufficiently smooth we define the tangential and normal derivative by $\nabla_\tau\psi:=\partial_s(\psi\circ\overline{X})\circ\overline{X}^{-1}\,\nabla s$ and $\partial_n\psi:=\partial_r(\psi\circ\overline{X})\circ\overline{X}^{-1}$, respectively. In the part of $\Gamma(2\delta)$ where the coordinate system is the orthogonal one, these definitions coincide with the ones in Abels, Liu \cite{ALiu}:
\[
\nabla_\tau\psi=\frac{\nabla s}{|\nabla s|}\frac{\nabla s}{|\nabla s|}\cdot\nabla\psi\quad\text{ and }\quad\partial_n\psi=\nabla r\cdot\nabla\psi\quad\text{ on }\Gamma(\eta,1-\delta_0).
\]
This follows from $\nabla\psi|_{\overline{X}}=\nabla r|_{\overline{X}}\partial_r(\psi\circ\overline{X})+\nabla s|_{\overline{X}}\partial_s(\psi\circ\overline{X})$. For $t\in[0,T]$ fixed and $\psi:\Gamma_t(\eta)\rightarrow\R$ smooth enough, we define $\nabla_\tau\psi$ and $\partial_n\psi$ analogously. In the orthogonal case similar identities as above hold.  
\item We define $J_t(r,s):=\det D_{(r,s)}X(r,s,t)$ for $(r,s,t)\in[-\delta,\delta]\times I\times[0,T]$. From the proof of Theorem \ref{th_coord} it follows that $J_t(r,s)^{-2}=\left[|\nabla r|^2|\nabla s|^2-(\nabla r\cdot\nabla s)^2\right]|_{\overline{X}(r,s,t)}$, in particular $J_t(0,s)=|\partial_sX_0(s,t)|$ for all $(s,t)\in I\times[0,T]$.
\item For convenience we choose $X_0$ such that 
\[
|\partial_sX_0(s,t)|=\frac{L(t)}{2}:=\frac{1}{2}\int_{-1}^1|\partial_sX_0(s,t)|\,ds\quad\text{ for all }(s,t)\in I\times[0,T],
\] 
in particular constant in $s$ for fixed $t$. Then Theorem \ref{th_coord} yields for $t\in[0,T]$:
\begin{align}\label{eq_coord_bem}
J_t(0,s)=\frac{L(t)}{2},\quad |\nabla s|\circ\overline{X}_0=\frac{2}{L(t)},\quad (N_{\partial\Omega}\cdot\nabla s)|_{\overline{p}^\pm(t)}=\pm\frac{2}{L(t)}.
\end{align}
\newline
\textit{Proof.} Given an arbitrary $X_0$ as above, this can be achieved as follows: we consider
\[
B:I\times[0,T]\rightarrow I:(s,t)\mapsto\frac{2}{L(t)}\int_{-1}^s|\partial_sX_0(\sigma,t)|\,d\sigma-1.
\]
Then $B$ is smooth and $\partial_sB>0$. Hence $B(.,t)$ is invertible for all $t\in[0,T]$ and the inverse mapping theorem applied to a smooth extension of $(B,\textup{pr}_t)$ on $I\times[0,T]$ yields the smoothness of the inverse in $(s,t)$ and hence $\tilde{X}_0(s,t):=X_0(B(.,t)^{-1}|_s,t)$ is the desired parametrization.\hfill$\square$
\end{enumerate}
\end{Remark}
As a starting point for the proof of Theorem \ref{th_coord} we show in the following lemma, that there are graph descriptions of $\partial\Omega$ viewed from the tangential planes at the contact points $p^\pm(t), t\in[0,T]$ in uniform neighbourhoods.
\begin{Lemma}\label{th_coord_lemma}
There is an $\eta>0$ and $w_\pm:B_\eta(0)\times[0,T]\rightarrow\R$ smooth such that
\[
(-\eta,\eta)=B_\eta(0)\ni r\mapsto p^\pm(t)+r\vec{n}(\pm1,t)+w_\pm(r,t)\vec{\tau}(\pm1,t)
\]
describes $\partial\Omega$ in the square $R^\pm_\eta(t):=p^\pm(t)+B_\eta(0)\vec{n}(\pm1,t)+B_\eta(0)\vec{\tau}(\pm1,t)$ and it holds $w_\pm(0,t)=\partial_rw_\pm(0,t)=0$ for all $t\in[0,T]$. Cf. Figure \ref{fig_coord1} for a sketch of the situation.
\begin{figure}[h]
\centering
\def\svgwidth{0.4\linewidth}
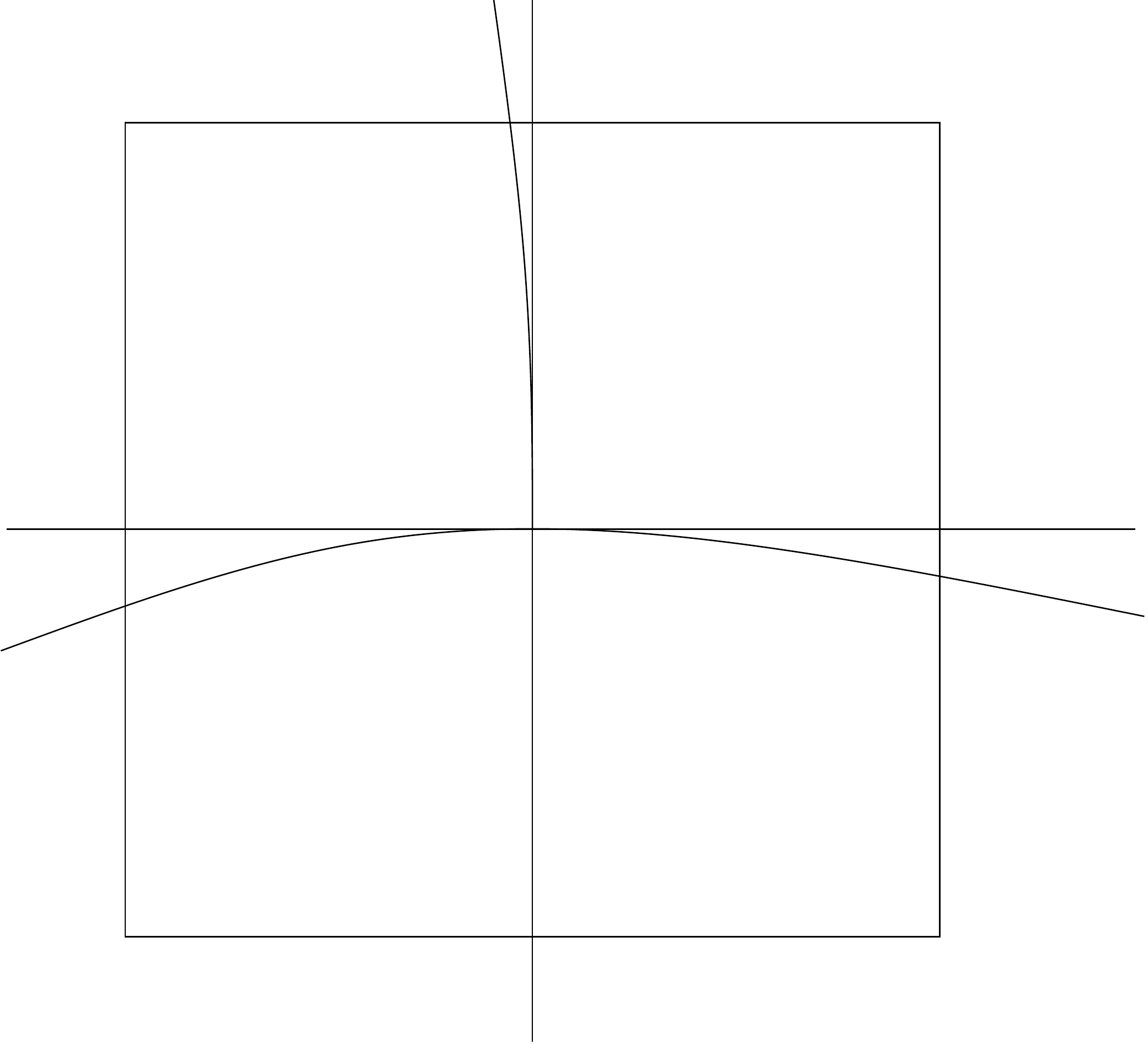
\caption{Lemma \ref{th_coord_lemma}}
\label{fig_coord1}
\end{figure}
\end{Lemma} 
\textit{Proof.} Let us fix $t_0\in[0,T]$. Then there is a graph parametrization of $\partial\Omega$ in a neighbourhood of $p^\pm(t_0)$ with a $\gamma:B_{\delta_0}(0)\subseteq\R\rightarrow\R$ in the $(\vec{n},\vec{\tau})(\pm1,t_0)$-coordinate system based at $p^\pm(t_0)$. The boundary points will stay nearby for small time variations: If $\varepsilon>0$ is small, then
\[
x^\pm(t):=(p^\pm(t)-p^\pm(t_0))\cdot\vec{n}(\pm 1,t_0)\in B_{\delta_0/2}(0)\quad\text{ for all }t\in \overline{B_{\varepsilon}(t_0)}\cap[0,T].
\]
The idea is to invert the projection of $\partial\Omega$ to the tangent space based at $p^\pm(t)$. Therefore we consider the smooth mapping $F^\pm:\overline{B_{\delta_0/2}(0)}\times(\overline{B_\varepsilon(t_0)}\cap[0,T])\rightarrow\R$ defined by
\[
F^\pm(x,t):=\left[(x+x^\pm(t))\vec{n}(\pm1,t_0)+\gamma(x+x^\pm(t))\vec{\tau}(\pm1,t_0)\right]\cdot\vec{n}(\pm1,t).
\]
Since $\partial_xF_\pm(0,t)\neq0$ for all $t\in\overline{B_\varepsilon(t_0)}\cap[0,T]$, the inverse mapping theorem applied to a smooth extension of $(F^\pm,\textup{pr}_t)$ and a compactness argument for the domain of the inverse yields that there is an $\eta>0$ such that for all $t\in\overline{B_\varepsilon(t_0)}\cap[0,T]$ there is an open neighbourhood $V_t$ of $0$ such that $F^\pm(.,t):V_t\rightarrow B_\eta(0)$ is a smooth diffeomorphism and
\[
B_\eta(0)\times(\overline{B_\varepsilon(t_0)}\cap[0,T])\rightarrow\R:(r,t)\mapsto F_\pm(.,t)^{-1}|_r
\]
is smooth. By construction
\[
w_\pm:B_\eta(0)\times(\overline{B_\varepsilon(t_0)}\cap[0,T])\rightarrow\R:(r,t)\mapsto(\textup{id},\gamma)(F_\pm(.,t)^{-1}|_r)\cdot\vec{\tau}(\pm1,t)
\]
has the claimed properties for $t\in \overline{B_\varepsilon(t_0)}\cap[0,T]$ after possibly shrinking $\eta$. Finally, compactness of $[0,T]$ implies the Lemma.\hfill$\square$\\
\newline
\textit{Proof of Theorem \ref{th_coord}.} The idea for the definition of $X$ is to extend the mapping in Lemma \ref{th_coord_lemma} such that it coincides with the usual tubular neighbourhood coordinate system outside a neighbourhood of the boundary and such that all the claimed properties are satisfied.\\
Therefore we first consider the construction of the latter coordinate system in Hildebrandt \cite{Hildebrandt_Analysis2}, Chapter 4.6. Observe that $X_0$ can be extended smoothly to an $\tilde{X}_0$ with analogous properties defined on an open neighbourhood $\tilde{I}\times(-\varepsilon,T+\varepsilon)$ of $I\times[0,T]$ and hence $\Gamma$ is compactly contained in a smooth manifold $\tilde{\Gamma}=\bigcup_{t\in[0,T]}\tilde{\Gamma}_t\times\{t\}$ that is parametrized via $\tilde{X}_0$. Similar ideas as in the proof of Lemma \ref{th_coord_lemma} above yield local graph parametrizations for $\tilde{\Gamma}_t$ as in Lemma 1 in Hildebrandt \cite{Hildebrandt_Analysis2}, Chapter 4.6, viewed from the tangent planes in squares of uniform width around every point in $\tilde{X}_0(I_0,t)$ for all $t\in[0,T]$ and an $I\subseteq I_0\subseteq \tilde{I}$ open. Therefore the construction in \cite{Hildebrandt_Analysis2}, Satz 1 in Chapter 4.6, yields that for all $\delta\in(0,\delta_0]$, where $\delta_0>0$ is small but independent of $t$, it holds that $(-\delta,\delta)\times I_0\ni (r,s)\mapsto \tilde{X}_0(s,t)+r\vec{n}(s,t)\in\R^2$ is a diffeomorphism onto its image (denote it by $S_\delta(t)$) and $S_\delta(t)\cap\overline{\Omega}=B_\delta(\tilde{\Gamma}_t)\cap\overline{\Omega}$ for all $t\in[0,T]$. Now we fix $\delta_0>0$ small.\\
We choose $\eta>0$ small such that $R^\pm_\eta(t)$ is contained in $S_{\delta_0}(t)$, the assertions of Lemma \ref{th_coord_lemma} are fulfilled and $\tilde{\Gamma}_t\cap R^\pm_\eta(t)$ has a graph parametrization in the square $R^\pm_\eta(t)$ viewed from the tangent plane for all $t\in[0,T]$. See figure \ref{fig_coord2} for a sketch of the situation.
\begin{figure}[h]
\centering
\def\svgwidth{0.5\linewidth}
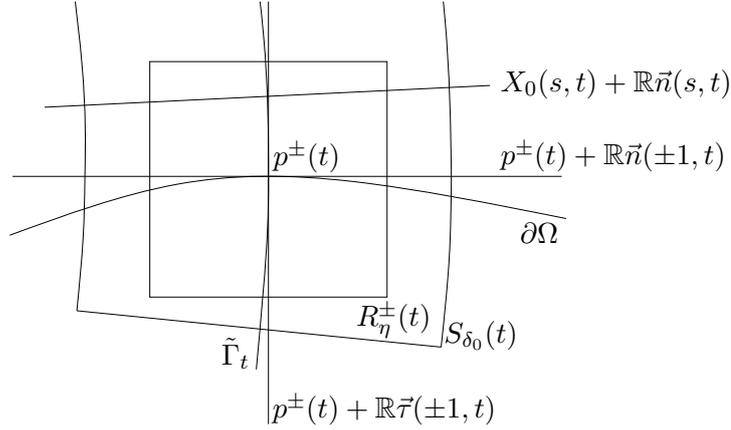
\caption{Construction of curvilinear coordinates}
\label{fig_coord2}
\end{figure}
Moreover, by shrinking $\eta$ we can assume\footnote{ Cf. Lemma 1 and Lemma 2 in Hildebrandt \cite{Hildebrandt_Analysis2}, Chapter 4.6.} that the intersections of $\partial\Omega$ and $\tilde{\Gamma}_t$ with the square lie in $\frac{\eta}{2}$-strips, respectively, i.e.
\[
\partial\Omega\cap R^\pm_{\eta}(t)\subseteq R^\pm_{\eta/4}(t)+\R\vec{n}(\pm1,t)
\quad\text{ and }\quad 
\tilde{\Gamma}_t\cap R^\pm_{\eta}(t)\subseteq R^\pm_{\eta/4}(t)+\R\vec{\tau}(\pm1,t)
\] 
for all $t\in[0,T]$ and the angles between the tangent vectors of $\partial\Omega\cap R^\pm_{\eta}(t)$ and $\tilde{\Gamma}_t\cap R^\pm_{\eta}(t)$, respectively, can be made smaller than an arbitrary fixed $\alpha>0$.\\ 
Now we define $X$. Because of uniform continuity we can choose an $\varepsilon>0$ such that for all $s\in I$ with $|s\mp 1|\leq\varepsilon$ and $t\in[0,T]$ it holds $X_0(s,t)\in R^\pm_{\eta/2}(t)$. Let $\chi:I\rightarrow[0,1]$ be a smooth cutoff-function with $\chi=1$ for $|s\mp 1|\leq\frac{\varepsilon}{2}$ and $\chi=0$ for $|s|\leq 1-\varepsilon$. Then we define
$\vec{T}:=\chi\vec{\tau}$ and for $\delta\in(0,\eta]$
\[
X(r,s,t):=X_0(s,t)+r\vec{n}(s,t)+w(r,s,t)\vec{T}(s,t)\quad\text{ for }(r,s,t)\in[-\delta,\delta]\times I\times[0,T],
\] 
where $w(r,s,t):=w_{\textup{sign}(s)}(r,t)$. In the following we show that the properties in the theorem are satisfied if $\delta>0$ is small and $\alpha>0$ above was chosen properly. $X$ is smooth and we compute
\begin{align}\begin{split}\label{eq_coord_X_abl}
\partial_rX(r,s,t)&=\vec{n}(s,t)+\partial_rw(r,s,t)\vec{T}(s,t),\\
\partial_sX(r,s,t)&=\partial_sX_0(s,t)+r\partial_s\vec{n}(s,t)+w(r,s,t)\partial_s\vec{T}(s,t)
\end{split}\end{align}
for $(r,s,t)\in[-\delta,\delta]\times I\times[0,T]$. If $|s|\leq 1-\varepsilon$, then $X$ coincides with the usual tubular neighbourhood coordinate system, $\partial_rX(r,s,t)=\vec{n}(s,t)$ and $\partial_sX(r,s,t)=\partial_sX_0(s,t)$. By uniform continuity and because of $w=\partial_rw=0$ for $r=0$, the perturbation terms for $|s\mp 1|\leq\varepsilon$ are estimated by an arbitrary small constant $c_0>0$ if $\delta>0$ is small. If $c_0$ is small enough, then $\measuredangle(\partial_rX(r,s,t),\vec{n}(\pm1,t))\leq 2\alpha$ and $\measuredangle(\partial_sX(r,s,t),\vec{\tau}(\pm1,t))\leq 2\alpha$ for all $s\in I$ such that $X_0(s,t)\in R^\pm_\eta(t)$. Hence, if $\alpha>0$ was chosen sufficiently small, $\overline{X}$ is injective by a contradiction argument using the Fundamental Theorem of Calculus. Moreover, $\overline{X}$ maps to $\overline{\Omega}\times[0,T]$, $D\overline{X}(r,s,t)$ is invertible for all $(r,s,t)\in[-\delta,\delta]\times I\times[0,T]$ and thus $\overline{X}$ is a homeomorphism onto its image. The inverse mapping theorem applied to a smooth extension of $X$ implies smoothness of the inverse. Next, we prove that the image is indeed a neighbourhood of $\Gamma$ in $\overline{\Omega}\times[0,T]$. Therefore we observe that $w_\pm(r,t)=\Oc(r^2)$ because of a Taylor expansion and hence 
\[
\textup{sdist}_{\tilde{\Gamma}_t}(X(r,s,t))=r+\Oc(r^2),
\] 
where $\textup{sdist}_{\tilde{\Gamma}_t}:S_{\delta_0}(t)\cap\overline{\Omega}\rightarrow\R$ is the signed distance function to $\tilde{\Gamma}_t$ with respect to $\vec{n}(.,t)$, cf. \cite{Hildebrandt_Analysis2}, Chapter 4.6. Thus for all $\delta>0$ small we have 
\[
B_{\delta/2}(\tilde{\Gamma}_t)\cap\overline{\Omega}\subseteq \overline{\Gamma_t(\delta)}=X([-\delta,\delta]\times I,t)\subseteq B_{3\delta/2}(\tilde{\Gamma}_t)\cap\overline{\Omega}.
\]
Therefore $\Gamma(\delta)$ is a neighbourhood of $\Gamma$ in $\overline{\Omega}\times[0,T]$ for all $\delta>0$ small and one can prove that it is open in $\overline{\Omega}\times[0,T]$.\\
It remains to prove the explicit identities in Theorem \ref{th_coord}. We have
\[
((D_{(r,s)}X)^\top D_{(r,s)}X)^{-1}=\left.
\begin{pmatrix}
|\nabla r|^2          & \nabla r\cdot\nabla s\\
\nabla r\cdot\nabla s & |\nabla s|^2
\end{pmatrix}\right|_{\overline{X}},\quad
\begin{pmatrix}
a & b\\
b & c
\end{pmatrix}
:=(D_{(r,s)}X)^\top D_{(r,s)}X.
\]
Using \eqref{eq_coord_X_abl} and that $\vec{T},\partial_sX_0,\partial_s\vec{n}$ are tangential, we obtain $a=1+(\partial_rw)^2|\vec{T}|^2$,
\begin{align*}
b&=w\partial_s\vec{T}\cdot\vec{n}+\partial_rw\partial_sX_0\cdot\vec{T}+r\partial_rw\partial_s\vec{n}\cdot\vec{T}+w\partial_rw\partial_s\vec{T}\cdot\vec{T},\\ c&=|\partial_sX_0|^2+r^2|\partial_s\vec{n}|^2+w^2|\partial_s\vec{T}|^2+2(r\partial_sX_0\cdot\partial_s\vec{n}+w\partial_sX_0\cdot\partial_s\vec{T}+rw\partial_s\vec{n}\cdot\partial_s\vec{T}). 
\end{align*}
The inverse can be computed explicitly. Since $w=\partial_rw=0$ for $r=0$, it follows that $a|_{(0,s,t)}=1, b|_{(0,s,t)}=0$ and $c|_{(0,s,t)}=|\partial_sX_0(s,t)|^2$ for $(s,t)\in I\times[0,T]$. Therefore
\[
|\nabla r|^2|_{\Gamma}=1,\quad \nabla r\cdot\nabla s|_{\Gamma}=0\quad\text{ and }\quad |\nabla s|\circ\overline{X}_0=1/|\partial_sX_0|. 
\]
Moreover, with $d:=ac-b^2$ we obtain $\partial_r(|\nabla r|^2\circ\overline{X})=(d\partial_rc-c\partial_rd)/d^2$,
\[
\partial_r((\nabla r\cdot\nabla s)\circ \overline{X})=-\frac{d\partial_rb-b\partial_rd}{d^2}\quad\text{ and }\quad 
\partial_r(|\nabla s|^2\circ \overline{X})=\frac{d\partial_ra-a\partial_rd}{d^2}.
\]
We have $\partial_ra|_{(0,s,t)}=0$, $\partial_rb|_{(0,s,t)}=\partial_r^2w|_{(0,s,t)}\partial_sX_0\cdot\vec{T}|_{(s,t)}$ and $\partial_rc|_{(0,s,t)}=2\partial_sX_0\cdot\partial_s\vec{n}|_{(s,t)}$ for $(s,t)\in I\times[0,T]$. Hence $\partial_r^kd|_{(0,s,t)}=\partial_r^kc|_{(0,s,t)}, k=0,1$ and $\partial_r(|\nabla r|^2\circ \overline{X})|_{(0,s,t)}=0$,
\[
\partial_r((\nabla r\cdot\nabla s)\circ \overline{X})|_{(0,s,t)}=\frac{-\partial^2_rw|_{(0,s,t)}\partial_sX_0\cdot\vec{T}|_{(s,t)}}{|\partial_sX_0|^2|_{(s,t)}}
\]
and $\partial_r(|\nabla s|^2\circ \overline{X})|_{(0,s,t)}=-2\partial_sX_0\cdot\partial_s\vec{n}|_{(s,t)}/|\partial_sX_0|^4|_{(s,t)}$ for all $(s,t)\in I\times[0,T]$.\\
Finally, we show that for the coordinate system constructed above the additional properties are satisfied. Because of \eqref{eq_coord_X_abl} we have for all $(s,t)\in I\times[0,T]$
\[
D_{(r,s)}X|_{(0,s,t)}=\begin{pmatrix}
\vec{n}, \partial_sX_0
\end{pmatrix}|_{(s,t)},\quad
(D_{(r,s)}X|_{(0,s,t)})^{-1}=
\begin{pmatrix}
\vec{n}^\top\\
\partial_sX_0^\top/|\partial_sX_0|^2
\end{pmatrix}|_{(s,t)}.
\]
This shows $\nabla r\circ\overline{X}_0=\vec{n}$ and $\nabla s\circ\overline{X}_0=\partial_sX_0/|\partial_sX_0|^2$. Furthermore, the chain rule applied to $r=r(X(r,s,t),t)$ yields
\[
V(s,t)=\partial_tX|_{(0,s,t)}\cdot\vec{n}(s,t)=-\partial_tr|_{(X(0,s,t),t)}\quad\text{ for }(s,t)\in I\times[0,T].
\]
Well-known formulas for the mean curvature, cf. e.g. Depner \cite{Depner}, Chapter 2.1, imply
\[
H_{\Gamma_t}(s)=-\diverg\left(\frac{\nabla r}{|\nabla r|}\right) |_{(X_0(s,t),t)}=-\left(\frac{\Delta r}{|\nabla r|}-\frac{1}{|\nabla r|^3}\nabla r^\top D^2r\nabla r\right)|_{(X_0(s,t),t)}
\]
for $(s,t)\in I\times[0,T]$. The second term vanishes because of $D^2r\nabla r=\frac{1}{2}\nabla(|\nabla r|^2)$ and
\[
\partial_r(|\nabla r|^2\circ \overline{X})|_{(0,s,t)}=\partial_s(|\nabla r|^2\circ \overline{X})|_{(0,s,t)}=0\quad\text{ for }(s,t)\in I\times[0,T].
\] 
With $|\nabla r|^2|_\Gamma=1$ we get $H(s,t)=-\Delta r|_{(X_0(s,t),t)}$ for $(s,t)\in I\times[0,T]$.\hfill$\square$\\
\subsection{Results from ODE-Theory}\label{sec_ODE}
In this section we review results needed from ODE-theory.
Let 
\[
f\in C^\infty(\R),\quad f'(\pm1)=0,\quad f''(\pm1)>0,\quad \int_{-1}^{u}f'=\int_1^u f'>0\quad \text{ for all }u\in(-1,1).
\] 
Schaubeck \cite{Schaubeck}, Lemma 2.6.1 yields
\begin{Lemma}\label{th_theta_0}
The system
\[
-w''+f'(w)=0,\quad w(0)=0,\quad \lim_{z\rightarrow\pm\infty}w(z)=\pm1
\]
has exactly one solution $\theta_0\in C^2(\R)$. Moreover, $\theta_0$ is smooth, $\theta_0'>0$ and it holds $D_z^k(\theta_0\mp 1)(z)=\Oc(e^{-\alpha|z|})$ for $z\rightarrow\pm\infty$ and all $k\in\N_0, \alpha\in(0,\min\{\sqrt{f''(\pm1)}\})$.
\end{Lemma}
We call $\theta_0$ the \textit{optimal profile}. A rescaled version will be the typical profile of the solutions across the interface. If $f$ is even, then $\theta_0$ is even, $\theta_0'$ is odd and $\theta_0''$ even etc. The linearization at $\theta_0$, i.e. $\Lc:=-\frac{d^2}{dz^2}+f''(\theta_0)$ will appear in the asymptotic expansion, too. In the next lemma we recall some properties of $\Lc$ viewed as an unbounded operator on $L^2(\R)$.
\begin{Lemma}\label{th_ODE_lin1}
The operator $\Lc:H^2(\R)\subseteq L^2(\R)\rightarrow L^2(\R)$ is selfadjoint, non-negative and the set $\sigma(\Lc)\cap(-\infty,\min\{f''(\pm1)\})$ lies in the discrete spectrum, i.e. consists only of isolated eigenvalues with finite algebraic multiplicity. The lowest eigenvalue is $0$ and $\ker \Lc$ is spanned by $\theta_0'$. Moreover, with $(\ker \Lc)^\perp:=\{w\in L^2(\R):(w,\theta_0')_{L^2(\R)}=0\}$ it holds
\[
0<\nu_2
:=\inf_{w\in H^2(\R)\cap(\ker\Lc)^\perp}(\Lc w,w)_{L^2(\R)}
=\inf_{w\in H^1(\R)\cap(\ker\Lc)^\perp}\int_{\R}|w'|^2+f''(\theta_0)|w|^2\,dz.
\]
\end{Lemma}
\textit{Proof.} Because $f''(\theta_0)$ is bounded from below, the Lax-Milgram Theorem and regularity theory imply that $\sigma(\Lc)$ is bounded from below, in particular $\rho(\Lc)\cap\R\neq\emptyset$. Since $\Lc$ is densely defined and symmetric, it follows that $\Lc$ is selfadjoint. Pankov \cite{Pankov}, Theorem 8.20 implies that $\sigma(\Lc)\cap(-\infty,\min\{f''(\pm1)\})$ is contained in the discrete spectrum. In particular $\lambda_0:=\inf\sigma(\Lc)$ is an isolated eigenvalue with finite algebraic multiplicity. Reed, Simon \cite{ReedSimonIV}, Theorem XIII, 48 yields that $\lambda_0$ is simple and corresponding eigenfunctions have a sign. Since $\Lc\theta_0'=0, \theta_0'>0$ by Lemma \ref{th_theta_0} and because eigenfunctions to distinct eigenvalues are orthogonal, by contradiction it follows that $\lambda_0=0$ and the eigenspace $\ker \Lc$ is spanned by $\theta_0'$. This yields that $\Lc$ is non-negative because of Hislop, Sigal \cite{HislopSigal}, Proposition 5.12 (or Werner \cite{Werner}, Aufgabe VII.5.24 on p.383). One can directly show that $\tilde{\Lc}:=\Lc|_{(\ker \Lc)^\perp}:H^2(\R)\cap(\ker \Lc)^\perp\rightarrow(\ker \Lc)^\perp$ is selfadjoint and $\sigma(\Lc)=\sigma(\tilde{\Lc})\cup\{0\}$. Since $0$ is not an eigenvalue of $\tilde{\Lc}$ and isolated points of $\sigma(\tilde{\Lc})$ are eigenvalues because of Hislop, Sigal \cite{HislopSigal}, Proposition 6.4, we obtain that $\sigma(\tilde{\Lc})=\sigma(\Lc)\textbackslash\{0\}$ and hence $\nu_2>0$ with \cite{HislopSigal}, Proposition 5.12. The last identity for $\nu_2$ follows with a density argument and integration by parts.\hfill$\square$\\
\newline
The following lemma is concerned with the solvability of the equation, which is obtained by linearization at $\theta_0'$. The result follows from the proof of Schaubeck \cite{Schaubeck}, Lemma 2.6.2.
\begin{Lemma}\label{th_ODE_lin2}
Let $A\in C^0(\R)$ with $\lim_{z\rightarrow\pm\infty}A(z)=A^\pm\in\R$. Then 
\begin{align}\label{eq_ODE}
-w''+f''(\theta_0)w=A\text{ in }\R,\quad w(0)=0
\end{align}
has a solution $w\in C^2(\R)\cap C_b^0(\R)$ if and only if $\int_{\R}A\theta_0'\,dz=0$. Moreover, if this is the case, then $w$ is unique and satisfies
\[
D_z^l\left[w-\frac{A^\pm}{f''(\pm1)}\right]=\Oc(e^{-\alpha|z|})\quad\text{ for }z\rightarrow\pm\infty,\quad l=0,1,2.
\]
Furthermore, if $U\subseteq\R^d$ and $A:\R\times U\rightarrow\R, A^\pm:U\rightarrow\R$ are smooth (i.e. locally smooth extendible) and it holds uniformly in $U$:
\[
D_x^m D_z^l\left[A(z,.)-A^\pm\right]=\Oc(e^{-\alpha|z|})\quad\text{ for }z\rightarrow\pm\infty,\quad m=0,...,M, l=0,...,L,
\] 
for some $M,N\in\N_0$, then $w:\R\times U\rightarrow\R$, where $w(.,x)$ is the solution of \eqref{eq_ODE} for $A(.,x)$ for all $x\in U$, is also smooth and it holds uniformly in $U$:
\[
D_x^m D_z^l\left[w(z,.)-\frac{A^\pm}{f''(\pm1)}\right]=\Oc(e^{-\alpha|z|})\quad\text{ for }z\rightarrow\pm\infty,\quad m=0,...,M, l=0,...,L+2.
\] 
\end{Lemma}
\subsection{Model Problem on $\R^2_+$}\label{sec_hp}
The arguments in this section are motivated by the brief sketch of the proof for Lemma 2.1 in Alikakos, Chen and Fusco \cite{ACF}. For the boundary layer expansion we need to solve the following model problem on $\R^2_+$: We consider
\begin{alignat}{2}\label{eq_hp_model1}
(-\Delta+f''(\theta_0(R)))u(R,H)&=G(R,H)& \quad &\text{ for }(R,H)\in\R^2_+,\\
-\partial_H u|_{H=0}(R)&=g(R)&\quad &\text{ for }R\in\R.
\label{eq_hp_model2}
\end{alignat}  
In Section \ref{sec_hp_weak_sol_reg} we show existence and uniqueness of weak solutions. The Lax-Milgram Theorem cannot be applied directly since coercivity fails. Therefore we split $L^2(\R^2_+)$ and $L^2(\R)$ into $\theta_0'$-$L^2(\R)$-orthogonal parts and use the Lax-Milgram Theorem for the orthogonal parts. For the parallel parts it turns out that for suitable $G,g$ satisfying the compatibility condition
\[
\int_{\R^2_+}G(R,H)\theta_0'(R)\,d(R,H)+\int_{\R}g(R)\theta_0'(R)\,dR=0
\] 
there is an explicit solution formula. For more regular data we apply regularity theory. Moreover, in Subsection \ref{sec_hp_exp_decay} we prove pointwise exponential decay of the solution $u$ and its derivatives for appropriate data. Here the basic idea is to consider $H\mapsto \|u(.,H)\|_{L^2(\R)}$ and $R\mapsto\|u(R,.)\|_{L^2(\R_+)}$ and to derive ordinary differential inequalities where possible. Then, if $G,g$ are suitable, we can show by contradiction estimates of type
\[
\|u(.,H)\|_{L^2(\R)}\leq C_u e^{-\nu H}\quad\text{ f.a.a. }H\in\R_+,\quad \|u(R,.)\|_{L^2(\R_+)}\leq C_u e^{-\nu|R|}\quad\text{ f.a.a. }R\in\R,
\]   
where $\nu\in(0,\sqrt{\nu_2})$ and $\nu_2$ is as in Lemma \ref{th_ODE_lin1}. Then by differentiating and rearranging the equations and by interpolation we get similar estimates for the derivatives if the data are appropriate. With embeddings we also obtain pointwise estimates. Finally, in Section \ref{sec_hp_time} we consider the case when $G,g$ also depend on $t\in[0,T]$.
\subsubsection{Weak Solutions and Regularity}\label{sec_hp_weak_sol_reg}
Let us start with the definition of weak solution:
\begin{Definition}\label{th_hp_weak_sol_def}\upshape
Let $G\in L^2(\R^2_+)$ and $g\in L^2(\R)$. Then $u\in H^1(\R^2_+)$ is called \textit{weak solution} of \eqref{eq_hp_model1}-\eqref{eq_hp_model2} if for all $\varphi\in H^1(\R^2_+)$ it holds that
\[
a(u,\varphi):=\int_{\R^2_+}\nabla u\cdot\nabla\varphi+f''(\theta_0(R))u\varphi\,d(R,H)=\int_{\R^2_+}G\varphi\,d(R,H)+\int_{\R}g(R)\varphi|_{H=0}(R)\,dR. 
\]
\end{Definition}
For the following we need to know how Lebesgue and Sobolev spaces behave on product sets. This is the content of
\begin{Lemma}\label{th_prod_sets_spaces} Let $\Omega_1\subseteq\R^n, \Omega_2\subseteq \R^m$ for $n,m\in\N$ be measurable.
\begin{enumerate}
\item  Let $f\in L^p(\Omega_1\times \Omega_2)$ for some $p\in[1,\infty)$. Then $f(x_1,.)\in L^p(\Omega_2)$ f.a.a. $x_1\in \Omega_1$, $Tf:\Omega_1\rightarrow L^p(\Omega_2):x_1\mapsto f(x_1,.)$ is strongly measurable and belongs to $L^p(\Omega_1, L^p(\Omega_2))$. Furthermore, $T:L^p(\Omega_1\times \Omega_2)\rightarrow L^p(\Omega_1,L^p(\Omega_2))$ is an isometric isomorphism.
\item If $\Omega_1,\Omega_2$ are open, then $W^1_p(\Omega_1\times\Omega_2)\cong L^p(\Omega_1,W^1_p(\Omega_2))\cap W^1_p(\Omega_1,L^p(\Omega_2))$ for $p\in(1,\infty)$ according to the above isomorphism. Moreover, the derivatives in the respective spaces are compatible and analogous assertions hold for higher orders. 
\item For $f\in W^k_p(\Omega_1)$ and $g\in W^k_p(\Omega_2)$, $p\in[1,\infty]$, $k\in\N_0$ it holds that
\[
(f\otimes g)(x_1,x_2):=f(x_1)g(x_2)
\]
is well-defined f.a.a. $(x_1,x_2)\in\Omega_1\times\Omega_2$, $f\otimes g\in W^k_p(\Omega_1\times\Omega_2)$, the derivatives are natural and $\|f\otimes g\|_{W^k_p(\Omega_1\times\Omega_2)}\leq C_{k,p}\|f\|_{W^k_p(\Omega_1)}\|g\|_{W^k_p(\Omega_2)}$. In case $\Omega_1=\R, \Omega_2=\R_+$ and $k=1, p\in[1,\infty)$ the trace is $\tr_{\partial\R^2_+}(f\otimes g)=g(0)f$.
\end{enumerate}
\end{Lemma}
\textit{Proof.} The first two assertions follow with ideas on p.40 in R\r{u}\v{z}i\v{c}ka \cite{Ruzicka} and p.199 in Schweizer \cite{Schweizer}. The third claim can be proven directly with the definitions.\hfill$\square$\\
\newline
Regarding weak solutions we have the following theorem:
\begin{Theorem}\label{th_hp_weak_sol} Let $G\in L^2(\R^2_+)$ and $g\in L^2(\R)$. Then it holds:
\begin{enumerate}
\item $a:H^1(\R^2_+)\times H^1(\R^2_+)\rightarrow\R$ is not coercive.
\item If $G(.,H),g\perp\theta_0'$ f.a.a. $H>0$ in $L^2(\R)$, then there is a weak solution $u$ with $u(.,H)\perp\theta_0'$ f.a.a. $H>0$ and it holds $\|u\|_{H^1(\R^2_+)}\leq C(\|G\|_{L^2(\R^2_+)}+\|g\|_{L^2(\R)})$.
\item Weak solutions are unique.
\item If $G\theta_0'\in L^1(\R^2_+)$ and $u$ is a weak solution with $\partial_Hu\,\theta_0'\in L^1(\R^2_+)$, then the following compatibility condition holds:
\begin{align}\label{eq_hp_intbed}
\int_{\R^2_+}G(R,H)\theta_0'(R)\,d(R,H)+\int_{\R}g(R)\theta_0'(R)\,dR=0.
\end{align}
\item If $G\theta_0'\in L^1(\R^2_+)$, then $\tilde{G}(H):=(G(.,H),\theta_0')_{L^2(\R)}$ is well defined f.a.a. $H>0$ and $\tilde{G}\in L^1(\R_+)\cap L^2(\R_+)$. Moreover, with $\beta:=\|\theta_0'\|_{L^2(\R)}^{-1}$ we have the decomposition
\begin{align}\label{eq_hp_orth_decomp}
G=\tilde{G}(H)\beta^2\theta_0'(R)+G^\perp(R,H),\quad g=(g,\theta_0')_{L^2(\R)}\beta^2\theta_0'(R)+g^\perp(R)
\end{align}
for some $G^\perp\in L^2(\R^2_+), g^\perp\in L^2(\R)$ with $G^\perp(.,H),g\perp \theta_0'$ in $L^2(\R)$ f.a.a. $H>0$.
\item If $\|G(.,H)\|_{L^2(\R)}\leq C e^{-\nu H}$ f.a.a. $H>0$ and a constant $\nu>0$, then $G\theta_0'\in L^1(\R^2_+)$. Let $\tilde{G}$ be defined as in 4. and the compatibility condition \eqref{eq_hp_intbed} hold. Then 
\begin{align}\label{eq_hp_sol_formula}
\tilde{u}(R,H):=-\beta^2\int_H^\infty\int_{\tilde{H}}^{\infty}\tilde{G}(\hat{H})\,d\hat{H}\,d\tilde{H}\,\theta_0'(R)
\end{align}
is well defined f.a.a. $(R,H)\in\R^2_+$, $\tilde{u}\in W^2_1(\R^2_+)\cap H^2(\R^2_+)$ and $\tilde{u}$ is a weak solution of \eqref{eq_hp_model1}-\eqref{eq_hp_model2} for $G-G^\perp,g-g^\perp$ in \eqref{eq_hp_orth_decomp} instead of $G,g$.
\end{enumerate}
\end{Theorem}
In Theorem \ref{th_hp_weak_sol}, 6. weaker conditions on $G$ should be enough, but later we are just interested in smooth and exponentially decaying data anyway. Altogether we obtain
\begin{Corollary}\label{th_hp_weak_sol_2}
\begin{enumerate}
\item Let $g\in L^2(\R)$, $G\in L^2(\R^2_+)$ with $\|G(.,H)\|_{L^2(\R)}\leq Ce^{-\nu H}$ f.a.a. $H>0$ and let \eqref{eq_hp_intbed} hold. Then there is a unique weak solution of \eqref{eq_hp_model1}-\eqref{eq_hp_model2}. 
\item Let $m\in\N_0$ and $u\in H^1(\R^2_+)$ be a weak solution of \eqref{eq_hp_model1}-\eqref{eq_hp_model2} for $G\in H^m(\R^2_+)$, $g\in H^{m+\frac{1}{2}}(\R)$, then $u\in H^{m+2}(\R^2_+)\hookrightarrow C^{m,\gamma}(\overline{\R^2_+})$ for all $\gamma\in(0,1)$ and it holds 
\[
\|u\|_{H^{m+2}(\R^2_+)}\leq C_m(\|G\|_{H^m(\R^2_+)}+\|g\|_{H^{m+\frac{1}{2}}(\R)}+\|u\|_{H^1(\R^2_+)}).
\]
\end{enumerate}
\end{Corollary}
\textit{Proof.} The first part directly follows from Theorem \ref{th_hp_weak_sol}. For the second assertion, we apply Triebel \cite{Triebel_Fct_SpacesI}, Theorem 2.7.2 to obtain a $\overline{g}\in H^{m+2}(\R^2_+)$ with $(-\partial_H\overline{g})|_{H=0}=g$ and $\|\overline{g}\|_{H^{m+2}(\R^2_+)}\leq C\|g\|_{H^{m+\frac{1}{2}}(\R)}$. Subtracting $\overline{g}$ from $u$ and using standard regularity theory, we get $u\in H^{m+2}(\R^2_+)$ and the estimate.\hfill$\square$\\
\newline
\textit{Proof of Theorem \ref{th_hp_weak_sol}. Ad 1.} We consider smooth cut-off functions $\chi_n:\R_+\rightarrow[0,1]$ for $n\in\N$ with $|\chi_n'|\leq C$, $\chi_n=1$ for $H\leq n$ and $\chi_n=0$ for $H\geq n+1$. If $a$ would be coercive on $H^1(\R^2_+)$, then for some $c>0$
\begin{align*}
a(\theta_0'(R)\chi_n(H),\theta_0'(R)\chi_n(H))&=\int_{\R^2_+}(\theta_0'\chi_n')^2+\left[(\theta_0'')^2+f''(\theta_0)\theta_0'^2\right]\chi_n^2\,d(R,H)=\\
&=\int_\R(\theta_0')^2\,dR\int_{\R_+}(\chi_n')^2\,dH\geq c\int_\R(\theta_0')^2\,dR\int_{\R_+}\chi_n^2\,dH,
\end{align*}
where we used integration by parts with respect to $R\in\R$ in the second term. This is a contradiction for $n\rightarrow\infty$.\hfill$\square_{1.}$\\
\newline
\textit{Ad 2.} Let $G\in L^2(\R^2_+), g\in L^2(\R)$ with $ G(.,H),g\perp\theta_0'$ in $L^2(\R)$ f.a.a. $H>0$. In order to show existence of a weak solution we apply the Lax-Milgram Theorem to $V:=\{u\in H^1(\R^2_+):u(.,H)\perp\theta_0'\textup{ in }L^2(\R)\textup{ f.a.a. }H>0\}$, the bilinear form
\[
a:V\times V\rightarrow\R:(u,v)\mapsto\int_{\R^2_+}\nabla u\cdot\nabla v+f''(\theta_0(R))uv\,d(R,H)
\]
and $x'\in V'$ defined by $x'(v):=\int_{\R^2_+}Gv+\int_{\R}g(R)v|_{H=0}(R)\,dR$ for all $v\in V$. First of all, $V$ is a Hilbert space as closed subspace of $H^1(\R^2_+)$. Here closedness follows from Lemma \ref{th_prod_sets_spaces}, 1. and linearity of $(.,\theta_0')_{L^2(\R)}:L^2(\R)\rightarrow\R$. Boundedness of $a$ can be shown directly and coercivity on $V$ follows from
\begin{align*}
a(v,v)&=\|\partial_Hv\|_{L^2(\R^2_+)}^2+\int_0^\infty\int_{\R}(\partial_Rv)^2+f''(\theta_0(R))v^2\,dR\,dH\geq\\
&\geq\|\partial_Hv\|_{L^2(\R^2_+)}^2+\delta(\|\partial_Rv\|_{L^2(\R^2_+)}^2-\max|f''(\theta_0)|\|v\|_{L^2(\R^2_+)}^2)+(1-\delta)\nu_2\|v\|_{L^2(\R^2_+)}^2\geq\\
&\geq c\|v\|_{H^1(\R^2_+)}^2
\end{align*}
for all $v\in V$ and a $c>0$ if $\delta>0$ is sufficiently small, where we used Fubini's Theorem and Lemma \ref{th_ODE_lin1}. Therefore the Lax-Milgram Theorem implies that there is a unique $u\in V$ such that $a(u,\varphi)=x'(\varphi)$ for all $\varphi\in V$ and that the estimate holds. Hence $u$ satisfies the definition of weak solution for all $\varphi\in V$. For $\varphi\in H^1(\R^2_+)$ let $\tilde{\varphi}(H):=(\varphi(.,H),\theta_0')_{L^2(\R)}$ f.a.a. $H>0$. By Lemma \ref{th_prod_sets_spaces} we have that $\tilde{\varphi}\in H^1(\R_+)$ and
\[
\varphi=\tilde{\varphi}(H)\beta^2\theta_0'(R)+\varphi^\perp\quad\textup{ with }\varphi^\perp\in H^1(\R^2_+)\textup{ such that }\varphi^\perp(.,H)\perp\theta_0'\textup{ f.a.a. }H>0.
\]
Since the definition of weak solution is linear in $\varphi$, we only have to verify it for the parallel part, i.e. we need to show
\begin{align*}
\int_{\R^2_+}\partial_Hu\,\partial_H\tilde{\varphi}(H)\theta_0'(R)+\partial_Ru\,\tilde{\varphi}(H)\theta_0''(R)+f''(\theta_0(R))u\,\tilde{\varphi}(H)\theta_0'(R)\,d(R,H)=\\
=\int_{\R^2_+}G\tilde{\varphi}(H)\theta_0'(R)\,d(R,H)+\int_{\R}g(R)\tilde{\varphi}(0)\theta_0'(R)\,dR.
\end{align*}
The right hand side is zero because of the orthogonality condition for $G,g$. The second and the last term on the left hand side cancel since we can apply integration by parts in $R$ for the second part and use $\theta_0'''=f''(\theta_0)\theta_0'$. Moreover, Lemma \ref{th_prod_sets_spaces} and linearity of $(.,\theta_0')_{L^2(\R)}:L^2(\R)\rightarrow\R$ imply
\[
0=\frac{d}{dH}(u(.,H),\theta_0')_{L^2(\R)}=(\partial_Hu(.,H),\theta_0')_{L^2(\R)}\quad\text{ f.a.a. }H>0.
\]
By Fubini's Theorem the first term vanishes. Hence $u$ is a weak solution.\hfill$\square_{2.}$\\
\newline
\textit{Ad 3.} Due to linearity it is enough to prove uniqueness for weak solutions $u$ to the data $G,g=0$. Given such a $u$, let us insert $\varphi=u$ in the Definition \ref{th_hp_weak_sol_def} of weak solution. This implies
\[
\|\partial_Hu\|_{L^2(\R^2_+)}+\int_{\R^2_+}(\partial_Ru)^2+f''(\theta_0(R))u^2\,dR\,dH=0.
\]
Because of Fubini's Theorem, Lemma \ref{th_ODE_lin1} and Lemma \ref{th_prod_sets_spaces} the second integral is non-negative. This yields $\partial_Hu=0$ and from Lemma \ref{th_prod_sets_spaces} we obtain that $u(.,H)\in L^2(\R)$ is constant in $H>0$. Thus $u=0$, otherwise we get a contradiction to $u\in L^2(\R_+,L^2(\R))$.\hfill$\square_{3.}$\\
\newline
\textit{Ad 4.} Let $G\theta_0'\in L^1(\R^2_+)$ and $u$ be a weak solution such that $\partial_Hu\,\theta_0'\in L^1(\R^2_+)$. Moreover, let $\chi_n:\R_+\rightarrow[0,1]$ be as in the proof of 1. Then by inserting $\varphi=\chi_n(H)\theta_0'(R)$ in the definition of weak solution we obtain
\begin{align*}
\int_{\R^2_+}(\partial_Ru\theta_0''+f''(\theta_0)u\theta_0')\chi_n+\partial_Hu\chi_n'\theta_0'\,d(R,H)=
\int_{\R^2_+}G\theta_0'\chi_n\,d(R,H)+\int_{\R}g\theta_0'\,dR.
\end{align*}
The first term on the left hand side vanishes since we can apply integration by parts in $R$ and the second term converges to $0$ for $n\rightarrow\infty$ because of the Dominated Convergence Theorem since $\partial_Hu\theta_0'\in L^1(\R^2_+)$. Because of $G\theta_0'\in L^1(\R^2_+)$ the latter theorem applied to the first integral on the right hand side yields the compatibility condition \eqref{eq_hp_intbed}.\hfill$\square_{4.}$\\
\newline
\textit{Ad 5.} Let $G\in L^2(\R^2_+)$ with $G\theta_0'\in L^1(\R^2_+)$, $g\in L^2(\R)$ and $\tilde{G}(H):=(G(.,H),\theta_0')_{L^2(\R)}$ for $H>0$. By Fubini's Theorem $\tilde{G}$ is well-defined a.e. on $\R_+$ and belongs to $L^1(\R_+)$. $\tilde{G}\in L^2(\R_+)$ follows from Lemma \ref{th_prod_sets_spaces}, 1. and linearity of $(.,\theta_0')_{L^2(\R)}:L^2(\R)\rightarrow\R$. We define $G^\perp$ and $g^\perp$ according to equations \eqref{eq_hp_orth_decomp}. The claimed properties can be directly verified.\hfill$\square_{5.}$\\
\newline
\textit{Ad 6.} Let $G\in L^2(\R^2_+)$ with $\|G(.,H)\|_{L^2(\R)}\leq Ce^{-\nu H}$ for $C,\nu>0$ and $g\in L^2(\R)$. First, we show $G\theta_0'\in L^1(\R^2_+)$. Because of Lemma \ref{th_prod_sets_spaces}, 1. and since multiplication with $\theta_0'$ gives a bounded, linear operator from $L^2(\R)$ to $L^1(\R)$, we know that $\R_+\ni H\mapsto G(.,H)\theta_0'\in L^1(\R)$ is strongly measurable. The estimate for $G$ ensures $G\theta_0'\in L^1(\R_+,L^1(\R))\cong L^1(\R^2_+)$.\\
Therefore we can define $\tilde{G}, G^\perp$ and $g^\perp$ as in 4. We show that $\tilde{u}$ defined via \eqref{eq_hp_sol_formula} is well-defined. Since $\tilde{G}\in L^1(\R_+)\cap L^2(\R_+)$ and $|\tilde{G}(H)|\leq Ce^{-\nu H}$ f.a.a. $H>0$, we obtain 
\[
\int_H^\infty\tilde{G}(\hat{H})\,d\hat{H}=\int_0^\infty\tilde{G}(\hat{H})\,d\hat{H}-\int_0^H\tilde{G}(\hat{H})\,d\hat{H}\in W^1_1(\R_+)\cap H^1(\R_+)\quad\text{ w.r.t. }H
\] 
with derivative $-\tilde{G}$. Analogously, 
\[
\int_H^\infty \int_{\tilde{H}}^\infty\tilde{G}(\hat{H})\,d\hat{H}\,d\tilde{H}\in W^2_1(\R_+)\cap H^2(\R_+)\quad\text{ w.r.t. }H
\] 
and the derivative is given by $-\int_H^\infty\tilde{G}(\hat{H})\,d\hat{H}$. Since $\theta_0'\in W^2_1(\R)\cap H^2(\R)$, Lemma \ref{th_prod_sets_spaces}, 3. yields $\tilde{u}\in W^2_1(\R^2_+)\cap H^2(\R^2_+)$ and we can explicitly calculate the derivatives. One can directly verify that $(-\Delta+f''(\theta_0))\tilde{u}=\beta^2\theta_0'\tilde{G}$. Moreover, let the compatibility condition \eqref{eq_hp_intbed} hold. Then
\[
-\partial_H\tilde{u}|_{H=0}=-\beta^2\theta_0'\int_0^\infty\tilde{G}(\hat{H})\,d\hat{H}=-\beta^2\theta_0'\int_{\R^2_+}G\theta_0'\,d(R,H)=\beta^2(g,\theta_0')_{L^2(\R)}\theta_0'.
\]
Therefore $\tilde{u}$ is a strong solution of \eqref{eq_hp_model1}-\eqref{eq_hp_model2} for $G-G^\perp,g-g^\perp$ instead of $G,g$ and thus a weak solution because of integration by parts.\hfill$\square_{6.}$\\
\subsubsection{Exponential Decay Estimates}\label{sec_hp_exp_decay}
We need the following auxiliary lemma:
\begin{Lemma}\label{th_hp_exp_decay_lemma}
Let $f,g\in H^1(\R^2_+)$ and define $h(H):=\int_{\R}(fg)(R,H)\,dR$ for $H>0$. Then $h$ is well-defined a.e. on $\R_+$, $h\in W^1_1(\R_+)$ and $h'(H)=\int_{\R}(\partial_Hfg+f\partial_Hg)(R,H)\,dR$. The analogous assertion holds for the other variable.
\end{Lemma}
\textit{Proof.} The Hölder inequality implies $fg\in L^1(\R^2_+)$ and Fubini's Theorem shows that $h\in L^1(\R_+)$. One can verify the definition of weak derivative via density of $C_{(0)}^\infty(\overline{\R^2_+})$-functions in $H^1(\R^2_+)$, integration by parts and Fubini's Theorem, where we define $C_{(0)}^\infty(\overline{\R^2_+}):=\{\varphi|_{\R^2_+}:\varphi\in C_0^\infty(\R^2) \}$. Similarly, this can be done for the other variable.\hfill$\square$\\
\begin{Theorem}\label{th_hp_exp_decay1}
Let $g\in H^{\frac{1}{2}}(\R), G\in L^2(\R^2_+)$ and $\|G(.,H)\|_{L^2(\R)}\leq Ce^{-\nu H}$ f.a.a. $H>0$, where $\nu\in(0,\sqrt{\nu_2})$. Moreover, let the compatibility condition \eqref{eq_hp_intbed} hold. Then by Corollary \ref{th_hp_weak_sol_2} there is a unique solution $u\in H^2(\R^2_+)$ of \eqref{eq_hp_model1}-\eqref{eq_hp_model2}. It holds
\begin{enumerate}
\item $\|u(.,H)\|_{L^2(\R)}\leq C_u e^{-\nu H}$ f.a.a. $H\in\R_+$.
\item If $\|G(R,.)\|_{L^2(\R_+)}, |g(R)|\leq Ce^{-\nu|R|}$ f.a.a. $R\in\R$, then $\|u(R,.)\|_{L^2(\R_+)}\leq C_ue^{-\nu|R|}$ f.a.a. $R\in\R$.
\end{enumerate}
For 2. the compatibility condition \eqref{eq_hp_intbed} is not needed, if we require the existence of a weak solution $u\in H^2(\R^2_+)$ instead.
\end{Theorem}
\textit{Proof. Ad 1.} We consider $w(H):=\|u(.,H)\|_{L^2(\R)}$ for $H>0$. Lemma \ref{th_hp_exp_decay_lemma} yields $w^2\in W^2_1(\R_+)\hookrightarrow C_b^1(\R_+)$, in particular $w\in C_b^0(\R_+)$ and $w$ is $C^1$ on $\{w>0\}$ with
\[
\partial_Hw(H)=\frac{1}{w(H)}\int_{\R}(\partial_Hu u)(R,H)\,dR,\quad\partial_H(\frac{1}{w})(H)=-\frac{\partial_Hw(H)}{w(H)^2}\quad\text{ for }w(H)>0,
\]
where the integral belongs to $W^1_1(\R_+)\hookrightarrow C_b^0(\R_+)$. Therefore for $c_0>0$ arbitrary, it follows that $w, \frac{1}{w}\in C_b^1(\{w>c_0\})$. This implies $w\in W^2_1((0,R_0)\cap\{w>c_0\})$ for all $R_0>0$ with
\[
\partial_H^2w(H)=\frac{w(H)\int_{\R}(\partial_H^2uu+(\partial_Hu)^2)(R,H)\,dR-\frac{1}{w(H)}(\int_{\R}(\partial_Huu)(R,H)\,dR)^2}{w^2(H)}
\]
f.a.a. $H\in(0,\infty)\cap\{w>c_0\}$. The Hölder inequality yields 
\[
\left|\int_{\R}(\partial_Huu)(R,H)\,dR\right|\leq\|\partial_Hu(.,H)\|_{L^2(\R)}w(H)\quad\text{ f.a.a. }H>0.
\] 
Hence we obtain
\[
\partial_H^2w(H)\geq\frac{1}{w(H)}\int_{\R}(\partial_H^2uu)(R,H)\,dR\quad\text{ f.a.a. }H\in(0,\infty)\cap\{w>c_0\}.
\]
We can assume $G(.,H), g\perp\theta_0'$ in $L^2(\R)$ f.a.a. $H>0$, because $G^\perp, g^\perp$ in Theorem \ref{th_hp_weak_sol}, 4. admit the same properties as $G,g$ in the formulation of Theorem \ref{th_hp_exp_decay1} and the parallel part $\tilde{u}$ of $u$ from Theorem \ref{th_hp_weak_sol}, 6. has the claimed exponential decay. Now we use $\partial_H^2u=-G+f''(\theta_0)u-\partial_R^2u$. By Theorem \ref{th_hp_weak_sol}, 2. we have $u(.,H)\perp\theta_0'$ f.a.a. $H\in\R_+$. Hence Lemma \ref{th_prod_sets_spaces} and Lemma \ref{th_ODE_lin1} imply
\[
\int_{\R}(f''(\theta_0(R))u-\partial_R^2u)u(R,H)\,dR\geq\nu_2 w^2(H)\quad\text{ f.a.a. }H>0.
\]
Moreover, the Hölder inequality yields $|\int_{\R}(Gu)(R,H)\,dR|\leq Ce^{-\nu H}w(H)$ f.a.a. $H>0$. Thus
\[
-\partial_H^2w(H)+\nu_2w(H)\leq Ce^{-\nu H}\quad\text{ f.a.a. }H\in(0,\infty)\cap\{w>c_0\},
\]
where $c_0>0$ was arbitrary.\\
To derive an estimate from this ordinary differential inequality, we compare it with a function that satisfies equality. One can directly verify that
\[
\tilde{w}(H):=\frac{C}{-\nu^2+\nu_2}e^{-\nu H}+w(0)e^{-\sqrt{\nu_2}H}
\]
is such a function and $\tilde{w}(0)=\frac{C}{-\nu^2+\nu_2}+w(0)>w(0)$. The following arguments are motivated by Bartels \cite{Bartels_spectrum}, Appendix B. For $E:=w-\tilde{w}\in C^0_b(\R_+)$ we have to show $E\leq0$. Assume there is an $x\geq0$ such that $E(x)>0$. Since $E(0)<0$, $x$ must be larger than zero. Let $(\alpha,\beta)$ be the maximal interval around $x$ such that $E|_{(\alpha,\beta)}>0$, where $0<\alpha<\beta\leq\infty$. By continuity $E(\alpha)=0$. For finite $\tilde{\beta}\in(\alpha,\beta]$ we have
\[
w|_{[\alpha,\tilde{\beta}]}\geq\tilde{w}|_{[\alpha,\tilde{\beta}]}\geq\frac{C}{-\nu^2+\nu_2}e^{-\nu\tilde{\beta}}=:2c_0(\tilde{\beta}).
\]
The above arguments with $c_0=c_0(\tilde{\beta})$ imply $w\in W^2_1((\alpha,\tilde{\beta}))$ and $-\partial_H^2w+\nu_2w\leq Ce^{-\nu H}$ on $(\alpha,\tilde{\beta})$. Therefore $E\in W^2_1((\alpha,\tilde{\beta}))$ with $E''\geq\nu_2 E\geq 0$ and hence $E'\geq E'|_\alpha\geq 0$ on $[\alpha,\tilde{\beta}]$ for all finite $\tilde{\beta}\in(\alpha,\beta]$, where we used that $E$ is $C^1$ in a neighbourhood of $\alpha$. This yields
\[
E\geq E(x)>0\quad\text{ on }[x,\tilde{\beta}]\text{ for all finite }\tilde{\beta}\in(x,\beta].
\]
Necessarily $\beta=\infty$ and $E\geq E(x)>0$ on $[x,\beta)$. Since $\tilde{w}(H)\overset{H\rightarrow\infty}{\longrightarrow}0$, this implies the existence of an $H_0>0$ such that $w\geq\frac{1}{2}E(x)>0$ on $[H_0,\infty)$ which is a contradiction to $w^2\in W^2_1(\R_+)$. Hence $w(H)\leq\tilde{w}(H)\leq C_ue^{-\nu H}$ for all $H\in\R_+$.\hfill$\square_{1.}$\\
\newline
\textit{Ad 2.} Let $\|G(R,.)\|_{L^2(\R_+)},|g(R)|\leq Ce^{-\nu|R|}$ f.a.a $R\in\R$ and $w(R):=\|u(R,.)\|_{L^2(\R_+)}$ for $R\in\R$. Analogously as in the first part it follows that $w^2\in W^2_1(\R)\hookrightarrow C_b^1(\R)$, in particular $w\in C_b^0(\R)$ and $w$ is $C^1$ on $\{w>0\}$ with
\[
\partial_Rw(R)=\frac{1}{w(R)}\int_{\R_+}\partial_Ruu(R,H)\,dH.
\] 
Similarly as before $w,\frac{1}{w}\in C_b^1(\{w>c_0\})$ for all $c_0>0$ and $w\in W^2_1((-R_0,R_0)\cap\{w>c_0\})$ for all $R_0>0$ and it holds
\[
\partial_R^2w(R)\geq\frac{1}{w(R)}\int_{\R_+}\partial_R^2uu(R,H)\,dH\quad\text{ f.a.a. }R\in\{w>c_0\}.
\]
Now we use the equation $\partial_R^2u=-G+f''(\theta_0)u-\partial_H^2u$. First of all, the Hölder inequality yields
\[
\left|\int_{\R_+}Gu(R,H)\,dH\right|\leq Ce^{-\nu|R|}w(R)\quad\text{ f.a.a. }R\in\R.
\]
Moreover, $\int_{\R_+}f''(\theta_0(R))u(R,H)^2\,dH=f''(\theta_0(R))w^2(R)$ f.a.a. $R\in\R$ and
\[
-\int_{\R_+}\partial_H^2uu(R,H)\,dH=\int_{\R_+}(\partial_Hu)^2(R,H)\,dH+(\partial_Hu\,u)|_{H=0}(R),
\]
where the latter identity can be shown by multiplying with $\eta=\eta(R)\in C_0^\infty(\R)$, integration over $R\in\R$, integration by parts and the Fundamental Theorem of Calculus of Variations. Altogether we obtain
\[
-\partial_R^2w(R)+f''(\theta_0(R))w(R)\leq Ce^{-\nu|R|}-\frac{1}{w(R)}\left(\int_{\R_+}(\partial_Hu)^2(R,H)\,dH+(\partial_Huu)|_{H=0}(R)\right)
\]
f.a.a. $R\in\{w>c_0\}$. To simplify the last two terms, we use
\[
|\tr_{H=0}u(R)|^2\leq 2w(R)\|\partial_Hu(R,.)\|_{L^2(\R_+)}\quad\text{ f.a.a. }R\in\R,
\] 
which follows from $u(R,0)^2=2\int_0^\infty\partial_Huu(R,H)\,dR$ for $u\in C_{(0)}^\infty(\overline{\R^2_+})$, Hölder's inequality and a density argument. Young's inequality and $-\partial_Hu|_{H=0}=g$ yields
\[
\frac{1}{w(R)}|\partial_Huu(R,0)|\leq \frac{1}{w(R)}e^{-\nu|R|}(C_\delta w(R)+\delta \|\partial_Hu(R,.)\|_{L^2(\R_+)})\quad\text{ for }\delta>0.
\]
Furthermore, $e^{-\nu|R|}\|\partial_Hu(R,.)\|_{L^2(\R_+)}\leq\frac{1}{2}(e^{-2\nu|R|}+\|\partial_Hu(R,.)\|_{L^2(\R_+)}^2)$. For small $\delta$ the $\|\partial_Hu(R,.)\|_{L^2(\R_+)}^2$-term is absorbed and we obtain
\[
-\partial_R^2w(R)+f''(\theta_0(R))w(R)\leq\tilde{c}e^{-\nu|R|}+\overline{c}\frac{e^{-2\nu|R|}}{w(R)}\quad\text{ f.a.a. }R\in\{w>c_0\}.
\]
To apply the comparison argument, let $R_0>0$ such that $f''(\theta_0(R))\geq\nu_2-\varepsilon>\nu^2$ for all $|R|\geq R_0$. We seek a function $\tilde{w}$ that fulfills equality in $\R\textbackslash\{0\}$ for $f''(\theta_0)$ replaced by $\nu_2-\varepsilon$ with the ansatz $\tilde{w}(R)=c_1e^{-\nu|R|}$. Equality is equivalent to
\[
\tilde{\nu}c_1:=\left[-\nu^2+\nu_2-\varepsilon\right]c_1=\tilde{c}+\frac{\overline{c}}{c_1}\quad\Leftrightarrow\quad\tilde{\nu}c_1^2-\tilde{c}c_1-\overline{c}=0,
\]
hence satisfied if we choose $c_1:=\frac{\tilde{c}+\sqrt{\tilde{c}^2+4\tilde{\nu}\overline{c}}}{2\tilde{\nu}}>0$. W.l.o.g. we can assume that $\tilde{w}|_{[-R_0,R_0]}>w|_{[-R_0,R_0]}$, otherwise we increase $\overline{c}$. Since $\tilde{w}(R)\leq w(R)$ implies $\frac{1}{w(R)}\leq\frac{1}{\tilde{w}(R)}$, for such $R$ we can estimate the right hand side and argue as in the first part to obtain $w\leq\tilde{w}$ on $\R\textbackslash(-R_0,R_0)$. Note that we did not use the compatibility condition \eqref{eq_hp_intbed} explicitly, just \eqref{eq_hp_model1}-\eqref{eq_hp_model2} are needed.\hfill$\square_{2.}$\\
\newline
The following theorem is concerned with similar estimates for the derivatives.
\begin{Theorem}\label{th_hp_exp_decay2}
Let $m\in\N$, $G\in H^m(\R^2_+), g\in H^{m+\frac{1}{2}}(\R)$ and for some $\nu\in(0,\sqrt{\nu_2})$ let $\|\partial_R^l\partial_H^kG(.,H)\|_{L^2(\R)}=\Oc(e^{-\nu H})$ for all $0\leq k+l\leq m$. Moreover, let the compatibility condition \eqref{eq_hp_intbed} hold and $u\in H^{m+2}(\R^2_+)$ be the solution of \eqref{eq_hp_model1}-\eqref{eq_hp_model2}. Then
\begin{enumerate}
\item $\|\partial_R^l\partial_H^ku(.,H)\|_{L^2(\R)}=\Oc(e^{-\nu H})$ for all $0\leq k+l\leq m-1$, if $l$ is odd, and  all $0\leq k+l\leq m$, if $l$ is even.
\item If $\|\partial_R^l\partial_H^kG(R,.)\|_{L^2(\R_+)}, |\partial_R^lg(R)|=\Oc(e^{-\nu|R|})$ holds for all $0\leq k+l\leq m$, then 
\[
\|\partial_R^l\partial_H^ku(R,.)\|_{L^2(\R_+)}=\Oc(e^{-\nu|R|})
\] 
for all $0\leq k+l\leq m-1$, if $k$ is odd, and all $0\leq k+l\leq m$ if $k$ is even.
\end{enumerate}
\end{Theorem}
Altogether we obtain
\begin{Corollary}\label{th_hp_exp_decay3}
Let $G\in C_b^\infty(\overline{\R^2_+}), g\in C_b^\infty(\R)$ with $|\partial^\alpha G(R,H)|=\Oc(e^{-\nu(|R|+H)})$ and $|\partial_R^\beta g(R)|=\Oc(e^{-\nu|R|})$ for all $\alpha\in\N_0^2, \beta\in\N_0$, where $\nu\in(0,\sqrt{\nu_2})$. Moreover, let the compatibility condition \eqref{eq_hp_intbed} hold. Then the unique solution $u\in H^1(\R^2_+)$ to \eqref{eq_hp_model1}-\eqref{eq_hp_model2} satisfies $u\in C_b^\infty(\overline{\R^2_+})$ and $|\partial^\gamma u(R,H)|=\Oc(e^{-\frac{\nu}{2}(|R|+H)})$ for all $\gamma\in\N_0^2$.
\end{Corollary}
\textit{Proof.} The assumptions ensure $G\in H^m(\R^2_+), g\in H^{m+\frac{1}{2}}(\R)$ for all $m\in\N_0$ and Corollary \ref{th_hp_weak_sol_2} implies the existence of a unique solution $u\in H^{m+2}(\R^2_+)\hookrightarrow C^{m,\gamma}(\overline{\R^2_+})$ of \eqref{eq_hp_model1}-\eqref{eq_hp_model2} for all $m\in\N_0, \gamma\in(0,1)$. Moreover, Theorem \ref{th_hp_exp_decay2} together with $H^1(Y)\hookrightarrow C_b^0(Y)$ for $Y\in\{\R,\R_+\}$ and $r=\sqrt{r}\sqrt{r}$ for $r\in\R$ yields the estimate.\hfill$\square$\\
\newline
\textit{Proof of Theorem \ref{th_hp_exp_decay2}. Ad 1.} First we consider the case $l=0$. Because of Theorem \ref{th_hp_exp_decay1}, 1. we can assume $1\leq k\leq m$. In order to apply Theorem \ref{th_hp_exp_decay1}, 1. to $\partial_H^ku$ we differentiate through the equations:
\begin{align*}
-\partial_R^2\partial_H^ku-\partial_H^{2+k}u+f''(\theta_0(R))\partial_H^ku&=\partial_H^kG,\\
-\partial_H\partial_H^ku|_{H=0}&=(\partial_R^2\partial_H^{k-1}u-f''(\theta_0)\partial_H^{k-1}u+\partial_H^{k-1}G)|_{H=0}=:g_k.
\end{align*} 
It holds $\partial_H^kG\in L^2(\R^2_+)$ and $\|\partial_H^kG(.,H)\|_{L^2(\R)}=\Oc(e^{-\nu H})$. Furthermore, $g_k\in H^{\frac{1}{2}}(\R)$ follows from $f''(\theta_0)\in C_b^1(\R)$ and the trace theorem. Therefore Theorem \ref{th_hp_exp_decay1}, 1. yields $\|\partial_H^ku(.,H)\|_{L^2(\R)}=\Oc(e^{-\nu H})$, if we show the associated compatibility condition
\begin{align}\label{eq_hp_exp_decay_intbed}
\int_{\R^2_+}\partial_H^kG(R,H)\theta_0'(R)\,d(R,H)+\int_{\R}g_k(R)\theta_0'(R)\,dR=0.
\end{align}
Let $\chi_n:\R_+\rightarrow[0,1]$ for $n\in\N$ be as in the proof of Theorem \ref{th_hp_weak_sol}, 1. Then integration by parts and the Dominated Convergence Theorem yields
\begin{align*}
\int_{\R^2_+}\partial_H^kG\,\theta_0'\,d(R,H)&=\lim_{n\rightarrow\infty}\int_{\R^2_+}\partial_H^kG\theta_0'\chi_n\,d(R,H)=\\
&=\lim_{n\rightarrow\infty}-\int_{\R^2_+}\partial_H^{k-1}G\,\theta_0'\chi_n'\,d(R,H)-\int_{\R}\partial_H^{k-1}G|_{H=0}\theta_0'\,dR.
\end{align*}
The first integral on the right hand side vanishes by the Dominated Convergence Theorem since $\partial_H^{k-1}G\,\theta_0'\in L^1(\R^2_+)$, cf. Theorem \ref{th_hp_weak_sol}, 6. It remains to show
\[
\int_{\R}(\partial_R^2\partial_H^{k-1}u-f''(\theta_0)\partial_H^{k-1}u)|_{H=0}(R)\theta_0'(R)\,dR=0.
\]
By density of $C_{(0)}^\infty(\overline{\R^2_+})$ in $H^3(\R^2_+)$ one can show $(\partial_R^2\partial_H^{k-1}u)|_{H=0}=\frac{d^2}{dR^2}(\partial_H^{k-1}u|_{H=0})$. Hence integration by parts yields \eqref{eq_hp_exp_decay_intbed} and the case $l=0$ is done.\\
For $m=0,1$ there is nothing left to prove. Now let $m\geq 2$ and $0\leq k\leq m-2$. Then differentiating the equation yields
\[
\partial_R^2\partial_H^ku=-\partial_H^{k+2}u+f''(\theta_0(R))\partial_H^ku-\partial_H^kG
\]
and this shows the case $l=2$. To cover $l=1$ we use the following interpolation inequality: Integration by parts and the Hölder inequality yield that there is a $C>0$ such that
\[
\|u'\|_{L^2(\R)}\leq C\|u\|_{L^2(\R)}^{\frac{1}{2}}\|u''\|_{L^2(\R)}^{\frac{1}{2}}\quad\text{ for all }u\in H^2(\R).
\]
Since $\partial_H^ku\in H^2(\R^2_+)\hookrightarrow L^2(\R_+,H^2(\R))$ and the derivatives are compatible because of Lemma \ref{th_prod_sets_spaces}, it follows that
\[
\|\partial_R\partial_H^ku(.,H)\|_{L^2(\R)}=\Oc(e^{-\nu H})\quad\text{ for }k=0,...,m-2.
\]
For $m=2,3$ we are done. The general case follows inductively by differentiating the equation for even $l$ and the interpolation inequality.\hfill$\square_{1.}$\\
\newline
\textit{Ad 2.} Let additionally $\|\partial_R^l\partial_H^kG(R,.)\|_{L^2(\R_+)}, |\partial_R^lg(R)|=\Oc(e^{-\nu|R|})$ f.a.a. $R\in\R$ hold for all $0\leq k+l\leq m$. First, we consider the case $k=0$. We can assume $1\leq l\leq m$ because of Theorem \ref{th_hp_exp_decay1}, 2. Differentiating the equations \eqref{eq_hp_model1}-\eqref{eq_hp_model2} yields
\begin{align}\label{eq_hp_exp_decay1}
-\partial_R^{2+l}u-\partial_R^l\partial_H^2u+f''(\theta_0)\partial_R^lu&=\partial_R^lG+\sum_{j=1}^l\begin{pmatrix}
l\\j\end{pmatrix}\partial_R^j(f''(\theta_0))\partial_R^{l-j}u=:G_k,\\
-\partial_H\partial_R^lu|_{H=0}&=\partial_R^lg.
\label{eq_hp_exp_decay2}
\end{align}
For this system we use Theorem \ref{th_hp_exp_decay1}. Note that the corresponding compatibility condition is not needed. Because of $\partial_R^lg\in H^{\frac{1}{2}}(\R)$ with $|\partial_R^lg(R)|=\Oc(e^{-\nu|R|})$ and $G_k\in L^2(\R^2_+)$, it remains to show
\[
\|G_k(.,H)\|_{L^2(\R)}=\Oc(e^{-\nu H})\quad\text{ and }\quad\|G_k(R,.)\|_{L^2(\R_+)}=\Oc(e^{-\nu|R|}).
\]
For $\partial_R^lG$ these estimates hold by assumption. Since $l-j\leq m-1$ in \eqref{eq_hp_exp_decay1}, the first part yields $\|\partial_R^{l-j}u(.,H)\|_{L^2(\R)}=\Oc(e^{-\nu H})$ and the $(.,H)$-estimate for $G_k$ follows. Now the $(R,.)$-estimates for $\partial_R^lu$ can be shown inductively with Theorem \ref{th_hp_exp_decay1}.\\
For $m=0,1$ there is nothing left to prove, thus let $m\geq 2$. Next we show all the estimates for $k$ even. The estimates for $k=2, 0\leq l\leq m-2$ are obtained using the differentiated equation \eqref{eq_hp_exp_decay1}. The cases $k=4,6,...$ can be shown inductively by further differentiating the equation in $H$.\\ 
Now consider the case $k=1$ and $0\leq l\leq m-2$. As in the proof of Theorem \ref{th_hp_exp_decay1}, 2. it follows that f.a.a. $R\in\R$
\[
\|\partial_H\partial_R^lu(R,.)\|_{L^2(\R_+)}^2=-\int_{\R_+}\partial_H^2\partial_R^lu\partial_R^lu(R,H)\,dH-\partial_H\partial_R^lu|_{H=0}(R)\partial_R^lu|_{H=0}(R).
\]
The first term on the right hand side is $\Oc(e^{-2\nu|R|})$ because of Hölder's inequality and the estimates proven above. Furthermore, as in the proof of Theorem \ref{th_hp_exp_decay1}, 2. we have
\[
|\partial_R^lu|_{H=0}(R)|\leq 2^{\frac{1}{2}}\|\partial_R^lu(R,.)\|_{L^2(\R_+)}^{\frac{1}{2}}\|\partial_H\partial_R^lu(R,.)\|_{L^2(\R_+)}^{\frac{1}{2}}\quad\text{ f.a.a. }R\in\R.
\]
Since $\partial_H\partial_R^lu|_{H=0}(R)=\partial_R^lg(R)=\Oc(e^{-\nu|R|})$ and $\|\partial_R^lu(R,.)\|_{L^2(\R_+)}=\Oc(e^{-\nu|R|})$, this yields
\[
\|\partial_H\partial_R^lu(R,.)\|_{L^2(\R_+)}^2\leq Ce^{-2\nu|R|}+Ce^{-\frac{3}{2}\nu|R|}\|\partial_H\partial_R^lu(R,.)\|_{L^2(\R_+)}^{\frac{1}{2}}.
\] 
Young's inequality with $p=4, p'=\frac{4}{3}$ applied to the last term shows the desired estimate.\\
For $m=2,3$ we are done, hence let $m\geq 4$. The estimates for $k=3, 0\leq l\leq m-4$ follow from \eqref{eq_hp_exp_decay1} after differentiating in $H$ using the previous estimates. The cases $k=5,7,...$ are obtained inductively by further differentiating in $H$.\hfill$\square_{2.}$\\
\subsubsection{Dependence on Time}\label{sec_hp_time}
In the following theorem we consider the situation when $G,g$ also depend on $t\in[0,T]$. 
\begin{Theorem}\label{th_hp_time}
Let $m,n\in\N, G\in C^n([0,T],H^m(\R^2_+)), g\in C^n([0,T],H^{m+\frac{1}{2}}(\R))$ with $\|\partial_t^jG(.,H,t)\|_{L^2(\R)}\leq Ce^{-\nu H}$ f.a.a. $H>0$ and all $j=0,...,n$ with $C>0$ independent of $t\in[0,T]$. Moreover, let the compatibility condition \eqref{eq_hp_intbed} hold. Then \eqref{eq_hp_model1}-\eqref{eq_hp_model2} has exactly one solution $u\in C^n([0,T],H^{m+2}(\R^2_+))$.\\
Furthermore, the assertions in Theorem \ref{th_hp_exp_decay1} and Theorem \ref{th_hp_exp_decay2} are satisfied for $u$ and $\partial_t^ju, j=1,...,n$ instead of $u$, when all the estimates are replaced by uniform ones in $t\in[0,T]$ and, if $j=1,...,n$, $G,g$ are replaced by $\partial_t^jG,\partial_t^jg$.
\end{Theorem} 
\textit{Proof.} For fixed $t\in[0,T]$ Corollary \ref{th_hp_weak_sol_2} yields a unique solution $u(t)$ in $H^{m+2}(\R^2_+)$ to \eqref{eq_hp_model1}-\eqref{eq_hp_model2} for $(G(t),g(t))$ instead of $(G,g)$. We consider the decomposition \eqref{eq_hp_orth_decomp} of $G,g$ in Theorem \ref{th_hp_weak_sol}, 5. and treat $t\in[0,T]$ as additional parameter. Then $\tilde{G}\in C^n([0,T],H^m(\R_+)), G^\perp\in C^n([0,T],H^m(\R^2_+))$ and $g^\perp\in C^n([0,T],H^{m+\frac{1}{2}}(\R))$ because of Lemma \ref{th_prod_sets_spaces}. By Theorem \ref{th_hp_weak_sol} it holds $u=\tilde{u}+u^\perp$, where $\tilde{u}$ is given by \eqref{eq_hp_sol_formula} and $u^\perp$ is a solution for the orthogonal parts in \eqref{eq_hp_orth_decomp}. Theorem \ref{th_hp_weak_sol}, 2. and Corollary \ref{th_hp_weak_sol_2} imply
\[
\|u^\perp(.,t)\|_{H^{m+2}(\R^2_+)}\leq C(\|G^\perp(.,t)\|_{H^m(\R^2_+)}+\|g^\perp(.,t)\|_{H^{m+\frac{1}{2}}(\R)})\quad\text{ for all }t\in[0,T].
\]
With difference quotients we obtain $u^\perp\in C^n([0,T],H^{m+2}(\R^2_+))$ and $\partial_t^j(u^\perp)$ solves \eqref{eq_hp_model1}-\eqref{eq_hp_model2} for $(\partial_t^jG^\perp,\partial_t^jg^\perp)$, $j=0,...,n$. The Dominated Convergence Theorem yields $\tilde{u}\in C^n([0,T],H^{m+2}(\R^2_+))$ and the derivatives are the canonical ones. Altogether we obtain $u\in C^n([0,T],H^{m+2}(\R^2_+))$ and that $\partial_t^ju$ solves \eqref{eq_hp_model1}-\eqref{eq_hp_model2} for $(\partial_t^jG,\partial_t^jg), j=0,...,n$.\\
Since the compatibility condition \eqref{eq_hp_intbed} is satisfied for $(\partial_t^jG,\partial_t^jg), j=0,...,n$, we can go into the proofs of Theorem \ref{th_hp_exp_decay1} and Theorem \ref{th_hp_exp_decay2} to get the additional estimates.\hfill$\square$\\
\begin{Corollary}\label{th_hp_time_cor}
Let $G\in C_b^\infty(\overline{\R^2_+}\times[0,T]), g\in C_b^\infty(\R\times[0,T])$ and for all $\alpha\in\N_0^3,\beta\in\N_0^2$ 
\[
|\partial^\alpha G(R,H,t)|\leq C_{\alpha}e^{-\nu(|R|+H)},\quad|\partial^\beta g(R,t)|\leq C_{\beta}e^{-\nu|R|}\quad\text{ for all }(R,H,t)\in\R^2_+\times[0,T],
\] 
where $\nu\in(0,\sqrt{\nu_2})$. Moreover, let the compatibility condition \eqref{eq_hp_intbed} hold. Then the unique solution $u$ to \eqref{eq_hp_model1}-\eqref{eq_hp_model2} fulfills $u\in C_b^\infty(\overline{\R^2_+}\times[0,T])$ and for all $\gamma\in\N_0^3$ 
\[
|\partial^\gamma u(R,H,t)|\leq C_{\gamma,u}e^{-\frac{\nu}{2}(|R|+H)}\quad\text{ for all }(R,H,t)\in\R^2_+\times[0,T].
\]
\end{Corollary}
\textit{Proof.} Because of Lemma \ref{th_prod_sets_spaces}, 2. and embedding theorems the conditions in Theorem \ref{th_hp_time} are satisfied for all $m,n\in\N$. Embedding theorems show the smoothness of the solution and the estimates are obtained as in the proof of Corollary \ref{th_hp_exp_decay3}.\hfill$\square$\\
\section{Asymptotic Expansion and Construction of the Approximate Solution}\label{sec_expansion}
Let $\Gamma:=(\Gamma_t)_{t\in[0,T]}$ be as in Section \ref{th_coord} and $\delta>0$ such that the assertions of Theorem \ref{th_coord} hold for $2\delta$ instead of $\delta$. Based on $\Gamma$ we construct an approximate solution $u_A^\varepsilon$ to \eqref{eq_AC1}-\eqref{eq_AC3} with $u_A^\varepsilon=\pm 1$ on $Q_T^\pm\textbackslash\Gamma(2\delta)$. Therefore we introduce height functions $h_1,h_2:I\times[0,T]\rightarrow\R$ and define a \enquote{stretched variable}
\[
\rho_\varepsilon(x,t):=\frac{r(x,t)}{\varepsilon}-h_\varepsilon(s(x,t),t)\quad\text{ for }(x,t)\in\Gamma(2\delta),\text{ where }h_\varepsilon:=h_1+\varepsilon h_2.
\]
The idea is that $\Gamma_t^\varepsilon:=\{x\in\Gamma_t(2\delta):\rho_\varepsilon(x,t)=0\}$ should approximate the zero-level set of the exact solution $u^\varepsilon(.,t)$ at time $t\in[0,T]$. We do not specify $h_j|_{t=0}$ at this point. These will be chosen later in a suitable way.
\subsection{Inner Expansion}\label{sec_inner_exp}
For the inner expansion we consider the following ansatz: Let $\varepsilon>0$ be small and
\[
u^I=\sum_{j=0}^3\varepsilon^ju_j^I,\quad u_j^I(x,t)=\hat{u}_j^I(\rho_\varepsilon(x,t),s(x,t),t)\quad\text{ for }(x,t)\in\Gamma(2\delta),
\]
where $\R\times I\times[0,T]\ni(\rho,s,t)\mapsto \hat{u}_j^I(\rho,s,t)\in\R$ for $j=0,...,3$ and we set $\hat{u}^I:=\sum_{j=0}^3\varepsilon^j\hat{u}_j^I$. By substituting $u^I$ into the Allen-Cahn equation while ignoring the Neumann boundary condition and expanding up to $\Oc(\varepsilon)$ (in order to fulfill the equation up to $\Oc(\varepsilon^2)$), we will derive $(s,t)$-dependent ODEs in $\rho$ of type as in Lemma \ref{th_ODE_lin2}. The corresponding solvability conditions will yield that $\Gamma$ should evolve according to mean curvature flow, the lowest order is given by the optimal profile $\theta_0$ and that the height functions should satisfy non-autonomous parabolic PDEs.\\
Therefore we need to know how the differential operators act on $u^I$:
\begin{Lemma}\label{th_trafo_inner}
Let $\hat{w}:\R\times I\times[0,T]\rightarrow\R$ be sufficiently smooth and $w:\Gamma(2\delta)\rightarrow\R$ be defined by $w(x,t):=\hat{w}(\rho_\varepsilon(x,t),s(x,t),t)$ for all $(x,t)\in\Gamma(2\delta)$. Then
\begin{align*}
\partial_tw&=\partial_\rho\hat{w}\left[\frac{\partial_tr}{\varepsilon}-(\partial_th_\varepsilon+\partial_ts\partial_sh_\varepsilon)\right]+\partial_s\hat{w}\,\partial_ts+\partial_t\hat{w},\\
\nabla w&=\partial_\rho\hat{w}\left[\frac{\nabla r}{\varepsilon}-\nabla s\partial_sh_\varepsilon\right]+\partial_s\hat{w}\,\nabla s,\\
\Delta w&=\partial_\rho\hat{w}\left[\frac{\Delta r}{\varepsilon}-(\Delta s\partial_sh_\varepsilon+|\nabla s|^2\partial_s^2h_\varepsilon)\right]+\partial_s\hat{w}\,\Delta s+\partial_s^2\hat{w}\,|\nabla s|^2+\\
&+2\partial_\rho\partial_s\hat{w}\,\nabla s\cdot\left[\frac{\nabla r}{\varepsilon}-\nabla s\partial_sh_\varepsilon\right]+\partial_\rho^2\hat{w}\left|\frac{\nabla r}{\varepsilon}-\nabla s\partial_sh_\varepsilon\right|^2,
\end{align*}
where the $w$-terms on the left hand side and derivatives of $r$ or $s$ are evaluated at $(x,t)\in\Gamma(2\delta)$, the $h_\varepsilon$-terms at $(s(x,t),t)$ and the $\hat{w}$-terms at $(\rho_\varepsilon(x,t),s(x,t),t)$.
\end{Lemma}
\textit{Proof.} This follows directly from the chain rule.\hfill$\square$\\
\newline
To derive suitable ODEs from the Allen-Cahn equation $\partial_tu-\Delta u+\frac{1}{\varepsilon^2}f'(u)=0$, we use Taylor expansions. First let us consider the $f'$-part. If the $u_j^I$ are uniformly bounded, then for $m\in\N$
\begin{align}\label{eq_inner_taylor_f}
f'(u^I)=f'(u_0^I)+\sum_{k=1}^m\frac{f^{(k+1)}(u_0^I)}{k!}\left[\sum_{j=1}^3 u_j^I\varepsilon^j\right]^k+\Oc(\varepsilon^{m+1})\quad\text{ on }\Gamma(2\delta).
\end{align}
The first terms in the $\varepsilon$-expansion are 
\begin{align*}
\Oc(1)&:\quad f'(u_0^I),\\
\Oc(\varepsilon)&:\quad f''(u_0^I)u_1^I,\\
\Oc(\varepsilon^2)&:\quad f''(u_0^I)u_2^I+\frac{f^{(3)}(u_0^I)}{2!}(u_1^I)^2,\\
\Oc(\varepsilon^3)&:\quad f''(u_0^I)u_3^I+\frac{f^{(3)}(u^I_0)}{2!}2u_1^Iu_2^I+\frac{f^{(4)}(u_0^I)}{3!}(u_1^I)^3.
\end{align*}
The derivatives of $r$ and $s$ are functions of $(x,t)\in\overline{\Gamma(2\delta)}$ and we will expand them with a Taylor expansion using $r(x,t)=\varepsilon(\rho_\varepsilon(x,t)+h_\varepsilon(s(x,t),t)$ for $(x,t)\in\overline{\Gamma(2\delta)}$. For smooth functions $g:\overline{\Gamma(2\delta)}\rightarrow\R$ the Taylor expansion yields uniformly in $(s,t)$: For $r\in[-2\delta,2\delta]$
\begin{align}\label{eq_inner_taylor}
\tilde{g}(r,s,t):=g(\overline{X}(r,s,t))=\sum_{k=0}^m\frac{\partial_r^k\tilde{g}|_{(0,s,t)}}{k!}r^k+\Oc(r^{m+1}).
\end{align}
The first few terms needed are
\begin{align*}
\Oc(1)&:\quad g|_{\overline{X}_0(s,t)},\\
\Oc(\varepsilon)&:\quad (\rho+h_1(s,t))\partial_r\tilde{g}|_{(0,s,t)},\\
\Oc(\varepsilon^2)&:\quad h_2(s,t)\partial_r\tilde{g}|_{(0,s,t)}+(\rho+h_1(s,t))^2\frac{\partial_r^2\tilde{g}|_{(0,s,t)}}{2},\\
\Oc(\varepsilon^3)&:\quad \frac{\partial_r^2\tilde{g}|_{(0,s,t)}}{2}2(\rho+h_1(s,t))h_2(s,t)+\frac{\partial_r^3\tilde{g}|_{(0,s,t)}}{3!}(\rho+h_1(s,t))^3,
\end{align*}
where we replaced $\rho_\varepsilon$ by an arbitrary $\rho\in\R$ in order to get ODEs on $\R$. But later we just use the expansion rigorously for $r=\varepsilon(\rho_\varepsilon+h_\varepsilon)\in[-2\delta,2\delta]$. The other terms without the $\Oc(r^{m+1})$-remainder are $\varepsilon^4$ times some polynomial in $\rho_\varepsilon$ with smooth coefficients in $(s,t)\in I\times[0,T]$, if the $h_j$ are smooth. Later, the $\Oc(r^m)$-terms will be multiplied with terms that decay exponentially in $|\rho|$. Then these remainder terms will become $\Oc(\varepsilon^m)$.\\ 
Now we expand the Allen-Cahn equation $\partial_tu-\Delta u+\frac{1}{\varepsilon^2}f'(u)=0$ for $u=u^I$ into $\varepsilon$-series. If we write down an equation or assertion for $(\rho,s,t)$ in the following, it is meant to hold for all $(\rho,s,t)\in\R\times I\times[0,T]$. Moreover, we often omit the argument $(s,t)$ in the $h_j$-terms.
\subsubsection{Inner Expansion: $\Oc(\frac{1}{\varepsilon^2})$}\label{sec_inner_1_eps2}
We obtain that the $\Oc(\frac{1}{\varepsilon^2})$-order is zero if
\[
\quad -|\nabla r|^2|_{\overline{X}_0(s,t)}\partial_\rho^2\hat{u}_0^I(\rho,s,t)+f'(\hat{u}_0^I(\rho,s,t))=0.
\]
Because of Lemma \ref{th_coord}, we have $|\nabla r|^2|_{\overline{X}_0(s,t)}=1$ and since $\{\rho_\varepsilon=0\}$ should approximate the zero level set of $u^I$, we require $\hat{u}_0^I(0,s,t)=0$. Together with $\lim_{\rho\rightarrow\pm\infty}\hat{u}_0^I(\rho,s,t)=\pm1$, we obtain $\hat{u}_0^I(\rho,s,t)=\theta_0(\rho)$ by Lemma \ref{th_theta_0}.
\subsubsection{Inner expansion: $\Oc(\frac{1}{\varepsilon})$}\label{sec_inner_1_eps}
We consider the $\Oc(\frac{1}{\varepsilon})$-order: From the $\partial_tu$-part we get $\frac{1}{\varepsilon}\partial_tr|_{\overline{X}_0(s,t)}\theta_0'(\rho)$ and from $\Delta u$:
\begin{align*}
\frac{1}{\varepsilon^2}\left[\partial_r(|\nabla r|^2\circ \overline{X})|_{(0,s,t)}\varepsilon(\rho+h_1(s,t))\theta_0''(\rho)+|\nabla r|^2|_{\overline{X}_0(s,t)}\varepsilon\partial_\rho^2\hat{u}_1^I(\rho,s,t)\right]+\\
+\varepsilon\left[\theta_0'(\rho)\Delta r|_{\overline{X}_0(s,t)}+2(\nabla r\cdot\nabla s)|_{\overline{X}_0(s,t)}(\partial_s\theta_0'(\rho)-\partial_sh_1(s,t)\theta_0''(\rho))\right]=\\
=\frac{1}{\varepsilon}\left[\partial_\rho^2\hat{u}_1^I(\rho,s,t)+\theta_0'(\rho)\Delta r|_{\overline{X}_0(s,t)}\right],
\end{align*}
where we used Theorem \ref{th_coord}. Therefore the $\frac{1}{\varepsilon}$-order cancels if  
\[
\Lc\hat{u}_1^I(\rho,s,t)+\theta_0'(\rho)(\partial_tr-\Delta r )|_{\overline{X}_0(s,t)}=0,\quad\text{ where }\Lc=-\partial_\rho^2+f''(\theta_0).
\]
Because of Lemma \ref{th_ODE_lin2} this parameter-dependent ODE together with $\hat{u}_1^I(0,s,t)=0$ and boundedness in $\rho$ has a (unique) solution $\hat{u}_1^I$ if and only if $(\partial_tr-\Delta r)|_{\overline{X}_0(s,t)}=0$ which is equivalent to mean curvature flow by Theorem \ref{th_coord}. Thus we define $\hat{u}_1^I=0$.
\subsubsection{Inner Expansion: $\Oc(1)$}\label{sec_inner_1}
At $\Oc(1)$-order we obtain from $\partial_tu$:
\begin{align*}
\frac{1}{\varepsilon}\left[\partial_tr|_{\overline{X}_0(s,t)}\varepsilon\partial_\rho\hat{u}_1^I+\partial_r(\partial_tr\circ\overline{X})|_{(0,s,t)}\varepsilon(\rho+h_1(s,t))\theta_0'(\rho)\right]+\\
+\theta_0'(\rho)(-\partial_th_1(s,t)-\partial_ts|_{\overline{X}_0(s,t)}\partial_sh_1(s,t))+\partial_s\theta_0\partial_ts|_{\overline{X}_0(s,t)}+\partial_t\theta_0(\rho)=\\
=\theta_0'(\rho)\left[(\rho+h_1(s,t))\partial_r(\partial_tr\circ \overline{X})|_{(0,s,t)}-\partial_th_1(s,t)-\partial_ts|_{\overline{X}_0(s,t)}\partial_sh_1(s,t)\right],
\end{align*}
and from $\Delta u$:
\begin{align*}
\frac{1}{\varepsilon^2}\theta_0''(\rho)\left[
\varepsilon^2\frac{1}{2}(\rho+h_1)^2\partial_r^2(|\nabla r|^2\circ\overline{X})|_{(0,s,t)}+\varepsilon^2 h_2\partial_r(|\nabla r|^2\circ\overline{X})|_{(0,s,t)}\right]+\\
+\frac{1}{\varepsilon^2}\partial_\rho^2\hat{u}_1^I\varepsilon^2(\rho+h_1)\partial_r(|\nabla r|^2\circ\overline{X})|_{(0,s,t)}+\frac{1}{\varepsilon^2}|\nabla r|^2|_{\overline{X}_0(s,t)}\varepsilon^2\partial_\rho^2\hat{u}_2^I+\\
+\frac{1}{\varepsilon}\left[\theta_0'(\rho)\varepsilon(\rho+h_1)\partial_r(\Delta r\circ\overline{X})|_{(0,s,t)}+\varepsilon\partial_\rho\hat{u}_1^I\Delta r|_{\overline{X}_0(s,t)}\right]+\\
+\frac{1}{\varepsilon}2\nabla r\cdot\nabla s|_{\overline{X}_0(s,t)}\left[\partial_s\partial_\rho\hat{u}_1^I\varepsilon-\partial_sh_1\varepsilon\partial_\rho^2\hat{u}_1^I-\varepsilon\partial_sh_2\theta_0''(\rho)\right]+\\
+\frac{1}{\varepsilon}2\partial_r((\nabla r\cdot\nabla s)\circ\overline{X})|_{(0,s,t)}\varepsilon(\rho+h_1)\left[\partial_s\theta_0'(\rho)-\partial_sh_1\theta_0''(\rho)\right]+\\
+\Delta s|_{\overline{X}_0(s,t)}\partial_s\theta_0(\rho)+|\nabla s|^2|_{\overline{X}_0(s,t)}\partial_s^2\theta_0(\rho)-2|\nabla s|^2|_{\overline{X}_0(s,t)}\partial_sh_1\partial_s\theta_0'(\rho)+\\
-\theta_0'(\rho)\left[\Delta s|_{\overline{X}_0(s,t)}\partial_sh_1+|\nabla s|^2|_{\overline{X}_0(s,t)}\partial_s^2h_1\right]+|\nabla s|^2|_{\overline{X}_0(s,t)}(\partial_sh_1)^2\theta_0''(\rho)=\\
=\theta_0''(\rho)\left[\frac{1}{2}(\rho+h_1)^2\partial_r^2(|\nabla r|^2\circ\overline{X})|_{(0,s,t)}+(\partial_sh_1)^2|\nabla s|^2|_{\overline{X}_0(s,t)}\right]+\partial_\rho^2\hat{u}_2^I+\\
+\theta_0''(\rho)2\partial_r((\nabla r\cdot\nabla s)\circ\overline{X})|_{(0,s,t)}(\rho+h_1)(-\partial_sh_1)+\\
-\theta_0'(\rho)\left[\Delta s|_{\overline{X}_0(s,t)}\partial_sh_1+|\nabla s|^2|_{\overline{X}_0(s,t)}\partial_s^2h_1\right]+\theta_0'(\rho)(\rho+h_1)\partial_r(\Delta r\circ\overline{X})|_{(0,s,t)}.
\end{align*}
Since $\hat{u}_1^I=0$, the contribution from the $f'$-part is $f''(\theta_0)\hat{u}_2^I$. Therefore for the cancellation of the $\Oc(1)$-term in the expansion for the Allen-Cahn equation we require
\begin{align*}
-\Lc\hat{u}_2^I(\rho,s,t)=R_1(\rho,& s,t):=\\
\theta_0'(\rho)&\left[(\rho+h_1)\partial_r((\partial_tr-\Delta r)\circ\overline{X})|_{(0,s,t)}-\partial_th_1+|\nabla s|^2|_{\overline{X}_0(s,t)}\partial_s^2h_1+\right.\\
&\left.-(\partial_ts-\Delta s)|_{\overline{X}_0(s,t)}\partial_sh_1\right]+\\
+\theta_0''(\rho)&\left[-\frac{1}{2}(\rho+h_1)^2\partial_r^2(|\nabla r|^2\circ\overline{X})|_{(0,s,t)}+\right.\\
&\left.+2(\rho+h_1)\partial_sh_1\partial_r((\nabla r\cdot\nabla s)\circ\overline{X})|_{(0,s,t)}-|\nabla s|^2|_{\overline{X}_0(s,t)}(\partial_sh_1)^2\right].
\end{align*}
Here because of Remark \ref{th_coord_bem}, 4. it holds $|\nabla s|^2|_{\overline{X}_0(s,t)}=4/L(t)^2$. If $h_1$ is smooth, then $R_1$ is smooth and together with all derivatives decays exponentially in $|\rho|$ uniformly in $(s,t)$. Hence Lemma \ref{th_ODE_lin2} yields that there is a (unique) solution $\hat{u}_2^I$ together with $\hat{u}_2^I(0,s,t)=0$ and boundedness in $\rho$ if and only if $\int_\R R_1(\rho,s,t)\theta_0'(\rho)\,d\rho=0$. Since $\int_\R\theta_0'(\rho)\theta_0''(\rho)\,d\rho=0$ due to integration by parts, the nonlinearities in $h_1$ drop out and we obtain a linear non-autonomous parabolic equation for $h_1$ with principal part $\partial_t-\frac{4}{L(t)^2}\partial_s^2$:
\begin{align}\label{eq_h1}
\partial_th_1-\frac{4}{L(t)^2}\partial_s^2h_1+b_1\partial_sh_1+c_1h_1=f_1\quad\text{ on }I\times[0,T].
\end{align} 
Here with
\begin{align}\begin{split}\label{eq_dj}
d_1:=\int_\R\theta_0'(\rho)^2\,d\rho,\quad d_2:=\int_\R\theta_0'(\rho)^2\rho\,d\rho,\quad
d_3:=\int_\R\theta_0'(\rho)^2\rho^2\,d\rho,\\
d_4:=\int_\R\theta_0'(\rho)\theta_0''(\rho)\rho\,d\rho,\quad d_5:=\int_\R\theta_0'(\rho)\theta_0''(\rho)\rho^2\,d\rho,\quad
d_6:=\int_\R\theta_0'(\rho)\theta_0''(\rho)\rho^3\,d\rho,\end{split}
\end{align}
we have set for all $(s,t)\in I\times[0,T]$
\begin{align}
b_1(s,t)&:=(\partial_ts-\Delta s)|_{\overline{X}_0(s,t)}-2\frac{d_4}{d_1}\partial_r((\nabla r\cdot\nabla s)\circ\overline{X})|_{(0,s,t)},\label{eq_b1}\\
c_1(s,t)&:=-\partial_r((\partial_tr-\Delta r)\circ\overline{X})|_{(0,s,t)}+\frac{d_4}{d_1}\partial_r^2(|\nabla r|^2\circ\overline{X})|_{(0,s,t)},\label{eq_c1}\\
f_1(s,t)&:=\frac{d_2}{d_1}\partial_r((\partial_tr-\Delta r)\circ\overline{X})|_{(0,s,t)}-\frac{d_5}{2d_1}\partial_r^2(|\nabla r|^2\circ\overline{X})|_{(0,s,t)}.
\end{align}
If $h_1$ is smooth and solves \eqref{eq_h1}, then $\hat{u}_2^I$ exists, is smooth and $\hat{u}_2^I$ together with all derivatives decays exponentially in $|\rho|$ uniformly in $(s,t)$ because of Lemma \ref{th_ODE_lin2}.
\begin{Remark}\label{th_bem_h1_feven1}\upshape
If $f$ is additionally even, then $\theta_0'$ is even and hence $d_2=d_5=0$. Thus $f_1=0$ and \eqref{eq_h1} is homogeneous in this case.
\end{Remark}
\subsubsection{Inner Expansion: $\Oc(\varepsilon)$}\label{sec_inner_eps}
At order $\varepsilon$ we get from $\partial_tu$:
\begin{align*}
\frac{1}{\varepsilon}\left[\varepsilon^2\partial_\rho\hat{u}_2^I\partial_tr|_{\overline{X}_0(s,t)}+\theta_0'(\rho)\varepsilon^2(h_2\partial_r(\partial_tr\circ\overline{X})|_{(0,s,t)}+\frac{(\rho+h_1)^2}{2}\partial_r^2(\partial_tr\circ\overline{X})|_{(0,s,t)})\right]+\\
-\varepsilon\theta_0'(\rho)\left[\partial_th_2+\partial_ts|_{\overline{X}_0(s,t)}\partial_sh_2+(\rho+h_1)\partial_sh_1\partial_r(\partial_ts\circ\overline{X})|_{(0,s,t)}\right],
\end{align*}
where we used $\hat{u}_1^I=0$ and $\partial_s\theta_0=0$. The contribution of $\Delta u$ is
\begin{align*}
\varepsilon\partial_\rho^2\hat{u}_3^I+\varepsilon\theta_0''(\rho)\left[0+\partial_r^2(|\nabla r|^2\circ\overline{X})|_{(0,s,t)}(\rho+h_1)h_2+\frac{\partial_r^3(|\nabla r|^2\circ\overline{X})|_{(0,s,t)}}{3!}(\rho+h_1)^3\right]+\\
+\frac{1}{\varepsilon}\theta_0'(\rho)\varepsilon^2\left[ h_2\partial_r(\Delta r\circ\overline{X})|_{(0,s,t)}+\frac{(\rho+h_1)^2}{2}\partial_r^2(\Delta r\circ\overline{X})|_{(0,s,t)}\right]+0+\frac{1}{\varepsilon}\varepsilon^2\partial_\rho\hat{u}_2^I\Delta r|_{\overline{X}_0(s,t)}+\\
+\frac{2}{\varepsilon}\left[0+\varepsilon(\rho+h_1)\partial_r((\nabla r\cdot\nabla s)\circ\overline{X})|_{(0,s,t)}(\varepsilon\partial_s\partial_\rho\hat{u}_1^I-\varepsilon\theta_0''(\rho)\partial_sh_2-\varepsilon\partial_\rho^2\hat{u}_1^I\partial_sh_1)\right]+\\
-\frac{2}{\varepsilon}\theta_0''(\rho)\partial_sh_1\left[\varepsilon^2h_2\partial_r((\nabla r\cdot\nabla s)\circ\overline{X})|_{(0,s,t)}+\varepsilon^2\frac{(\rho+h_1)^2}{2}\partial_r^2((\nabla r\cdot\nabla s)\circ\overline{X})|_{(0,s,t)}\right]+\\
+\varepsilon(\rho+h_1)\partial_r(\Delta s\circ\overline{X})|_{(0,s,t)}\partial_s\theta_0+\Delta s|_{\overline{X}_0(s,t)}\partial_s\hat{u}_1^I+0+\\
-\theta_0'(\rho)\left[\Delta s|_{\overline{X}_0(s,t)}\varepsilon\partial_sh_2+\varepsilon(\rho+h_1)\partial_sh_1\partial_r(\Delta s\circ\overline{X})|_{(0,s,t)}\right]+\\
-\theta_0'(\rho)\left[|\nabla s|^2|_{\overline{X}_0(s,t)}\varepsilon\partial_s^2h_2+\varepsilon(\rho+h_1)\partial_s^2h_1\partial_r(|\nabla s|^2\circ\overline{X})|_{(0,s,t)}\right]+0+\\
+\theta_0''(\rho)\left[\varepsilon(\rho+h_1)\partial_r(|\nabla s|^2\circ\overline{X})|_{(0,s,t)}(\partial_sh_1)^2+2|\nabla s|^2|_{\overline{X}_0(s,t)}\varepsilon\partial_sh_1\partial_sh_2\right]+0.
\end{align*}
Since $u_1^I=0$, we have $f''(\theta_0)\hat{u}_3^I$ from the $f'$-part. Altogether the $\varepsilon$-order cancels if
\begin{align*}
-\Lc\hat{u}_3^I(\rho,s,t)=R_2(\rho,s,t):=\\
\theta_0'(\rho)\left[h_2\partial_r((\partial_tr-\Delta r)\circ\overline{X})|_{(0,s,t)}+\frac{(\rho+h_1)^2}{2}\partial_r^2((\partial_tr-\Delta r)\circ\overline{X})|_{(0,s,t)}-\partial_th_2+\right.\\
-(\partial_ts-\Delta s)|_{\overline{X}_0(s,t)}\partial_sh_2-(\rho+h_1)\partial_sh_1\partial_r((\partial_ts-\Delta s)\circ\overline{X})|_{(0,s,t)}+\\
\left.+\frac{4}{L(t)^2}\partial_s^2h_2+(\rho+h_1)\partial_s^2h_1\partial_r(|\nabla s|^2\circ\overline{X})|_{(0,s,t)}\right]+\\
+\theta_0''(\rho)\left[-(\rho+h_1)h_2\partial_r^2(|\nabla r|^2\circ\overline{X})|_{(0,s,t)}-\frac{(\rho+h_1)^3}{3!}\partial_r^3(|\nabla r|^2\circ\overline{X})|_{(0,s,t)}-\frac{8\partial_sh_1\partial_sh_2}{L(t)^2}+\right.\\
+2(\rho+h_1)\partial_sh_2\partial_r((\nabla r\cdot\nabla s)\circ\overline{X})|_{(0,s,t)}+2\partial_sh_1(h_2\partial_r((\nabla r\cdot\nabla s)\circ\overline{X})|_{(0,s,t)}+\\
\left.-\frac{(\rho+h_1)^2}{2}\partial_r^2((\nabla r\cdot\nabla s)\circ\overline{X})|_{(0,s,t)})-(\rho+h_1)(\partial_sh_1)^2\partial_r(|\nabla s|^2\circ\overline{X})|_{(0,s,t)}\right].
\end{align*}
If $h_1,h_2$ are smooth, then $R_2$ is smooth and together with all derivatives decays exponentially in $|\rho|$ uniformly in $(s,t)$. Therefore Lemma \ref{th_ODE_lin2} implies the existence of a unique smooth solution $\hat{u}_3^I$ that fulfills $\hat{u}_3^I(0,s,t)=0$ and is bounded in $\rho$ if and only if $\int_\R R_2(\rho,s,t)\theta_0'(\rho)\,d\rho=0$. This yields a linear non-autonomous parabolic equation for $h_2$ with principal part $\partial_t-\frac{4}{L(t)^2}\partial_s^2$:
\begin{align}\label{eq_h2}
\partial_th_2-\frac{4}{L(t)^2}\partial_s^2h_2+b_1\partial_sh_2+c_1h_2=f_2\quad\text{ on }I\times[0,T],
\end{align}
where $b_1$ and $c_1$ are defined in \eqref{eq_b1} and \eqref{eq_c1}, respectively, and with the $d_j, j=1,...,6$ as in \eqref{eq_dj} we set for $(s,t)\in I\times[0,T]$
\begin{align*}
f_2(s,t)&:=\frac{1}{2d_1}\left[d_3+2d_2h_1+h_1^2\right]\partial_r^2((\partial_tr-\Delta r)\circ\overline{X})|_{(0,s,t)}+\\
&+(\frac{d_2}{d_1}+h_1)\left[-\partial_sh_1\partial_r((\partial_ts-\Delta s)\circ\overline{X})|_{(0,s,t)}+\partial_s^2h_1\partial_r(|\nabla s|^2\circ\overline{X})|_{(0,s,t)}\right]+\\
&-\frac{1}{3!d_1}[d_6+2d_5h_1+2d_4h_1^2+d_3h_1^3]\partial_r^3(|\nabla r|^2\circ\overline{X})|_{(0,s,t)}+\\
&-\frac{1}{d_1}[d_5+2d_4h_1+d_3h_1^2]\frac{\partial_r^2((\nabla r\cdot\nabla s)\circ\overline{X})|_{(0,s,t)}}{2}-\frac{d_4}{d_1}(\partial_sh_1)^2\partial_r(|\nabla s|^2\circ\overline{X})|_{(0,s,t)}.
\end{align*}
If $h_2$ is smooth and satisfies \eqref{eq_h2}, then we obtain from Lemma \ref{th_ODE_lin2} that $\hat{u}_3^I$ exists, is smooth and $\hat{u}_3^I$ together with all derivatives decays exponentially in $|\rho|$ uniformly in $(s,t)$.
\subsection{Boundary Layer Expansion}\label{sec_bdry_exp}
In the boundary layer expansion we make the ansatz $u=u^I+u^{B\pm}$ in $\Gamma(2\delta)$ near the contact points $p^\pm(t), t\in[0,T]$. For $u^{B\pm}$ we combine the stretched-variable ansatz with stretching the tangential variable: With $s^\pm:=\mp(s\mp 1)$ and $H^\pm_\varepsilon:=\frac{s^\pm}{\varepsilon}$ we set
\[
u^{B\pm}=u_1^{B\pm}\varepsilon+u_2^{B\pm}\varepsilon^2,\quad u_j^{B\pm}(x,t)=\hat{u}_j^{B\pm}(\rho_\varepsilon(x,t),H^\pm_\varepsilon(x,t),t)\quad\text{ for }(x,t)\in\Gamma(2\delta)
\]
and $\hat{u}^{B\pm}:=\hat{u}_1^{B\pm}\varepsilon+\hat{u}_2^{B\pm}\varepsilon^2$, where $\R^2_+\times[0,T]\ni(\rho,H,t)\mapsto \hat{u}_j^{B\pm}(\rho,H,t)\in\R$ for $j=1,2$.\\
To simplify the asymptotic expansion, we remark that later $u^I$ should solve the equation $\partial_tu^I-\Delta u^I+f'(u^I)/\varepsilon^2=0$ approximately. Therefore instead of the Allen-Cahn equation for $u=u^I+u^{B\pm}$, we expand the difference 
\begin{align}\label{eq_bdry_layer_equ}
\partial_tu^{B\pm}-\Delta u^{B\pm}+\frac{1}{\varepsilon^2}\left[f'(u^I+u^{B\pm})-f'(u^I)\right]=0
\end{align}
into $\varepsilon$-series with coefficients in $(\rho_\varepsilon,H^\pm_\varepsilon,t)$. Here compared to the inner expansion we only expand up to $\Oc(1)$ which later turns out to be sufficient. Moreover, we expand the Neumann boundary condition into $\varepsilon$-series with coefficients in $(\rho_\varepsilon,t)$ up to $\Oc(\varepsilon)$. The cancellation of terms in the expansions yield equations on $\R^2_+$ of type as in Subsection \ref{sec_hp} (up to a $t$-dependent scaling in $H$). The corresponding solvability conditions will give us the boundary conditions for the height functions.\\
In the following lemma we compute how the differential operators act on $u^{B\pm}$.
\begin{Lemma}\label{th_trafo_bdry}
Let $\R^2_+\times[0,T]\ni(\rho,H,t)\mapsto\hat{w}(\rho,H,t)\in\R$ be sufficiently smooth and $w:\Gamma(2\delta)\rightarrow\R$ defined by $w(x,t):=\hat{w}(\rho_\varepsilon(x,t),\frac{s^\pm(x,t)}{\varepsilon},t)$ for all $(x,t)\in\Gamma(2\delta)$. Then
\begin{align*}
\partial_tw&=\partial_\rho\hat{w}\left[\frac{\partial_tr}{\varepsilon}-(\partial_th_\varepsilon+\partial_ts\partial_sh_\varepsilon)\right]\mp\frac{1}{\varepsilon}\partial_H\hat{w}\,\partial_ts+\partial_t\hat{w},\\
\nabla w&=\partial_\rho\hat{w}\left[\frac{\nabla r}{\varepsilon}-\nabla s\partial_sh_\varepsilon\right]\mp\frac{1}{\varepsilon}\partial_H\hat{w}\,\nabla s,\\
\Delta w&=\partial_\rho\hat{w}\left[\frac{\Delta r}{\varepsilon}-(\Delta s\partial_sh_\varepsilon+|\nabla s|^2\partial_s^2h_\varepsilon)\right]\mp\frac{1}{\varepsilon}\partial_H\hat{w}\,\Delta s+\frac{1}{\varepsilon^2}\partial_H^2\hat{w}\,|\nabla s|^2+\\
&\mp 2\partial_\rho\partial_H\hat{w}\,\frac{\nabla s}{\varepsilon}\cdot\left[\frac{\nabla r}{\varepsilon}-\nabla s\partial_sh_\varepsilon\right]+\partial_\rho^2\hat{w}\left|\frac{\nabla r}{\varepsilon}-\nabla s\partial_sh_\varepsilon\right|^2,
\end{align*}
where the $w$-terms on the left hand side and derivatives of $r$ or $s$ are evaluated at $(x,t)$, the $h_\varepsilon$-terms at $(s(x,t),t)$ and the $\hat{w}$-terms at $(\rho_\varepsilon(x,t),\frac{s^\pm(x,t)}{\varepsilon},t)$.
\end{Lemma}
\textit{Proof.} This can be directly shown using the chain rule.\hfill$\square$\\
\newline
We have to expand the $f'$-part: If the $u_j^I, u_j^{B\pm}$ are uniformly bounded, we apply a Taylor expansion to obtain on $\Gamma(2\delta)$ with $u_3^{B\pm}:=0$
\begin{align}\label{eq_bdry_taylor_f}
f'(u^I+u^{B\pm})=f'(\theta_0)+\sum_{k=1}^m\frac{1}{k!}f^{(k+1)}(\theta_0)\left[\sum_{j=1}^3(u_j^I+u_j^{B\pm})\varepsilon^j\right]^k+\Oc(\varepsilon^{m+1}).
\end{align}
Combining this with the expansion for $f'(u^I)$ in \eqref{eq_inner_taylor_f}, the first few terms in the asymptotic expansion for $f'(u^I+u^{B\pm})-f'(u^I)$ are
\begin{align*}
\Oc(1)&:\quad 0,\\
\Oc(\varepsilon)&:\quad f''(\theta_0)u_1^{B\pm},\\
\Oc(\varepsilon^2)&:\quad f''(\theta_0)u_2^{B\pm}+\frac{f^{(3)}(\theta_0)}{2!}(u_1^{B\pm})^2,\\
\Oc(\varepsilon^3)&:\quad f''(\theta_0)u_3^{B\pm}+2\frac{f^{(3)}(\theta_0)}{2!}u_1^{B\pm}(u_2^I+u_2^{B\pm})+\frac{f^{(4)}(\theta_0)}{3!}(u_1^{B\pm})^3.
\end{align*}
\begin{Remark}\upshape\label{th_bdry_taylor_f_bem}
Also in higher orders every term contains a boundary term as factor.
\end{Remark}
To derive an $\varepsilon$-expansion with coefficients in $(\rho,H,t)$, we have to expand terms depending on $(s,t)$, $(\rho,s,t)$ and $(x,t)$. Let $g:I\times[0,T]\rightarrow\R$ or $g:\R\times I\times[0,T]\rightarrow\R$ smooth with uniformly bounded derivatives in $s$. Since $s=\pm1\mp\varepsilon H^\pm_\varepsilon$, we apply a Taylor expansion to a smooth extension to get for $m\in\N$
\begin{align}\label{eq_bdry_layer_taylor2}
g|_{s=\pm1\mp\varepsilon H}=g|_{s=\pm1}+\sum_{j=1}^m(\mp\varepsilon H)^j\frac{\partial_s^j g|_{s=\pm1}}{j!}+\Oc((\varepsilon H)^{m+1})\quad\text{ for }H\in[0,\varepsilon].
\end{align}
Moreover, if $g:\overline{\Gamma(2\delta)}\rightarrow\R$ is smooth, then a Taylor expansion yields for $m\in\N$
\[
\tilde{g}(r,s,t):=g(\overline{X}(r,s,t))=\sum_{k+l=0}^m\frac{\partial_r^k\partial_s^l\tilde{g}(0,\pm1,t)}{k!\,l!}r^k(s\mp 1)^l+\Oc(|(r,s\mp 1)|^{m+1})
\]
uniformly in $(r,s,t)\in[-2\delta,2\delta]\times I\times[0,T]$. Later we insert 
\[
r=\varepsilon(\rho_\varepsilon(x,t)+h_\varepsilon(s,t)),\quad s=\pm1\mp\varepsilon H^\pm_\varepsilon(x,t)\quad\text{ for }(x,t)\in\Gamma(2\delta)
\] 
and expand $h_\varepsilon$ with \eqref{eq_bdry_layer_taylor2}. Then $(\rho_\varepsilon,H^\pm_\varepsilon)$ are replaced by arbitrary $(\rho,H)\in\R^2_+$ in order to derive suitable equations. The first terms needed in the resulting expansion are
\begin{align*}
\Oc(1):\quad & g|_{\overline{p}^\pm(t)},\\
\Oc(\varepsilon):\quad &\partial_r\tilde{g}|_{(0,\pm 1,t)}(\rho+h_1|_{(\pm1,t)})+\partial_s\tilde{g}|_{(0,\pm 1,t)}(\mp H).
\end{align*} 
The remainder term without $\Oc(|(r,s\mp 1)|^{m+1})$ can be estimated by $\varepsilon^2$ times some polynomial in $(|\rho_\varepsilon|,H^\pm_\varepsilon)$.
\subsubsection{Boundary Layer Expansion: Equation \eqref{eq_bdry_layer_equ}}\label{sec_bdry_exp_equ}
We expand \eqref{eq_bdry_layer_equ} with the above identities. Let $\Delta_t:=\partial_\rho^2+\frac{4}{L(t)^2}\partial_H^2$. Because of Remark \ref{th_coord_bem} and Theorem \ref{th_coord}, the lowest order $\Oc(\frac{1}{\varepsilon})$ vanishes if 
\begin{align}\label{eq_bdry_eq1}
-\Delta_t\hat{u}_1^{B\pm}(\rho,H,t)+f''(\theta_0(\rho))\hat{u}_1^{B\pm}(\rho,H,t)+0\cdot\partial_H\partial_\rho u_1^{B\pm}(\rho,H,t)=0.
\end{align}
Next, we consider the $\Oc(1)$-order. From $\partial_tu$ we get $\partial_tr|_{\overline{p}^\pm(t)}\partial_\rho \hat{u}_1^{B\pm}\mp\partial_H\hat{u}_1^{B\pm}\partial_ts|_{(\pm 1,t)}$ and from $\Delta u$
\begin{align*}
0\cdot\partial_\rho^2\hat{u}_1^{B\pm}+\Delta_t\hat{u}_2^{B\pm}+\left[\partial_r(|\nabla s|^2\circ\overline{X})|_{(0,\pm 1,t)}(\rho+h_1|_{(\pm 1,t)})+0\right]\partial_H^2 \hat{u}_1^{B\pm}+\\
\mp 0\cdot \partial_H\partial_\rho \hat{u}_2^{B\pm}\mp 2\left[ \partial_r((\nabla r\cdot\nabla s)\circ\overline{X})|_{(0,\pm 1,t)}(\rho+h_1|_{(\pm 1,t)})+0\right]\partial_H\partial_\rho \hat{u}_1^{B\pm}+\\
+\Delta r|_{\overline{p}^\pm(t)}\partial_\rho\hat{u}_1^{B\pm}+0\pm\frac{8}{L(t)^2}\partial_sh_1|_{(\pm1,t)}\partial_H\partial_\rho\hat{u}_1^{B\pm}\mp\Delta s|_{\overline{p}^\pm(t)}\partial_H\hat{u}_1^{B\pm}.
\end{align*}
Since $(\partial_tr-\Delta r)\circ\overline{X}_0=0$, the $\Oc(1)$-order in the expansion for \eqref{eq_bdry_layer_equ} is zero if 
\begin{align}\begin{split}\label{eq_bdry_eq2}
(-\Delta_t+f''(\theta_0))\hat{u}_2^{B\pm}=G_2(\rho,H,t):=-\frac{f^{(3)}(\theta_0)}{3!}(\hat{u}_1^{B\pm})^2+\\
\pm\partial_H\hat{u}_1^{B\pm}(\partial_ts-\Delta s)|_{\overline{p}^\pm(t)}+(\rho+h_1|_{(\pm 1,t)})\partial_H^2\hat{u}_1^{B\pm}\partial_r(|\nabla s|^2\circ\overline{X})|_{(0,\pm1,t)}+\\
\pm 2\partial_H\partial_\rho\hat{u}_1^{B\pm}\left[\frac{4}{L(t)^2}\partial_sh_1|_{(\pm 1,t)}-(\rho+h_1|_{(\pm 1,t)})\partial_r((\nabla r\cdot\nabla s)\circ\overline{X})|_{(0,\pm 1,t)}\right].\end{split}
\end{align}
To complete these equations we need boundary conditions. These will be obtained from
\subsubsection{Boundary Layer Expansion: The Neumann Boundary Condition}\label{sec_bdry_exp_neum}
Next we consider the Neumann boundary condition, i.e. $N_{\partial\Omega}\cdot\nabla(u^I+u^{B\pm})|_{\partial\Omega}=0$. Lemma \ref{th_trafo_inner} and Lemma \ref{th_trafo_bdry} yield on $\Gamma(2\delta)$
\begin{align*}
\nabla u^I|_{(x,t)}&=\partial_\rho\hat{u}^I|_{(\rho,s,t)}\left[\frac{\nabla r|_{(x,t)}}{\varepsilon}-\nabla s|_{(x,t)}\partial_sh_\varepsilon|_{(s,t)}\right]+\partial_s\hat{u}^I|_{(\rho,s,t)}\,\nabla s|_{(x,t)},\\
\nabla u^{B\pm}|_{(x,t)}&=\partial_\rho\hat{u}^{B\pm}|_{(\rho,H,t)}\left[\frac{\nabla r|_{(x,t)}}{\varepsilon}-\nabla s|_{(x,t)}\partial_sh_\varepsilon|_{(s,t)}\right]\mp\frac{1}{\varepsilon}\partial_H\hat{u}^{B\pm}|_{(\rho,H,t)}\,\nabla s|_{(x,t)},
\end{align*}
where $\rho=\rho_\varepsilon(x,t), H=H^\pm_\varepsilon(x,t)$ and $s=s(x,t)$. We evaluate at $x=X(r,\pm 1,t)$, i.e. $H=0$ and $s=\pm 1$. For $g:\overline{\Gamma(2\delta)}\cap\partial\Omega\rightarrow\R$ smooth we use the analogous Taylor expansion \eqref{eq_inner_taylor} for $s=\pm 1$. As before we use $r=\varepsilon(\rho_\varepsilon+h_\varepsilon)$ and replace $\rho_\varepsilon$ by an arbitrary $\rho\in\R$.\\
Using this we expand the Neumann boundary condition. At lowest order $\Oc(\frac{1}{\varepsilon})$ we have
$(N_{\partial\Omega}\cdot\nabla r)|_{\overline{p}^\pm(t)}\theta_0'(\rho)$. This vanishes because we required a $90$°-contact angle.\\
At $\Oc(1)$ we obtain
\begin{align*}
& 0\cdot\partial_\rho\hat{u}_1^{B\pm}|_{H=0}+0\cdot\partial_\rho\hat{u}_1^I|_{s=\pm 1}+(N_{\partial\Omega}\cdot\nabla s)|_{\overline{p}^\pm(t)}\left[\mp\partial_H\hat{u}_1^{B\pm}|_{H=0}+\partial_s\theta_0(\rho)\right]+\\
&+\theta_0'(\rho)\left[-(N_{\partial\Omega}\cdot\nabla s)|_{\overline{p}^\pm(t)}\partial_sh_1|_{(\pm 1,t)}+\partial_r((N_{\partial\Omega}\cdot\nabla r)\circ\overline{X})|_{(0,\pm 1,t)}(\rho+h_1|_{(\pm 1,t)})\right].
\end{align*}
Here because of Remark \ref{th_coord_bem}, 4. it holds $(N_{\partial\Omega}\cdot\nabla s)|_{\overline{p}^\pm(t)}=\pm 1/|\partial_sX_0(s,t)|=\pm 2/L(t)$. Therefore the $\Oc(1)$-order cancels if 
\begin{align}\begin{split}\label{eq_bdry_rb1}
-\frac{2}{L(t)}\partial_H\hat{u}_1^{B\pm}(\rho,0,t)=g_1^\pm(\rho,t),\\
g_1^\pm(\rho,t):=\theta_0'(\rho)\left[\pm\frac{2}{L(t)}\partial_sh_1|_{(\pm 1,t)}-\partial_r((N_{\partial\Omega}\cdot\nabla r)\circ\overline{X})|_{(0,\pm 1,t)}(\rho+h_1|_{(\pm 1,t)})\right].
\end{split}
\end{align}
We define $\overline{u}_j^{B\pm}:\R^2_+\times[0,T]\rightarrow\R:(\rho,H,t)\mapsto\hat{u}_j^{B\pm}(\rho,\frac{2H}{L(t)},t)$ for $j=1,2$. Then equations \eqref{eq_bdry_eq1} and \eqref{eq_bdry_rb1} for $\hat{u}_1^{B\pm}$ are equivalent to
\begin{align}\label{eq_u1B_1}
-\Delta\overline{u}_1^{B\pm}+f''(\theta_0)\overline{u}_1^{B\pm}&=0,\\
-\partial_{H}\overline{u}_1^{B\pm}|_{H=0}&=g_1^\pm(\rho,t).\label{eq_u1B_2}
\end{align}
The corresponding solvability condition \eqref{eq_hp_intbed} is
$\int_\R g_1^\pm(\rho,t)\theta_0'(\rho)\,d\rho=0$.
This gives a linear boundary condition for $h_1$ of the form
\[
\pm\frac{2}{L(t)}\partial_sh_1|_{(\pm 1,t)}-\partial_r((N_{\partial\Omega}\cdot\nabla r)\circ\overline{X})|_{(0,\pm 1,t)}h_1|_{(\pm 1,t)}=\tilde{g}_1^\pm(t),
\]
where $\tilde{g}_1^\pm$ is smooth. Together with the parabolic equation for $h_1$ from Subsection \ref{sec_inner_1}, we have a time-dependent parabolic boundary value problem, where $h_1|_{t=0}$ is not specified yet. To obtain a smooth solution, certain compatibility conditions have to be fulfilled, cf. Lunardi, Sinestrari and von Wahl \cite{LunardiSvW}, Chapter 9. To solve this problem, we extend the coefficients and right hand sides smoothly to $[-T,T]$ such that the coefficient in front of $\partial_s^2$ is bounded below by a $c_0>0$ and the right hand sides are zero for $t\leq-\frac{1}{2}T$. Then for the initial value zero at $t=-T$ all compatibility conditions are fulfilled and we obtain a smooth solution on $[-T,T]$ by \cite{LunardiSvW}, Theorem 9.1. Restriction to $[0,T]$ yields a smooth solution $h_1$ on $[0,T]$.
\begin{Remark}\label{th_bem_h1_feven2}\upshape
If $f$ is even, then so is $\theta_0'$ and the boundary condition for $h_1$ is homogeneous. Hence because of Remark \ref{th_bem_h1_feven1} we can take $h_1=0$ in this case.
\end{Remark}
Having determined $h_1$, we know that $g_1^\pm$ is smooth and decays together with all derivatives exponentially in $|\rho|$ uniformly in $t\in[0,T]$. Corollary \ref{th_hp_time_cor} yields a unique smooth solution $\overline{u}_1^{B\pm}$ to \eqref{eq_u1B_1}-\eqref{eq_u1B_2} such that $\overline{u}_1^{B\pm}$ and all derivatives decay exponentially in $|\rho|+H$ uniformly in $t\in[0,T]$.\\
For the $\varepsilon$-order we get
\begin{align*}
\frac{1}{\varepsilon}\left[0\cdot\varepsilon^2\partial_\rho\hat{u}_2^{B\pm}|_{H=0}+\varepsilon\partial_\rho\hat{u}_1^{B\pm}|_{H=0}\varepsilon(\rho+h_1|_{(\pm 1,t)})\partial_r((N_{\partial\Omega}\cdot\nabla r)\circ\overline{X})|_{(0,\pm 1,t)}\right]+\\
\mp\frac{1}{\varepsilon}\left[\pm\varepsilon^2\frac{2}{L(t)}\partial_H\hat{u}_2^{B\pm}|_{H=0}+\varepsilon\partial_H\hat{u}_1^{B\pm}|_{H=0}\varepsilon(\rho+h_1|_{(\pm 1,t)})\partial_r((N_{\partial\Omega}\cdot\nabla s))|_{(0,\pm 1,t)}\right]+\\
\mp\frac{2}{L(t)}\varepsilon\partial_\rho\hat{u}_1^{B\pm}|_{H=0}\partial_sh_1|_{(\pm 1,t)}+0+0+\\
+\varepsilon\theta_0'\left[\frac{(\rho+h_1|_{(\pm 1,t)})^2}{2}\partial_r^2((N_{\partial\Omega}\cdot\nabla r)\circ\overline{X})|_{(0,\pm 1,t)}+h_2|_{(\pm 1,t)}\partial_r((N_{\partial\Omega}\cdot\nabla r)\circ\overline{X})|_{(0,\pm 1,t)}\right]+\\
+0-\theta_0'(\rho)\left[\pm\frac{2}{L(t)}\varepsilon\partial_sh_2|_{(\pm 1,t)}+\varepsilon(\rho+h_1|_{(\pm 1,t)})\partial_sh_1|_{(\pm 1,t)}\partial_r((N_{\partial\Omega}\cdot\nabla s)\circ\overline{X})|_{(0,\pm 1,t)}\right].
\end{align*}
The requirement that this term vanishes leads to an equation of the type 
\[
-\frac{2}{L(t)}\partial_H\hat{u}_2^{B\pm}|_{H=0}=g_2^\pm(\rho,t)\quad\text{ for }(\rho,t)\in\R\times[0,T],
\] 
where $g_2^\pm$ is smooth and decays together with all derivatives exponentially in $|\rho|$ uniformly in $t$. Together with \eqref{eq_bdry_eq2} and $\overline{G}_2:\R^2\times[0,T]\rightarrow\R:(\rho,H,t)\mapsto G_2(\rho,\frac{2H}{L(t)},t)$ this yields
\begin{align*}
(-\Delta+f''(\theta_0))\overline{u}_2^{B\pm}&=\overline{G}_2,\\
-\partial_H\overline{u}_2^{B\pm}|_{H=0}&=g_2^\pm.
\end{align*}
The corresponding compatibility condition \eqref{eq_hp_intbed}
\[
\int_{\R^2_+}\overline{G}_2(\rho,H,t)\theta_0'(\rho)\,d(\rho,H)+\int_\R g_2^\pm(\rho,t)\,d\rho=0
\]
leads to a boundary condition for $h_2$ of the form
\[
\pm\frac{2}{L(t)}\partial_sh_2|_{(\pm 1,t)}-\partial_r((N_{\partial\Omega}\cdot\nabla r)\circ\overline{X})|_{(0,\pm 1,t)}h_2|_{(\pm 1,t)}=\tilde{g}_2^\pm(t),
\]
where $\tilde{g}_2^\pm$ is smooth. Similarly as above we solve this equation together with the parabolic one in Section \ref{sec_inner_eps} and get a smooth solution $h_2$. Since $G_2, g_2^\pm$ are smooth and (together with all derivatives) decay exponentially in $|\rho|+H$, $|\rho|$, respectively, uniformly for $t\in[0,T]$, we obtain from Corollary \ref{th_hp_time_cor} a unique smooth solution $\overline{u}_2^{B\pm}$ with analogous decay.
\subsection{The Approximate Solution}\label{sec_approx_sol}
Let $\Gamma$ be as before and a smooth solution to mean curvature flow (with $90$°-contact angle). Let $\chi\in C^\infty(\R)$ with $\chi(r)=1$ for $|r|\leq 1$ and $\chi(r)=0$ for $|r|\geq 2$. Then we set for $\varepsilon>0$
\begin{align*}
u_A^\varepsilon:=
\begin{cases}
\chi(\frac{r}{\delta})\left[u^I+\sum_\pm u^{B\pm} \chi(\frac{s^\pm}{\delta_0})\right]+(1-\chi(\frac{r}{\delta}))\textup{sign}(r)&\quad\text{ on }\overline{\Gamma(2\delta)},\\
\pm 1&\quad\text{ on }Q_T^\pm\textbackslash\Gamma(2\delta),
\end{cases}
\end{align*}
where $u^I$ and $u^{B\pm}$ were constructed in Sections \ref{sec_inner_exp} and \ref{sec_bdry_exp} and depend on $\varepsilon>0$. This is an approximate solution in the following sense:
\begin{Lemma}\label{th_approx_lsq}
$u_A^\varepsilon$ is smooth and for the remainder $r_A^\varepsilon:=\partial_t u_A^\varepsilon-\Delta u_A^\varepsilon+\frac{1}{\varepsilon^2}f'(u_A^\varepsilon)$ in the Allen-Cahn equation and $\varepsilon>0$ small it holds
\begin{alignat*}{2}\textstyle
|r_A^\varepsilon|&\leq C(\varepsilon e^{-c(|\rho_\varepsilon|+H_\pm^\varepsilon)}+\varepsilon^2e^{-c|\rho_\varepsilon|}+\varepsilon^3)&\quad&\text{ in }\Gamma^\pm(2\delta,1),\\
r_A^\varepsilon&=0&\quad&\text{ in }\Omega\textbackslash\Gamma(2\delta),\\
|\partial_\nu u_A^\varepsilon|&\leq C\varepsilon^2e^{-c|\rho_\varepsilon|}&\quad&\text{ on }\partial\Omega\cap\Gamma(2\delta),\\
\partial_\nu u_A^\varepsilon&=0&\quad&\text{ on }\partial\Omega\textbackslash\Gamma(2\delta).
\end{alignat*}
\end{Lemma}
\textit{Proof.} The second and the last equation are evident from the construction. Moreover, one can directly verify the above Taylor expansions rigorously to obtain 
\begin{alignat*}{2}\textstyle
|\partial_tu^I-\Delta u^I+\frac{1}{\varepsilon^2}f'(u^I)|&\textstyle\leq C(\varepsilon^2 e^{-c|\rho_\varepsilon|}+\varepsilon^3)&\textstyle\text{ in }\Gamma(2\delta),\\
\textstyle|\partial_tu^{B\pm}-\Delta u^{B\pm}+\frac{f'(u^I+u^{B\pm})-f'(u^I)}{\varepsilon^2}|&\textstyle\leq C(\varepsilon e^{-c(|\rho_\varepsilon|+H_\pm^\varepsilon)}+\varepsilon^3)&\textstyle\text{ in }\Gamma^\pm(2\delta,1),\\\textstyle
|\partial_\nu(u^I+u^{B\pm})|&\textstyle\leq C\varepsilon^2e^{-c|\rho_\varepsilon|}&\textstyle\text{ on }\Gamma(2\delta)\cap\partial\Omega,
\end{alignat*}
where we used the expansions \eqref{eq_inner_taylor_f} and \eqref{eq_bdry_taylor_f} for $m=4$ (i.e. up to $\Oc(\varepsilon^5)$) and that for all $\alpha>0$ any polynomial $p=p(\rho)$ or $p=p(\rho,H)$ can be estimated by $C_{p,\alpha}e^{\alpha|\rho|}$ and $C_{p,\alpha}e^{\alpha(|\rho|+H)}$, respectively. Then one uses the product rule in the definition of $u_A^\varepsilon$, the asymptotics of $\theta_0$ from Lemma \ref{th_theta_0} and the decay properties of $\hat{u}_j^I$ and $\hat{u}_j^{B\pm}$, $j=1,2,3$.\hfill$\square$\\
\section{Spectral Estimate}\label{sec_SE}
Let $T>0$ and $\Gamma=(\Gamma_t)_{t\in[0,T]}$ be as in Section \ref{sec_coord} and $\delta>0$ such that Theorem \ref{th_coord} holds for $2\delta$ instead of $\delta$. Throughout this section we assume the following structure: We define the stretched variables
\begin{align*}
\rho_\varepsilon:=\frac{r-\varepsilon h_\varepsilon(s,t)}{\varepsilon},\quad s^\pm:=\mp(s\mp 1),\quad H_\pm^\varepsilon:=\frac{s^\pm}{\varepsilon}\quad\text{ in }\overline{\Gamma(2\delta)},
\end{align*}
where $h_j\in B([0,T],C^0(I)\cap C^2(\tilde{I})), j=1,2$ with $\tilde{I}:=I\textbackslash(-1+2\delta_0,1-2\delta_0)$ and $h_\varepsilon:=h_1+\varepsilon h_2$. Then for $\varepsilon>0$ small it holds
\begin{align*}
u_A^\varepsilon=
\begin{cases}
\theta_0(\rho_\varepsilon)+\Oc(\varepsilon^2)&\quad\text{ in }\Gamma(2\delta,1-\delta_0),\\
\theta_0(\rho_\varepsilon)+\varepsilon u_1^{B\pm}+\Oc(\varepsilon^2)&\quad\text{ in }\Gamma^\pm(2\delta,2\delta_0),\\
\pm 1+\Oc(\varepsilon^2)&\quad\text{ in }Q_T^\pm\textbackslash\Gamma(\delta),
\end{cases}
\end{align*}
where $\Oc(\varepsilon^2)$ denotes measurable functions bounded by $C\varepsilon^2$ and with $\hat{u}_1^{B\pm}:\R^2_+\times[0,T]\rightarrow\R$ it holds
\[
u_1^{B\pm}(x,t)=\hat{u}_1^{B\pm}(\rho_\varepsilon(x,t),H_\varepsilon^\pm(x,t),t)\quad\text{ for }(x,t)\in\Gamma^\pm(2\delta,2\delta_0).
\] 
Here $\hat{u}_1^{B\pm}$ is smooth in $(\rho,H)$ and together with all derivatives in $(\rho,H)$ decays exponentially in $|\rho|+H$ uniformly in $t\in[0,T]$.\\
This structure is enough\footnote{ In particular the $h_j$ and $u_1^{B\pm}$ are not specified as in the asymptotic expansion and less regularity is needed. One can also allow an additional term of the form $\varepsilon\theta_1(\rho)p(s,t)$ in $u_A^\varepsilon$ with suitable $\theta_1, p$, cf. Remark \ref{th_bem_chen_zusatz}, 1. below.} to prove a spectral estimate for the (at $u_A^\varepsilon(.,t)$ linearized) Allen-Cahn operator for $t\in[0,T]$, which will be important when we estimate the difference of the exact and approximate solution. 
\begin{Theorem}\label{th_SE}
Let $M>0$ and $\|h_j\|_{B([0,T],C^0(I)\cap C^2(\tilde{I}))}\leq M$ for $j=1,2$. There are $\varepsilon_0,c_0>0$, independent of the $h_j$ for $M$ fixed, such that
\[
\int_\Omega|\nabla\psi|^2+\frac{1}{\varepsilon^2}f''(u_A^\varepsilon(.,t))\psi^2\,dx\geq -C\|\psi\|_{L^2(\Omega)}^2+\|\nabla\psi\|_{L^2(\Omega\textbackslash\Gamma_t(\delta))}^2+c_0\|\nabla_\tau\psi\|_{L^2(\Gamma_t(\delta))}^2
\]
for all $\varepsilon\in(0,\varepsilon_0], t\in[0,T]$ and $\psi\in H^1(\Omega)$.
\end{Theorem}
The main new difficulty compared to Theorem 2.3 in Chen \cite{ChenSpectrums} and Theorem 2.13 in Abels, Liu \cite{ALiu} is to prove a spectral estimate on $\Omega^{B\pm}_t:=\Gamma^\pm_t(\delta,2\delta_0)$, $t\in[0,T]$. This is the content of
\begin{Theorem}\label{th_SE_bdry}
Let the assumptions of Theorem \ref{th_SE} hold. There are $\tilde{\varepsilon}_0, \tilde{c}_0>0$, independent of the $h_j$ for $M$ fixed, such that
\[
\int_{\Omega^{B\pm}_t}|\nabla\psi|^2+\frac{1}{\varepsilon^2}f''(u_A^\varepsilon(.,t))\psi^2\,dx\geq -C\|\psi\|_{L^2(\Omega_t^{B\pm})}^2+\tilde{c}_0\|\nabla_\tau\psi\|_{L^2(\Omega_t^{B\pm})}^2
\]
for all $\varepsilon\in(0,\tilde{\varepsilon}_0], t\in[0,T]$ and $\psi\in H^1(\Omega^{B\pm}_t)$ with $\psi|_{X(.,s,t)}=0$ for $s\in(\frac{3}{2}\delta_0,2\delta_0)$.
\end{Theorem}
The additional assumption on $\psi$ is not needed but simplifies the proof. This version is enough to show Theorem \ref{th_SE}:\\
\newline
\textit{Proof of Theorem \ref{th_SE}.}
For $\varepsilon_0>0$ small and all $\varepsilon\in(0,\varepsilon_0]$ we have $f''(u_A^\varepsilon)\geq 0$ on $Q_T^\pm\textbackslash\Gamma(\delta)$. Therefore it is enough to prove the estimate for $\Gamma_t(2\delta)$ instead of $\Omega$. The estimate holds for $\Gamma_t(2\delta,1-\delta_0)$ instead of $\Omega$ with $c_0=1$ because there our curvilinear coordinate system coincides with the orthogonal one, cf. Theorem \ref{th_coord}, and hence we can use the arguments in Abels, Liu \cite{ALiu}, proof of Theorem 2.13. Moreover, using Theorem \ref{th_SE_bdry} we obtain the estimate for $\Gamma^\pm(2\delta,2\delta_0)$ instead of $\Omega$. Now we put all estimates together with a suitable partition of unity for 
\begin{align}\label{eq_se_proof1}
\Gamma_t(2\delta)\subseteq\overline{\Gamma_t(2\delta,1-\delta_0)}\cup\overline{\Gamma_t^-(2\delta,2\delta_0)}\cup\overline{\Gamma_t^+(2\delta,2\delta_0)}.
\end{align}
Therefore let $\eta_0,\eta_\pm:I\rightarrow[0,1]$ be a partition of unity subordinated to 
\[
[-1,-1+\frac{3}{2}\delta_0],[-1+\delta_0,1-\delta_0],[1-\frac{3}{2}\delta_0,1].
\]
W.l.o.g. $\sqrt{\eta_j}\in C_b^\infty(I)$ for $j=0,\pm$, otherwise we replace $\eta_j$ by $\eta_j^2/\sum_{j=0,\pm}\eta_j^2$. Then 
\[
\tilde{\eta}_j(.,t):\Gamma_t(2\delta)\rightarrow[0,1]:x\mapsto\eta_j(s(x,t))\quad\text{ for }j=0,\pm
\] 
defines a partition of unity for \eqref{eq_se_proof1} and $\chi_j(.,t):=\sqrt{\tilde{\eta}_j(.,t)}\in C_b^\infty(\overline{\Gamma_t(2\delta)}), j=0,\pm$. For any $\psi\in H^1(\Gamma_t(2\delta))$ it holds $\psi^2=\sum_{j=0,\pm}(\chi_j(.,t)\psi)^2$ and
\begin{align*}
\nabla(\chi_j(.,t)\psi)&=\nabla\chi_j(.,t)\psi+\chi_j(.,t)\nabla\psi,\\
|\nabla(\chi_j(.,t)\psi)|^2&=|\nabla\chi_j(.,t)|^2\psi^2+2\chi_j(.,t)\nabla\chi_j(.,t)\cdot\psi\nabla\psi+\chi_j(.,t)^2|\nabla\psi|^2.
\end{align*} 
Since $\sum_{j=0,\pm}\chi_j(.,t)^2=1$, we have $\sum_{j=0,\pm}\chi_j(.,t)\nabla\chi_j(.,t)=0$ and therefore
\[
\sum_{j=0,\pm}|\nabla(\chi_j(.,t)\psi)|^2=|\nabla\psi|^2+\psi^2\sum_{j=0,\pm}|\nabla\chi_j(.,t)|^2.
\]
This identity also holds for $\nabla_\tau$ instead of $\nabla$ which can be proven similarly. We write 
\begin{align*}
\int_{\Gamma_t(2\delta)}|\nabla\psi|^2+\frac{1}{\varepsilon^2}f''(u_A^\varepsilon(.,t))\psi^2\,dx=-\int_{\Gamma_t(2\delta)}\psi^2\sum_{j=0,\pm}|\nabla\chi_j(.,t)|^2\,dx+\\
+\sum_{j=0,\pm}\int_{\Gamma_t(2\delta)}|\nabla(\chi_j(.,t)\psi)|^2+\frac{1}{\varepsilon^2}f''(u_A^\varepsilon(.,t))(\chi_j(.,t)\psi)^2\,dx.
\end{align*}
Using the spectral estimates on $\Gamma_t(2\delta,1-\delta_0)$ and $\Gamma_t^\pm(2\delta,2\delta_0)$ and that $|\nabla\chi_j(.,t)|$ and $|\nabla_\tau\chi_j(.,t)|$ are bounded on $\Gamma_t(2\delta)$ uniformly in $t\in[0,T]$ due to $\sqrt{\eta_j}\in C_b^\infty(I)$, we obtain the spectral estimate in Theorem \ref{th_SE}.\hfill$\square$\\
\subsection{Proof of Theorem \ref{th_SE_bdry}.}
W.l.o.g. let $h_j(s,t)=0$ for $|s\mp 1|\in(\frac{7}{4}\delta_0,2\delta_0)$, otherwise we cut-off suitably. The additional assumption on $\psi$ ensures that this does not change the integral. Because of a Taylor expansion, it is enough to prove the estimate for 
\[
\frac{1}{\varepsilon^2}f''(\theta_0(\rho_\varepsilon(.,t)))+\frac{1}{\varepsilon}f'''(\theta_0(\rho_\varepsilon(.,t)))u_1^{B\pm}(.,t)
\] 
instead of $\frac{1}{\varepsilon^2}f''(u_A^\varepsilon(.,t))$. Similarly as in Alikakos, Chen, Fusco \cite{ACF}, Section 3 we seek an approximation $\phi_A^\varepsilon(.,t)$ to the first eigenfunction of 
\[
\Lc_{\varepsilon,t}^\pm:=-\Delta+\frac{1}{\varepsilon^2}f''(\theta_0(\rho_\varepsilon(.,t)))+\frac{1}{\varepsilon}f'''(\theta_0(\rho_\varepsilon(.,t)))u_1^{B\pm}(.,t)\quad\text{ on }\Omega_t^{B\pm}
\] 
with homogeneous Neumann-boundary condition and decompose $H^1(\Omega^{B\pm}_t)$ along the subspace of tangential variations of $\phi_A^\varepsilon(.,t)$. Therefore we make a suitable ansatz. This is motivated as follows:\\ 
\textbf{Motivation.} First consider the simpler operator $L=-\Delta+\frac{1}{\varepsilon^2}f''(\theta_0(\frac{\rho}{\varepsilon}))$ on a rectangle $[-\delta,\delta]\times[0,\eta]$ with homogeneous Neumann boundary condition. One can obtain all eigenvalues and eigenfunctions with a separation ansatz. Formally, because of Lemma \ref{th_ODE_lin1} and a scaling argument for small $\varepsilon>0$ the eigenfunctions corresponding to the lowest eigenvalues should have the form $a(s)\theta_0'(\frac{\rho}{\varepsilon})$ with $a:[0,\eta]\rightarrow\R$ and there should be a spectral gap. In $\Lc_{\varepsilon,t}^\pm$ there is $\rho_\varepsilon$ instead of $\rho$. Moreover, we have to deal with the $u_1^{B\pm}$-term at order $\frac{1}{\varepsilon}$ and we have to take into account the curvilinear structure of $\Omega^{B\pm}_t$.\\ 
Therefore we make the ansatz
\begin{align*}
\phi_A^\varepsilon(.,t)&:=\frac{1}{\sqrt{\varepsilon}}[v^I(.,t)+\varepsilon v^{B\pm}(.,t)],\\
v^I(.,t)&:=\theta_0'(\rho_\varepsilon(.,t))q^\pm(s^\pm(.,t),t),\quad v^{B\pm}(.,t):=\hat{v}^{B\pm}(\rho_\varepsilon(.,t),H_\varepsilon^\pm(.,t),t)\quad\text{ on }\Omega^{B\pm}_t,
\end{align*}
where $[0,2\delta_0]\times[0,T]\ni(\sigma,t)\mapsto q^\pm(\sigma,t)\in\R$ and $\hat{v}^{B\pm}:\R^2_+\times[0,T]\rightarrow\R$. Here the $\frac{1}{\sqrt{\varepsilon}}$-factor is multiplied for a certain normalization later.\\ 
In Subsection \ref{sec_SE_phi_A} we expand $\Lc_{\varepsilon,t}^\pm\phi_A^\varepsilon(.,t)$ and $\partial_\nu\phi_A^\varepsilon(.,t)$ similarly as before and choose $q^\pm$ and $\hat{v}^{B\pm}$ such that there is some cancellation. The $q^\pm$-term was introduced in order to fulfill the compatiblity condition for the equations for $\hat{v}^{B\pm}$. Then in Subsection \ref{sec_SE_BLF} we split $H^1(\Omega_t^{B\pm})$ orthogonally in $L^2(\Omega_t^{B\pm})$ with
\[
V_{\varepsilon,t}^\pm
:=\{\phi=a(s^\pm(.,t))\phi_A^\varepsilon(.,t)|a\in H^1(0,2\delta_0)\}
\] 
and analyze the bilinear form $B_{\varepsilon,t}^\pm$ corresponding to $\Lc_{\varepsilon,t}^\pm$ on every subspace. Again, there are some analogies to Alikakos, Chen, Fusco \cite{ACF}, Section 3. Here for $\phi,\psi\in H^1(\Omega_t^{B\pm})$
\begin{align}\label{eq_Bepst}
B_{\varepsilon,t}^\pm(\phi,\psi)=\int_{\Omega_t^{B\pm}}\nabla\phi\cdot\nabla\psi+\left[\frac{1}{\varepsilon^2}f''(\theta_0(\rho_\varepsilon(.,t)))+\frac{1}{\varepsilon}f'''(\theta_0(\rho_\varepsilon(.,t)))u_1^{B\pm}(x,t)\right]\phi\psi\,dx.
\end{align}
\subsubsection{Asymptotic Expansion for $\phi_A^\varepsilon$}\label{sec_SE_phi_A}
First we consider $\sqrt{\varepsilon}\Lc_{\varepsilon,t}^\pm\phi_A^\varepsilon(.,t)$. We expand $\Delta v^I$ as in the inner expansion in Section \ref{sec_inner_exp}: At lowest order $\Oc(\frac{1}{\varepsilon^2})$ we have $\frac{1}{\varepsilon^2}|\nabla r|^2|_{\overline{X}_0(s,t)}\theta_0'''(\rho)q^\pm(s^\pm,t)$ which cancels with $\frac{1}{\varepsilon^2}f''(\theta_0(\rho))\theta_0'(\rho)q^\pm(s^\pm,t)$ in $\Lc_{\varepsilon,t}^\pm\phi_A^\varepsilon(.,t)$. For the $\frac{1}{\varepsilon}$-order of $\Delta v^I$ we get 
\begin{align*}
\frac{1}{\varepsilon}\theta_0'''(\rho)q^\pm(s^\pm,t)\left[(\rho+h_1)\partial_r(|\nabla r|^2\circ\overline{X})|_{(0,s,t)}-2(\nabla r\cdot\nabla s)|_{\overline{X}_0(s,t)}\partial_sh_1\right]+\\
+\frac{1}{\varepsilon}\theta_0''(\rho)\left[\Delta r|_{\overline{X}_0(s,t)}q^\pm(s^\pm,t)\mp 2\partial_\sigma q^\pm(s^\pm,t)\nabla r\cdot\nabla s|_{\overline{X}_0(s,t)}\right]=\frac{1}{\varepsilon}\theta_0''(\rho)\Delta r|_{\overline{X}_0(s,t)}q^\pm(s^\pm,t).
\end{align*} 
As in Section \ref{sec_bdry_exp}, but without applying a Taylor expansion\footnote{ Since the $h_j$ are less regular and here we just need the lowest order.} for the $h_j$, we expand $\varepsilon\Delta v^{B\pm}$ and obtain at lowest order $\Oc(\frac{1}{\varepsilon})$:
\begin{align*}
\varepsilon\partial_H^2\hat{v}^{B\pm}\frac{|\nabla s|^2|_{\overline{p}^\pm(t)}}{\varepsilon^2}\mp\varepsilon 2\partial_\rho\partial_H\hat{v}^{B\pm}\frac{\nabla r\cdot\nabla s|_{\overline{p}^\pm(t)}}{\varepsilon^2}
+\varepsilon\partial_\rho^2\hat{v}^{B\pm}\frac{|\nabla r|^2|_{\overline{p}^\pm(t)}}{\varepsilon^2}=\frac{1}{\varepsilon}\Delta_t\hat{v}^{B\pm},
\end{align*}
where $\Delta_t:=\partial_\rho^2+\frac{4}{L(t)^2}\partial_H^2$. From the $f$-parts we have $\frac{f''(\theta_0(\rho))}{\varepsilon}\hat{v}^{B\pm}+\frac{1}{\varepsilon}f'''(\theta_0(\rho))\hat{u}_1^{B\pm}v^I$. To obtain an equation for $\hat{v}^{B\pm}$ in $(\rho,H,t)$, we use a Taylor expansion for $q^\pm$: 
\[
q^\pm(\varepsilon H,t)=q^\pm(0,t)+\Oc(\varepsilon H)\quad\text{ for }H\in[0,\frac{2\delta_0}{\varepsilon}].
\]
Therefore we require
\begin{align}\label{eq_se_bdry_vB_1}
(-\Delta_t+f''(\theta_0(\rho)))\hat{v}^{B\pm}=-f'''(\theta_0(\rho))\theta_0'\hat{u}_1^{B\pm}(\rho,H,t)q^\pm(0,t)\quad\text{ in }\R^2_+\times[0,T].
\end{align}
The remainder term $\frac{1}{\varepsilon}\Delta r|_{\overline{X}_0(s,t)}q^\pm(s^\pm,t)\theta_0''(\rho)$ seems bad, but there is still hope to gain a power in $\varepsilon$ later since $\int_\R\theta_0''\theta_0'\,d\rho=0$.\\
Now we expand $\sqrt{\varepsilon}\partial_\nu\phi_A^\varepsilon$ as in Section \ref{sec_bdry_exp_neum}. We obtain $\frac{1}{\varepsilon}(N_{\partial\Omega}\cdot\nabla r)|_{\overline{p}^\pm(t)}\theta_0''(\rho)q^\pm(0,t)$ for the $\frac{1}{\varepsilon}$-order. This is $0$ because of the $90$°-contact angle condition. The $\Oc(1)$-order is
\begin{align*}
q^\pm(0,t)\theta_0''(\rho)\left[(\rho+h_1|_{(\pm 1,t)})\partial_r((N_{\partial\Omega}\cdot\nabla r)\circ\overline{X})|_{(0,\pm 1,t)}-(N_{\partial\Omega}\cdot\nabla s)|_{\overline{p}^\pm(t)}\partial_sh_1|_{(\pm 1,t)}\right]+\\
\mp(N_{\partial\Omega}\cdot\nabla s)|_{\overline{p}^\pm(t)} \partial_\sigma q^\pm(0,t)\theta_0'(\rho)+0\cdot\partial_\rho\hat{v}^{B\pm}\mp(N_{\partial\Omega}\cdot\nabla s)|_{\overline{p}^\pm(t)}\partial_H\hat{v}^{B\pm}|_{H=0}.
\end{align*}
Because of Remark \ref{th_coord_bem}, 4. the cancellation is equivalent to
\begin{align*}
-\frac{2}{L(t)}\partial_H\hat{v}^{B\pm}|_{H=0}=g^{B\pm}(\rho,t)+\frac{2}{L(t)}\partial_\sigma q^\pm(0,t)\theta_0'(\rho),\\
g^{B\pm}(\rho,t):=q^\pm|_{(0,t)}\theta_0''(\rho)\left[\pm\frac{2}{L(t)}\partial_sh_1|_{(\pm 1,t)}-(\rho+h_1|_{(\pm 1,t)})\partial_r((N_{\partial\Omega}\cdot\nabla r)\circ\overline{X})_{(0,\pm 1,t)}\right].
\end{align*}
This equation together with \eqref{eq_se_bdry_vB_1} is equivalent to
\begin{alignat}{2}\label{eq_se_vBbar_1}
-\Delta\overline{v}^{B\pm}+f''(\theta_0(\rho))\overline{v}^{B\pm}&=-f'''(\theta_0(\rho))\theta_0'(\rho)\overline{u}_1^{B\pm}q^\pm(0,t)&\quad&\text{ in }\R^2_+\times[0,T],\\
-\partial_H\overline{v}^{B\pm}|_{H=0}&=g^{B\pm}+\frac{2}{L(t)}\partial_\sigma q^\pm(0,t)\theta_0'(\rho)&\quad&\text{ in }\R\times[0,T],\label{eq_se_vBbar_2}
\end{alignat}
where $\overline{v}^{B\pm},\overline{u}_1^{B\pm}:\R^2_+\times[0,T]\rightarrow\R$ are related to $\hat{v}^{B\pm}$ and $\hat{u}_1^{B\pm}$ as in Section \ref{sec_bdry_exp}. The right hand sides are smooth for fixed $t$ and decay together with all derivatives in $(\rho,H)$ exponentially in $|\rho|+H$ and $|\rho|$, respectively, uniformly in $t\in[0,T]$. We can take $q^\pm\in B([0,T],C^2([0,2\delta_0]))$ such that the solvability condition \eqref{eq_hp_intbed} corresponding to \eqref{eq_se_vBbar_1}-\eqref{eq_se_vBbar_2} holds and such that $q^\pm(0,t)=1, q^\pm(.,t)=1$ on $[\delta_0,2\delta_0]$ for all $t\in[0,T]$ and $c\leq q^\pm\leq C$ for some $c,C>0$. Theorem \ref{th_hp_time} yields a unique solution $\overline{v}^{B\pm}$ with the same regularity and decay as $\overline{u}_1^{B\pm}$.
\begin{Remark}\label{th_bem_chen_zusatz}\upshape
\begin{enumerate}
\item It is possible to add an additional term $\varepsilon\theta_1(\rho_\varepsilon)p(s,t)$ in $u_A^\varepsilon$, where 
\[
\theta_1\in C^\infty(\R)\cap L^\infty(\R)\quad\text{ with }\quad\int_\R f'''(\theta_0)\theta_1(\theta_0')^2\,d\rho=0
\] 
and $p\in B([0,T],C^0(I)\cap C^2(\tilde{I}))$. The spectral estimate can be shown with similar arguments. On the \enquote{orthogonal} part of the tubular neighbourhood this is basically done in Chen \cite{ChenSpectrums}, Theorem 2.3 and one can apply the same arguments as in the proof of Theorem 2.13 in Abels, Liu \cite{ALiu} to get the additional $\nabla_\tau$-term. Moreover, one adds $\varepsilon v_1(\rho_\varepsilon(.,t))q_1^\pm(s^\pm(.,t),t)$ in the ansatz for $\phi_A^\varepsilon$. Then in the expansion there is another term at order $\frac{1}{\varepsilon}$ we need to take care of, namely
\[
\left[-\frac{1}{\varepsilon^2}\varepsilon\partial_\rho^2v_1(\rho)+\frac{f''(\theta_0(\rho))}{\varepsilon^2}\varepsilon v_1(\rho)\right]q_1^\pm(s^\pm,t)+\frac{1}{\varepsilon}\left[f'''(\theta_0)\theta_1\theta_0'\right](\rho)p(s,t)q^\pm(s^\pm,t).
\]
Therefore we set $q_1^\pm(\sigma,t):=p(\pm1\mp\sigma)q^\pm(\sigma,t)$ and solve
\[
\left[-\partial_\rho^2+f''(\theta_0)\right]v_1=f'''(\theta_0)\theta_1\theta_0'
\] 
together with good decay for $v_1$ using Lemma \ref{th_ODE_lin2}. This is possible because of the integral condition for $\theta_1$. Then the same arguments apply. 
\item Consider the situation of Section \ref{sec_expansion}. Then $\overline{u}_1^{B\pm}$ solves \eqref{eq_u1B_1}-\eqref{eq_u1B_2}. By differentiating these equations with respect to $\rho$ we obtain
\begin{alignat*}{2}
-\Delta\partial_\rho\overline{u}_1^{B\pm}+f''(\theta_0)\partial_\rho\overline{u}_1^{B\pm}&=-f'''(\theta_0)\theta_0'\overline{u}_1^{B\pm}&\quad&\text{ in }\R^2_+\times[0,T],\\
-\partial_H(\partial_\rho\overline{u}_1^{B\pm})|_{H=0}&=g^{B\pm}-\theta_0'\partial_r((N_{\partial\Omega}\cdot\nabla r)\circ\overline{X})|_{(0,\pm 1,t)}&\quad&\text{ in }\R\times[0,T].
\end{alignat*}
Therefore in this case we can choose 
\[
\hat{v}^{B\pm}:=\partial_\rho\hat{u}_1^{B\pm}\quad\text{ and }\quad\partial_\sigma q^\pm(0,t)=-\frac{L(t)}{2}\partial_r((N_{\partial\Omega}\cdot\nabla r)\circ\overline{X})|_{(0,\pm 1,t)}.
\]\end{enumerate}
\end{Remark}
\begin{Lemma}\label{th_SE_phi_A}
The function $\phi_A^\varepsilon$ is $C^2(\overline{\Omega_t^{B\pm}})$ for fixed $t$ and satisfies uniformly in $t\in[0,T]$:
\begin{alignat*}{2}
|\sqrt{\varepsilon}\Lc_{\varepsilon,t}^\pm\phi_A^\varepsilon(.,t)+\frac{1}{\varepsilon}\Delta r|_{\overline{X}_0(s(.,t),t)}q^\pm|_{(s^\pm(.,t),t)}\theta_0''(\rho_\varepsilon(.,t))|&\leq Ce^{-c|\rho_\varepsilon(.,t)|}&\quad&\text{ in }\Omega_t^{B\pm},\\
|\sqrt{\varepsilon}\partial_\nu\phi_A^\varepsilon(.,t)|&\leq C\varepsilon e^{-c|\rho_\varepsilon(.,t)|}&\quad&\text{ on }\partial\Omega_t^{B\pm}\cap\partial\Omega,\\
|\sqrt{\varepsilon}\partial_\nu\phi_A^\varepsilon(.,t)|&\leq Ce^{-c/\varepsilon}&\quad&\text{ on }\partial\Omega_t^{B\pm}\textbackslash\partial\Omega.
\end{alignat*}
\end{Lemma}
\textit{Proof.} One can rigorously verify the Taylor expansions to get
\begin{align*}
|\sqrt{\varepsilon}\Lc_{\varepsilon,t}^\pm v^I(.,t)+\frac{1}{\varepsilon}\Delta r|_{\overline{X}_0(s,t)}q^\pm|_{(s^\pm,t)}\theta_0''(\rho_\varepsilon)|_{(.,t)}-\frac{1}{\varepsilon}f'''(\theta_0(\rho_\varepsilon))u_1^{B\pm}v^I|_{(.,t)}|\leq Ce^{-c|\rho_\varepsilon(.,t)|},\\
|\sqrt{\varepsilon}\Lc_{\varepsilon,t}^\pm\phi_A^\varepsilon(.,t)+\frac{1}{\varepsilon}\Delta r|_{\overline{X}_0(s,t)}q^\pm|_{(s^\pm,t)}\theta_0''(\rho_\varepsilon)|_{(.,t)}|\leq Ce^{-c|\rho_\varepsilon(.,t)|}+Ce^{-c(|\rho_\varepsilon(.,t)|+H_\varepsilon^\pm(.,t))}.
\end{align*} 
This shows the first estimate. The second one also directly follows from a Taylor expansion. On $\Omega_t^{B\pm}\textbackslash\Gamma_t(\frac{\delta}{2})$ we know that $\phi_A^\varepsilon(.,t)$, together with all derivatives, are $\Oc(e^{-\frac{c}{\varepsilon}})$. Therefore we only have to consider $\partial_\nu\phi_A^\varepsilon(.,t)$ on $(\partial\Omega_t^{B\pm}\cap\Gamma_t(\delta))\textbackslash\partial\Omega$: Here our coordinate system coincides with the orthogonal one and hence $\nabla r\cdot\nabla s|_{(.,t)}=0$ and the normal is given by $\mp\nabla s/|\nabla s||_{(.,t)}$ on this part of the boundary. Therefore because of $q^\pm(.,t)=1$ on $(\delta_0,2\delta_0)$
\[
\partial_\nu\theta_0'(\rho_\varepsilon(.,t))=\frac{2}{L(t)}\left[\pm\theta_0''(\rho_\varepsilon(.,t))\partial_sh_\varepsilon(s(.,t),t)+\theta_0'(\rho_\varepsilon(.,t))\right]\quad\text{ on }(\partial\Omega_t^{B\pm}\cap\Gamma_t(\delta))\textbackslash\partial\Omega.
\]
This vanishes since we assumed w.l.o.g. that $h_j(s,t)=0$ for $|s\mp 1|\in(\frac{7}{4}\delta_0,2\delta_0)$. Using the decay in $H$ of $\hat{v}^{B\pm}$ and its derivatives in $(\rho,H)$ the claim follows.\hfill$\square$\\
\subsubsection{Analysis of the Bilinear Form}\label{sec_SE_BLF}
We consider $V_{\varepsilon,t}^\pm:=
\{\phi=a(s^\pm(.,t))\phi_A^\varepsilon(.,t)|a\in H^1(0,2\delta_0)\}$, decompose $H^1(\Omega_t^{B\pm})$ orthogonally in $L^2(\Omega_t^{B\pm})$ and analyze $B_{\varepsilon,t}^\pm$ defined in \eqref{eq_Bepst} on every part. For simplicity we introduce 
\[
X^\pm:[-\delta,\delta]\times[0,2\delta_0]\times[0,T]\rightarrow\bigcup_{t\in[0,T]}\overline{\Omega_t^{B\pm}}\times\{t\}:(r,\sigma,t)\mapsto X(r,\pm1\mp\sigma,t),
\]
$\overline{X}^\pm:=(X^\pm,\textup{pr}_t)$ and similarly define $X_0^\pm$ and $\overline{X}_0^\pm$. Moreover, $J_t^\pm(r,\sigma):=J_t(r,\pm1\mp\sigma)$ for $(r,\sigma)\in[-\delta,\delta]\times[0,2\delta_0]$ and $t\in[0,T]$, where $J_t$ is defined in Remark \ref{th_coord_bem}, 3. Furthermore,
\[
h_j^\pm(\sigma,t):=h_j(\pm1\mp\sigma,t)\quad\text{ for all }\sigma\in[0,2\delta_0],t\in[0,T], j=1,2
\] 
and $h_\varepsilon^\pm:=h_1^\pm+\varepsilon h_2^\pm$. For transformation arguments we set 
\[
F_{\varepsilon,t}^\pm(z,\sigma):=\varepsilon(z+h_\varepsilon^\pm(\sigma,t))\quad\text{ and }\quad\tilde{J}_{\varepsilon,t,\sigma}^\pm(z):=J_t^\pm(F_{\varepsilon,t}^\pm(z,\sigma),\sigma)
\] 
for $|z+h_\varepsilon^\pm(\sigma,t)|\leq\frac{\delta}{\varepsilon}$ and $\sigma\in[0,2\delta_0],t\in[0,T]$. Moreover, we use
\begin{Lemma}\label{th_remainder}
Let $M>0$ and $\|h_j\|_{B([0,T],C^0(I)\cap C^2(\tilde{I}))}\leq M$ for $j=1,2$. Furthermore, let $\varepsilon_0\in(0,1]$ and $R_\varepsilon:\R\times[-2\delta,2\delta]\rightarrow\R$ be integrable for $\varepsilon\in(0,\varepsilon_0]$. Moreover, let $J\subseteq[-2\delta,2\delta]$ be an interval. Then
\begin{enumerate}
\item It holds 
\[
\int_J R_\varepsilon(\rho_\varepsilon|_{\overline{X}^\pm(r,\sigma,t)},r)\,dr=\varepsilon\int_{\frac{J}{\varepsilon}-h_\varepsilon^\pm(\sigma,t)} R_\varepsilon(z,F_{\varepsilon,t}^\pm(z,\sigma))\,dz
\]
for all $\varepsilon\in(0,\varepsilon_0], \sigma\in[0,2\delta_0]$ and $t\in[0,T]$.
\item If additionally $|R_\varepsilon(\rho,r)|\leq C_0|r|^k e^{-\alpha|\rho|}$ for all $(\rho,r)\in\R\times[-2\delta,2\delta]$ and some $k\in\N_0$, $C_0,\alpha>0$, then for all $\varepsilon\in(0,\varepsilon_0]$, $\sigma\in[0,2\delta_0]$ and $t\in[0,T]$
\[
\int_J |R_\varepsilon(\rho_\varepsilon|_{\overline{X}^\pm(r,\sigma,t)},r)|\,dr\leq C_0 C(\alpha,M,k)\varepsilon^{k+1}.
\] 
\end{enumerate}
\end{Lemma}
\textit{Proof.} The first assertion follows from the transformation rule. Using this formula for the second part, we obtain for all $\varepsilon\in(0,\varepsilon_0]$, $\sigma\in[0,2\delta_0]$ and $t\in[0,T]$
\[
\int_J|R_\varepsilon(\rho_\varepsilon|_{\overline{X}^\pm(r,\sigma,t)},r)|\,dr\leq C_0\varepsilon^{k+1}\int_\R|z+h_\varepsilon^\pm(\sigma,t)|^k e^{-\alpha|z|}\,dz.
\]
Because of $|h_\varepsilon^\pm(\sigma,t)|\leq 2M$ for all $\sigma\in[0,2\delta_0]$, $t\in[0,T]$, this yields the claim.\hfill$\square$\\
\newline
Let us characterize the splitting:
\begin{Lemma}\label{th_SE_split_L2}\begin{enumerate}
\item $V_{\varepsilon,t}^\pm$ is a subspace of $H^1(\Omega_t^{B\pm})$ and for $\varepsilon_0>0$ small there are constants $c_1,C_1>0$ with $c_1\|a\|_{L^2(0,2\delta_0)}\leq\|\phi\|_{L^2(\Omega_t^{B\pm})}\leq C_1\|a\|_{L^2(0,2\delta_0)}$ for all $\phi=a(s^\pm(.,t))\phi_A^\varepsilon(.,t)\in V_{\varepsilon,t}^\pm$ and $\varepsilon\in(0,\varepsilon_0], t\in[0,T]$.
\item We have $H^1(\Omega_t^{B\pm})=V_{\varepsilon,t}^\pm\oplus (V_{\varepsilon,t}^\pm)^\perp$ and for every $\psi\in H^1(\Omega_t^{B\pm})$ it holds
\[
\psi\in(V_{\varepsilon,t}^\pm)^\perp\quad\Leftrightarrow\quad\int_{-\delta}^\delta(\phi_A^\varepsilon(.,t)\psi)|_{X^\pm(r,\sigma,t)}J_t^\pm(r,\sigma)\,dr=0\quad\text{ f.a.a. }\sigma\in(0,2\delta_0).
\]
\end{enumerate}
\end{Lemma}
\textit{Proof. Ad 1.} We know that $\phi_A^\varepsilon(.,t)$ belongs to $C^2(\overline{\Omega_t^{B\pm}})$ for fixed $t\in[0,T]$. Moreover, it holds $a(s^\pm(.,t))\in H^1(\Omega_t^{B\pm})$ with Lemma \ref{th_prod_sets_spaces}, 3. and since smooth, bounded diffeomorphisms induce bounded linear operators between Sobolev spaces, cf. e.g. Adams and Fournier \cite{Adams}, Theorem 3.41. Therefore $V_{\varepsilon,t}^\pm$ is a subspace of $H^1(\Omega_t^{B\pm})$. Now we show the norm equivalence if $\varepsilon_0>0$ is small. The transformation rule and Fubini's Theorem yield
\begin{align}\label{eq_SE_L2}
\|\phi\|_{L^2(\Omega_t^{B\pm})}^2=\int_0^{2\delta_0}a(\sigma)^2\int_{-\delta}^\delta\phi_A^\varepsilon(\overline{X}^\pm(r,\sigma,t))^2 J_t^\pm(r,\sigma)\,dr\,d\sigma.
\end{align}
We compute the inner integral. The leading order term is $q^\pm(\sigma,t)$ times
\begin{align}\label{eq_SE_trafo}
\frac{1}{\varepsilon}\int_{-\delta}^\delta\theta_0' \left({\frac{r-\varepsilon h_\varepsilon^\pm(\sigma,t)}{\varepsilon}}\right)^2 J_t^\pm(r,\sigma)\,dr=\int_{-\frac{\delta}{\varepsilon}-h_\varepsilon^\pm(\sigma,t)} ^{\frac{\delta}{\varepsilon}-h_\varepsilon^\pm(\sigma,t)}\theta_0'(z)^2 \tilde{J}_{\varepsilon,t,\sigma}^\pm(z)\,dz,
\end{align}
where we used Lemma \ref{th_remainder}, 1. For $\varepsilon_0=\varepsilon_0(M)>0$ small it holds $|\varepsilon h_\varepsilon^\pm|\leq\frac{\delta}{2}$. Moreover, there are constants $c,C>0$ independent of $t\in[0,T]$ such that $c\leq J_t,q^\pm\leq C$. Therefore the integral in \eqref{eq_SE_trafo} can be estimated above and below by constants $\tilde{c},\tilde{C}>0$ independent of $t\in[0,T],\varepsilon\in(0,\varepsilon_0]$. For the remainder in the inner integral in \eqref{eq_SE_L2} we use Lemma \ref{th_remainder} and obtain an estimate of the absolute value to $C\varepsilon$. For $\varepsilon_0>0$ small this shows the claim.\hfill$\square_{1.}$\\  
\newline
\textit{Ad 2.} By definition
\[
(V_{\varepsilon,t}^\pm)^\perp=\left\{\psi\in H^1(\Omega_t^{B\pm}):\int_{\Omega_t^{B\pm}}\psi a(s^\pm(.,t))\phi_A^\varepsilon(.,t)\,dx=0\text{ for all }a\in H^1(0,2\delta_0)\right\}.
\]
Here the integral equals $\int_0^{2\delta_0}a(\sigma)\int_{-\delta}^\delta(\phi_A^\varepsilon(.,t)\psi)|_{X^\pm(r,\sigma,t)}J^\pm_t(r,\sigma)\,dr\,d\sigma$. Hence the Fundamental Theorem of Calculus of Variations yields the characterization. Moreover, the first part implies $V_{\varepsilon,t}^\pm\cap(V_{\varepsilon,t}^\pm)^\perp=\{0\}$ and it remains to prove $V_{\varepsilon,t}^\pm+(V_{\varepsilon,t}^\pm)^\perp=H^1(\Omega_t^{B\pm})$. Therefore we set
\[
w:[0,2\delta_0]\rightarrow\R:\sigma\mapsto\int_{-\delta}^\delta(\phi_A^\varepsilon|_{\overline{X}^\pm(r,\sigma,t)})^2J_t^\pm(r,\sigma)\,dr.
\]
Let $\psi\in H^1(\Omega_t^{B\pm})$ be arbitrary. Then we define
\[
a:[0,2\delta_0]\rightarrow\R:\sigma\mapsto\frac{1}{w(\sigma)}\int_{-\delta}^\delta(\phi_A^\varepsilon(.,t)\psi)|_{X^\pm(r,\sigma,t)}J_t^\pm(r,\sigma)\,dr.
\] 
It holds $w\in C^1([0,2\delta_0])$ and using the proof of 1. we have $w\geq c>0$. Because of Lemma \ref{th_prod_sets_spaces}, 2. and since integration gives a bounded linear operator on $L^2(-\delta,\delta)$, it follows that $a\in H^1(0,2\delta_0)$. For $\psi^\perp:=\psi-a(s^\pm(.,t))\phi_A^\varepsilon(.,t)\in H^1(\Omega_t^{B\pm})$ we have
\[
\int_{-\delta}^\delta(\phi_A^\varepsilon(.,t)\psi^\perp)|_{X^\pm(r,\sigma,t)}J_t^\pm(r,\sigma)\,dr=a(\sigma)w(\sigma)-a(\sigma)w(\sigma)=0\quad\text{ f.a.a. }\sigma\in(0,2\delta_0).
\]
The characterization above shows $\psi^\perp\in(V_{\varepsilon,t}^\pm)^\perp$.\hfill$\square$\\
\newline
First we consider $B_{\varepsilon,t}^\pm$ on $V_{\varepsilon,t}^\pm\times V_{\varepsilon,t}^\pm$.
\begin{Lemma}\label{th_SE_VxV}
There is an $\varepsilon_0>0$ such that for all $\varepsilon\in(0,\varepsilon_0],t\in[0,T]$ 
\[
B_{\varepsilon,t}^\pm(\phi,\phi)\geq-C\|\phi\|_{L^2(\Omega_t^{B\pm})}^2+c\|a\|_{H^1(0,2\delta_0)}^2\quad\text{ for all }\phi=a(s^\pm(.,t))\phi_A^\varepsilon(.,t)\in V_{\varepsilon,t}^\pm
\]
with $C,c>0$ independent of $\varepsilon,t$.
\end{Lemma}
\textit{Proof.} Let $\phi$ be as in the lemma. Then $\nabla\phi=\nabla(a(s^\pm(.,t)))\phi_A^\varepsilon(.,t)+a(s^\pm(.,t))\nabla\phi_A^\varepsilon(.,t)$, hence
\[
|\nabla\phi|^2=|\nabla(a(s^\pm))\phi_A^\varepsilon|^2|_{(.,t)}+a^2(s^\pm)|\nabla\phi_A^\varepsilon|^2|_{(.,t)}+\nabla(a^2(s^\pm))\cdot\nabla\phi_A^\varepsilon\phi_A^\varepsilon|_{(.,t)}.
\]
Integration by parts shows 
\begin{align*}
\int_{\Omega_t^{B\pm}}\left[\nabla(a^2(s^\pm))\cdot\nabla\phi_A^\varepsilon\phi_A^\varepsilon\right]|_{(.,t)}\,dx&=-\int_{\Omega_t^{B\pm}}\left[a^2(s^\pm)(\Delta\phi_A^\varepsilon\,\phi_A^\varepsilon+|\nabla\phi_A^\varepsilon|^2)\right]|_{(.,t)}\,dx+\\
&+\int_{\partial\Omega_t^{B\pm}}\left[\partial_\nu\phi_A^\varepsilon\,\tr(a^2(s^\pm)\phi_A^\varepsilon|_{(.,t)})\right]\,d\Hc^1.
\end{align*}
Therefore we obtain
\begin{align*}
B_{\varepsilon,t}^\pm(\phi,\phi)&=\int_{\Omega_t^{B\pm}}|\nabla(a(s^\pm))\phi_A^\varepsilon|^2|_{(.,t)}\,dx+\int_{\Omega_t^{B\pm}}(a^2(s^\pm)\phi_A^\varepsilon)|_{(.,t)}\Lc_{\varepsilon,t}^\pm\phi_A^\varepsilon|_{(.,t)}\,dx+\\
&+\int_{\partial\Omega_t^{B\pm}}\left[\partial_\nu\phi_A^\varepsilon\,\tr(a^2(s^\pm)\phi_A^\varepsilon|_{(.,t)})\right]\,d\Hc^1=:(I)+(II)+(III).
\end{align*}
\textit{Ad $(I)$}: It holds $|\nabla(a(s^\pm(.,t)))|^2=\left[|\nabla s|^2|(a')^2(s^\pm)\right]|_{(.,t)}$ and therefore 
\[
(I)=\int_0^{2\delta_0}(a')^2(\sigma)\int_{-\delta}^\delta\left[|\nabla s|^2(\phi_A^\varepsilon)^2\right]|_{\overline{X}^\pm(r,\sigma,t)}J_t^\pm(r,\sigma)\,dr\,d\sigma.
\]
Since $|\nabla s|,J_t\geq c>0$, we obtain similarly as in the proof of Lemma \ref{th_SE_split_L2}, 1. that $(I)\geq c_0\|a'\|_{L^2(0,2\delta_0)}^2$ for a $c_0>0$ and all $\varepsilon\in(0,\varepsilon_0], t\in[0,T]$, if $\varepsilon_0>0$ is small.\\
\newline
\textit{Ad $(II)$}: Analogously we write
\[
(II)=\int_0^{2\delta_0}a^2(\sigma)\int_{-\delta}^\delta \phi_A^\varepsilon|_{\overline{X}^\pm(r,\sigma,t)}(\Lc_{\varepsilon,t}^\pm\phi_A^\varepsilon(.,t))|_{X^\pm(r,\sigma,t)}J_t^\pm(r,\sigma)\,dr\,d\sigma
\]
and estimate the inner integral. Lemma \ref{th_SE_phi_A} implies
\[
|\sqrt{\varepsilon}\Lc_{\varepsilon,t}^\pm\phi_A^\varepsilon(.,t)+\frac{1}{\varepsilon}\Delta r|_{\overline{X}_0(s,t)}q^\pm(s^\pm(.,t),t)\theta_0''(\rho_\varepsilon(.,t))|\leq Ce^{-c|\rho_\varepsilon(.,t)|}\quad\text{ in }\Omega_t^{B\pm}.
\]
Then the most delicate term in the inner integral in $(II)$ is
\begin{align*}
\frac{1}{\varepsilon^2}\Delta r|_{\overline{X}^\pm_0(\sigma,t)}q^\pm(\sigma,t)\int_{-\delta}^\delta(\theta_0''\theta_0')(\rho_\varepsilon)|_{\overline{X}^\pm(r,\sigma,t)}J_t^\pm(r,\sigma)\,dr.
\end{align*}
Here because of a Taylor expansion and Remark \ref{th_coord_bem}, 4. it holds $J_t^\pm(r,\sigma)=\frac{L(t)}{2}+\phi^\pm(r,\sigma,t)$, where $|\phi^\pm(r,\sigma,t)|\leq C|r|$ uniformly in $(r,\sigma,t)$. Using Lemma \ref{th_remainder}, 1. and $\int_\R\theta_0''\theta_0'\,dz=0$, the $\frac{L(t)}{2}$-term is estimated by a constant $C>0$ independent of $t\in[0,T],\varepsilon\in(0,\varepsilon_0]$ if $\varepsilon_0>0$ is small. The remaining terms in $(II)$ can be estimated with Lemma \ref{th_remainder}. Altogether $|(II)|\leq C\|a\|_{L^2(0,2\delta_0)}^2$ for all $\varepsilon\in(0,\varepsilon_0], t\in[0,T]$ if $\varepsilon_0>0$ is small.\\
\newline
\textit{Ad $(III)$}: It holds
\begin{align*}
(III)&=\sum_\pm\int_0^{2\delta_0}a^2(\sigma)\left[\phi_A^\varepsilon\partial_\nu\phi_A^\varepsilon\right]|_{\overline{X}^\pm(\pm\delta,\sigma,t)}|\partial_\sigma X^\pm(\pm\delta,\sigma,t)|\,d\sigma+\\
&+\sum_{\sigma=0,2\delta_0}a^2(\sigma)\int_{-\delta}^\delta \left[\phi_A^\varepsilon\partial_\nu\phi_A^\varepsilon\right]|_{\overline{X}^\pm(r,\sigma,t)}|\partial_rX^\pm(r,\sigma,t)|\,dr.
\end{align*}
Using Lemma \ref{th_SE_phi_A} and Lemma \ref{th_remainder} for the integral where $\sigma=0$ is fixed, we obtain
\[
|(III)|\leq Ce^{-c/\varepsilon}\|a\|_{L^2(0,2\delta_0)}^2+Ce^{-c/\varepsilon}a^2(2\delta_0)+C\varepsilon a^2(0).
\]
We use $H^1(0,2\delta_0)\hookrightarrow C_b^0([0,2\delta_0])$ and put all estimates together. Then by Lemma \ref{th_SE_split_L2}, 1. the claim follows.\hfill$\square$\\
\newline
Next we analyze $B^\pm_{\varepsilon,t}$ on $(V_{\varepsilon,t}^\pm)^\perp\times(V_{\varepsilon,t}^\pm)^\perp$.
\begin{Lemma}\label{th_SE_VperpxVperp}
There are $\varepsilon_0,\nu_0>0$ such that for all $\varepsilon\in(0,\varepsilon_0],t\in[0,T]$ and $\psi\in(V^\pm_{\varepsilon,t})^\perp$ it holds
\[
B^\pm_{\varepsilon,t}(\psi,\psi)\geq\nu_0\left[\frac{1}{\varepsilon^2}\|\psi\|_{L^2(\Omega_t^{B\pm})}^2+\|\nabla\psi\|_{L^2(\Omega_t^{B\pm})}^2\right].
\]  
\end{Lemma}
\textit{Proof.} It is enough to prove the existence of $\tilde{\varepsilon}_0,\tilde{\nu}_0>0$ such that for all $\varepsilon\in(0,\tilde{\varepsilon}_0], t\in[0,T]$ and $\psi\in(V^\pm_{\varepsilon,t})^\perp$ it holds
\begin{align}\label{eq_SE_BLF_orth1}
\tilde{B}^\pm_{\varepsilon,t}(\psi,\psi):=\int_{\Omega^{B\pm}_t}|\nabla\psi|^2+\frac{f''(\theta_0(\rho_\varepsilon(.,t)))}{\varepsilon^2}\psi^2\,dx\geq\frac{\tilde{\nu}_0}{\varepsilon^2}\|\psi\|_{L^2(\Omega^{B\pm}_t)}^2.
\end{align}
Then the claim follows because for $\alpha>0$ small
\begin{align*}
B_{\varepsilon,t}^\pm(\psi,\psi)\geq\tilde{B}_{\varepsilon,t}^\pm(\psi,\psi)-\frac{C}{\varepsilon}\|\psi\|_{L^2(\Omega^{B\pm}_t)}=(1-\alpha+\alpha)\tilde{B}_{\varepsilon,t}^\pm(\psi,\psi)-\frac{C}{\varepsilon}\|\psi\|_{L^2(\Omega^{B\pm}_t)}\geq\\
\geq{\textstyle\frac{(1-\alpha)\tilde{\nu}_0-C(\alpha+\varepsilon)}{\varepsilon^2}}\|\psi\|_{L^2(\Omega^{B\pm}_t)}^2+\alpha\|\nabla\psi\|_{L^2(\Omega^{B\pm}_t)}^2\geq\nu_0\left[\frac{1}{\varepsilon^2}\|\psi\|_{L^2(\Omega^{B\pm}_t)}^2+\|\nabla\psi\|_{L^2(\Omega^{B\pm}_t)}^2\right],
\end{align*}
if $\nu_0,\varepsilon_0>0$ are small and $\varepsilon\in(0,\varepsilon_0], t\in[0,T]$ and $\psi\in(V_{\varepsilon,t}^\pm)^\perp$.\\
In the following we prove \eqref{eq_SE_BLF_orth1} by reducing to Neumann-boundary problems in normal direction which is also the idea in Chen \cite{ChenSpectrums}, proof of Theorem 2.3. Therefore we define $\tilde{\psi}_t^\pm:=\psi|_{X^\pm(.,t)}$ for $t\in[0,T]$. It holds
\[
|\nabla\psi|^2|_{X^\pm(.,t)}=(\nabla_{(r,s)}\tilde{\psi}_t^\pm)^\top
\begin{pmatrix}
|\nabla r|^2 & \mp\nabla r\cdot\nabla s\\
\mp\nabla r\cdot\nabla s & |\nabla s|^2
\end{pmatrix}|_{\overline{X}^\pm(.,t)}
\nabla_{(r,s)}\tilde{\psi}_t^\pm.
\]
Theorem \ref{th_coord}, a Taylor expansion and Young's inequality imply 
\[
|\nabla\psi|^2|_{X^\pm(.,t)}\geq(1-Cr^2)(\partial_r\tilde{\psi}_t^\pm)^2
\] 
for a $C>0$. To get $Cr^2$ small enough (which will be precise later), we fix $\tilde{\delta}>0$ small and estimate for $r$ on
\[
I_{\varepsilon,t,\sigma}^\pm:=(-\tilde{\delta},\tilde{\delta})+\varepsilon h_\varepsilon^\pm(\sigma,t)\quad\text{ and }\quad \tilde{I}_{\varepsilon,t,\sigma}^\pm:=(-\delta,\delta)\textbackslash I_{\varepsilon,t,\sigma}^\pm.
\]
If $\varepsilon_0=\varepsilon_0(\tilde{\delta},M)>0$ is small, then $f''(\theta_0(\rho_\varepsilon(\overline{X}^\pm(r,\sigma,t))))\geq c>0$ for $r\in\tilde{I}_{\varepsilon,t,\sigma}^\pm$ and $|r|\leq \tilde{\delta}+\varepsilon|h_\varepsilon^\pm(r,\sigma)|\leq 2\tilde{\delta}$ for $r\in I_{\varepsilon,t,\sigma}^\pm$ for all $\varepsilon\in(0,\varepsilon_0], t\in[0,T]$ and $\sigma\in[0,2\delta_0]$. Therefore with $\tilde{c}:=2C\tilde{\delta}^2$ we obtain 
\begin{align*}
\tilde{B}_{\varepsilon,t}^\pm(\psi,\psi)&\geq\int_0^{2\delta_0}\int_{\tilde{I}_{\varepsilon,t,\sigma}^\pm}\frac{c}{\varepsilon^2}(\tilde{\psi}_t^\pm)^2 J_t^\pm|_{(r,\sigma)}\,dr\,d\sigma+\\
&+\int_0^{2\delta_0}\int_{I_{\varepsilon,t,\sigma}^\pm}\left[(1-\tilde{c})(\partial_r\tilde{\psi}_t^\pm)^2+\frac{f''(\theta_0(\rho_\varepsilon))}{\varepsilon^2}|_{\overline{X}^\pm(.,t)}(\tilde{\psi}_t^\pm)^2\right] J_t^\pm|_{(r,\sigma)}\,dr\,d\sigma.
\end{align*}
Using Lemma \ref{th_remainder} for the second part, Lemma \ref{th_SE_VperpxVperp} follows if we show that for 
\[
\Psi_{\varepsilon,t,\sigma}^\pm:=\sqrt{\varepsilon}\tilde{\psi}_t^\pm(F_{\varepsilon,t}^\pm(.,\sigma),\sigma)\quad\text{ on }I_\varepsilon:=(-\frac{\tilde{\delta}}{\varepsilon},\frac{\tilde{\delta}}{\varepsilon})
\] 
and f.a.a. $\sigma\in[0,2\delta_0]$, all $t\in[0,T]$ and $\varepsilon\in(0,\varepsilon_0]$, $\varepsilon_0>0$ small it holds
\begin{align}\begin{split}\label{eq_SE_VperpxVperp_proof}
B_{\varepsilon,t,\sigma}^{\pm,\tilde{c}}(\Psi_{\varepsilon,t,\sigma}^\pm,\Psi_{\varepsilon,t,\sigma}^\pm)&:=\int_{I_\varepsilon}\left[(1-\tilde{c})(\partial_z\Psi_{\varepsilon,t,\sigma}^\pm)^2+f''(\theta_0(z))(\Psi_{\varepsilon,t,\sigma}^\pm)^2\right]\tilde{J}_{\varepsilon,t,\sigma}^\pm\,dz\\
&\geq c\|\Psi_{\varepsilon,t,\sigma}^\pm\|^2_{L^2(I_\varepsilon,\tilde{J}_{\varepsilon,t,\sigma}^\pm)},\end{split}
\end{align}
where $c>0$ is independent of $\varepsilon\in(0,\varepsilon_0], t\in[0,T]$ and $\sigma\in[0,2\delta_0]$ and $L^2(I_\varepsilon,\tilde{J}_{\varepsilon,t,\sigma}^\pm)$ is the space of $L^2$-functions on $I_\varepsilon$ with respect to the weight $\tilde{J}_{\varepsilon,t,\sigma}^\pm$. We denote the scalar-product in $L^2(I_\varepsilon,\tilde{J}_{\varepsilon,t,\sigma}^\pm)$ by $(.,.)_{\varepsilon,t,\sigma}$, the norm with $\|.\|_{\varepsilon,t,\sigma}$ and the orthogonal relation by $\perp_{\varepsilon,t,\sigma}$. For the proof of \eqref{eq_SE_VperpxVperp_proof} we need results for $B^{\pm,0}_{\varepsilon,t,\sigma}$. With respect to $(.,.)_{\varepsilon,t,\sigma}$, $B_{\varepsilon,
t,\sigma}^{\pm,0}$ is the bilinear form corresponding to
\[
\Lc_{\varepsilon,t,\sigma}^{\pm,0}:=-(\tilde{J}_{\varepsilon,t,\sigma}^\pm)^{-1}\frac{d}{dz}\left(\tilde{J}_{\varepsilon,t,\sigma}^\pm\frac{d}{dz}\right)+f''(\theta_0)
\]
on $H^2(I_\varepsilon)$ with homogeneous Neumann boundary condition.
\begin{Lemma}\label{th_Chenspectrums}
For $\tilde{\delta}\in(0,\frac{3\delta}{4}]$ fixed, $\varepsilon>0$ small, $t\in[0,T]$ and $\sigma\in[0,2\delta_0]$ it holds
\begin{enumerate}
\item The first eigenvalue $\lambda_1^\pm=\lambda_1^\pm(\varepsilon,t,\sigma)$ of $\Lc_{\varepsilon,t,\sigma}^{\pm,0}$ is simple, it admits a positive, normalized eigenfunction $\Psi_1^\pm=\Psi_1^\pm(\varepsilon,t,\sigma)$ and 
satisfies uniformly in $t,\sigma$
\[
\lambda_1^\pm=\inf_{\Psi\in H^1(I_\varepsilon),\|\Psi\|_{\varepsilon,t,\sigma}=1}B_{\varepsilon,t,\sigma}^{\pm,0}(\Psi,\Psi)=B_{\varepsilon,t,\sigma}^{\pm,0}(\Psi_1^\pm,\Psi_1^\pm)=\Oc(\varepsilon^2).
\]
Moreover, $\Psi_{1,R}^\pm:=\Psi_1^\pm-J_t^\pm(0,\sigma)^{-\frac{1}{2}}\beta_\varepsilon\theta_0'$, where $\beta_\varepsilon:=\|\theta_0'\|_{L^2(I_\varepsilon)}^{-1}$, fulfills uniformly in $t,\sigma$
\[
\|\Psi_{1,R}^\pm\|_{\varepsilon,t,\sigma}+\|\partial_z\Psi_{1,R}^\pm\|_{\varepsilon,t,\sigma}=\Oc(\varepsilon).
\]
\item There is a $v_2>0$ independent of $\tilde{\delta}$ such that for all $\varepsilon\in(0,\varepsilon_2]$, $\varepsilon_2=\varepsilon_2(\tilde{\delta})>0$ small
\[
\lambda_2^\pm=\lambda_2^\pm(\varepsilon,t,\sigma)=\inf_{\Psi\in H^1(I_\varepsilon),\|\Psi\|_{\varepsilon,t,\sigma}=1, \Psi\perp_{\varepsilon,t,\sigma}\Psi_1}B_{\varepsilon,t,\sigma}^{\pm,0}(\Psi,\Psi)\geq v_2.
\]
\end{enumerate}
\end{Lemma}
\textit{Proof.} This can be shown analogously to Lemma 2.2 in Chen \cite{ChenSpectrums}. Note that in Chen \cite{ChenSpectrums} it is assumed w.l.o.g. that $\tilde{\delta}=1$ and here we additionally introduced some height functions. Moreover, the coordinate system for the interface in Chen \cite{ChenSpectrums} is the orthogonal one. Nevertheless, $\tilde{J}_{\varepsilon,t,\sigma}^\pm$ behaves similar as $J(\varepsilon.,s)$ in Chen \cite{ChenSpectrums} and therefore an analogous proof works. Additionally, we obtain that $v_2$ does not depend on $\tilde{\delta}$.\hfill$\square$\\
\newline
\textit{Proof of \eqref{eq_SE_VperpxVperp_proof}.} Using $\psi\in (V_{\varepsilon,t}^\pm)^\perp$, Lemma \ref{th_SE_split_L2}, 2. and $0<c\leq q^\pm\leq C$, we obtain 
\[
\left|\int_{I_{\varepsilon,t,\sigma}^\pm}\frac{1}{\sqrt{\varepsilon}}\theta_0'(\rho_\varepsilon)|_{\overline{X}^\pm(r,\sigma,t)}[\tilde{\psi}_tJ_t^\pm]|_{(r,\sigma)}\,dr\right|\leq
C\int_{I_{\varepsilon,t,\sigma}^\pm}|\sqrt{\varepsilon}v^{B\pm}|_{\overline{X}^\pm(r,\sigma,t)}[\tilde{\psi}_tJ_t^\pm]|_{(r,\sigma)}|\,dr
\]
f.a.a. $\sigma\in[0,2\delta_0]$, all $\varepsilon\in(0,\varepsilon_0]$ and $t\in[0,T]$. Hölder's inequality and Lemma \ref{th_remainder} imply that the right hand side is estimated by $C\varepsilon\|\tilde{\psi}_t^\pm(.,\sigma)\|_{L^2(I_{\varepsilon,t,\sigma}^\pm,J_t^\pm(.,\sigma))}$. Therefore Lemma \ref{th_remainder}, 1. yields
\[
|(\Psi_{\varepsilon,t,\sigma}^\pm,\theta_0')_{\varepsilon,t,\sigma}|\leq C\varepsilon\|\Psi_{\varepsilon,t,\sigma}^\pm\|_{\varepsilon,t,\sigma}.
\]
Hence Lemma \ref{th_Chenspectrums}, 1. and boundedness of $\beta_\varepsilon$ for small $\varepsilon$ imply
\begin{align}\label{eq_SE_VperpxVperp_proof2}
|(\Psi_{\varepsilon,t,\sigma}^\pm,\Psi_1^\pm)_{\varepsilon,t,\sigma}|\leq C\varepsilon\|\Psi_{\varepsilon,t,\sigma}^\pm\|_{\varepsilon,t,\sigma}
\end{align}
f.a.a. $\sigma\in[0,2\delta_0]$, all $t\in[0,T]$ and $\varepsilon\in(0,\varepsilon_0]$, if $\varepsilon_0>0$ is small. Now we decompose $\Psi_{\varepsilon,t,\sigma}^\pm$ orthogonally in $L^2(I_\varepsilon,\tilde{J}_{\varepsilon,t,\sigma}^\pm)$: With some $\Psi_\perp^\pm=\Psi_\perp^\pm(\varepsilon,t,\sigma)$
\[
\frac{\Psi_{\varepsilon,t,\sigma}^\pm}{\|\Psi_{\varepsilon,t,\sigma}^\pm\|_{\varepsilon,t,\sigma}}=\left(\Psi_1^\pm,\frac{\Psi_{\varepsilon,t,\sigma}^\pm}{\|\Psi_{\varepsilon,t,\sigma}^\pm\|_{\varepsilon,t,\sigma}}\right)_{\varepsilon,t,\sigma}\Psi_1^\pm+\Psi_\perp^\pm.
\]
Taking $\|.\|_{\varepsilon,t,\sigma}^2$ in this identity yields $|1-\|\Psi_\perp^\pm\|_{\varepsilon,t,\sigma}^2|\leq C\varepsilon^2$ because of \eqref{eq_SE_VperpxVperp_proof2}. Then
\begin{align*}
B_{\varepsilon,t,\sigma}^{\pm,\tilde{c}}(\Psi_{\varepsilon,t,\sigma}^\pm,\Psi_{\varepsilon,t,\sigma}^\pm)=
(1-\tilde{c})B_{\varepsilon,t,\sigma}^{\pm,0}(\Psi_{\varepsilon,
t,\sigma}^\pm,\Psi_{\varepsilon,
t,\sigma}^\pm)+\tilde{c}\int_{I_\varepsilon}f''(\theta_0)(\Psi_{\varepsilon,
t,\sigma}^\pm)^2\tilde{J}_{\varepsilon,t,\sigma}^\pm\,dz.
\end{align*}
The last part is estimated by $\tilde{c}\,\sup_{z\in\R}|f''(\theta_0(z))|\,\|\Psi_{\varepsilon,t,\sigma}^\pm\|_{\varepsilon,t,\sigma}^2$. Moreover, with \eqref{eq_SE_VperpxVperp_proof2}
\[
\frac{1}{\|\Psi_{\varepsilon,t,\sigma}^\pm\|^2_{\varepsilon,t,\sigma}}B_{\varepsilon,t,\sigma}^{\pm,0}(\Psi_{\varepsilon,t,\sigma}^\pm,\Psi_{\varepsilon,t,\sigma}^\pm)=\lambda_1^\pm\Oc(\varepsilon^2)+2B_{\varepsilon,t,\sigma}^{\pm,0}(\Psi_1^\pm,\Psi_\perp^\pm)\Oc(\varepsilon)+B_{\varepsilon,t,\sigma}^{\pm,0}(\Psi_\perp^\pm,\Psi_\perp^\pm).
\]
Lemma \ref{th_Chenspectrums} implies $\lambda_1^\pm=\Oc(\varepsilon^2)$ and $B_{\varepsilon,t,\sigma}^{\pm,0}(\Psi_\perp^\pm,\Psi_\perp^\pm)\geq v_2\|\Psi_\perp^\pm\|_{\varepsilon,t,\sigma}^2\geq v_2(1-C\varepsilon^2)$. Since 
\[
B_{\varepsilon,t,\sigma}^{\pm,0}(\Psi_1^\pm,\Psi_\perp^\pm)=(\Lc_{\varepsilon,t,\sigma}^{\pm,0}\Psi_1^\pm,\Psi_\perp^\pm)_{\varepsilon,t,\sigma}=\lambda_1^\pm(\Psi_1^\pm,\Psi_\perp^\pm)_{\varepsilon,t,\sigma}=0
\]
due to integration by parts, it follows that
\[
\frac{B_{\varepsilon,t,\sigma}^{\pm,\tilde{c}}(\Psi_{\varepsilon,t,\sigma}^\pm,\Psi_{\varepsilon,t,\sigma}^\pm)}{\|\Psi_{\varepsilon,t,\sigma}^\pm\|^2_{\varepsilon,t,\sigma}}\geq(1-\tilde{c})\left[-C\varepsilon^4+v_2(1-C\varepsilon^2)\right]-\tilde{c}\,\sup_{\R}|f''(\theta_0)|\geq c>0
\]
f.a.a. $\sigma\in[0,2\delta_0]$, all $t\in[0,T]$ and $\varepsilon\in(0,\varepsilon_0]$ if $\varepsilon_0=\varepsilon_0(\tilde{\delta},M)>0$ is small and $\tilde{\delta}>0$ (and therefore $\tilde{c}=2C\tilde{\delta}^2$) was chosen small enough before. Here it is important that $v_2$ in Lemma \ref{th_Chenspectrums} is independent of $\tilde{\delta}$. This shows \eqref{eq_SE_VperpxVperp_proof} and thus Lemma \ref{th_SE_VperpxVperp}.\hfill$\square$\\ 
\newline
For $B_{\varepsilon,t}^\pm$ on $V_{\varepsilon,t}^\pm\times(V_{\varepsilon,t}^\pm)^\perp$ we have
\begin{Lemma}\label{th_SE_VxVperp}
There is an $\varepsilon_0>0$ such that for all $\varepsilon\in(0,\varepsilon_0], t\in[0,T]$
\[
|B_{\varepsilon,t}^\pm(\phi,\psi)|\leq\frac{C}{\varepsilon}\|\phi\|_{L^2(\Omega_t^{B\pm})}\|\psi\|_{L^2(\Omega_t^{B\pm})}+\frac{1}{4}B_{\varepsilon,t}^\pm(\psi,\psi)+C\varepsilon^2\|a\|_{H^1(0,2\delta_0)}^2
\]
for all $\phi\in V_{\varepsilon,t}^\pm$ and $\psi\in (V_{\varepsilon,t}^\pm)^\perp$.
\end{Lemma}
\textit{Proof.} We have $\phi=a(s^\pm(.,t))\phi_A^\varepsilon(.,t)$ with some $a\in H^1(0,2\delta_0)$. By product rule $\nabla\phi=\nabla(a(s^\pm(.,t)))\phi_A^\varepsilon+a(s^\pm(.,t))\nabla\phi_A^\varepsilon(.,t)$ and integration by parts yields
\begin{align*}
\int_{\Omega_t^{B\pm}}a(s^\pm)\nabla\phi_A^\varepsilon|_{(.,t)}\cdot\nabla\psi\,dx=-\int_{\Omega_t^{B\pm}}\nabla(a(s^\pm))\cdot\nabla\phi_A^\varepsilon|_{(.,t)}\psi\,dx+\\
-\int_{\Omega_t^{B\pm}}a(s^\pm)\Delta\phi_A^\varepsilon|_{(.,t)}\psi\,dx+\int_{\partial\Omega_t^{B\pm}}\partial_\nu\phi_A^\varepsilon|_{(.,t)}\tr\left[a(s^\pm(.,t))\psi\right]\,d\Hc^1.
\end{align*}
Therefore we obtain
\begin{align*}
B_{\varepsilon,t}^\pm(\phi,\psi)=\int_{\Omega_t^{B\pm}}a(s^\pm)|_{(.,t)}\Lc_{\varepsilon,t}^\pm\phi_A^\varepsilon|_{(.,t)}\psi\,dx+\int_{\partial\Omega_t^{B\pm}}\partial_\nu\phi_A^\varepsilon|_{(.,t)}\tr\left[a(s^\pm(.,t))\psi\right]\,d\Hc^1+\\
+\int_{\Omega_t^{B\pm}}\nabla(a(s^\pm))|_{(.,t)}\cdot(\phi_A^\varepsilon|_{(.,t)}\nabla\psi-\nabla\phi_A^\varepsilon|_{(.,t)}\psi)\,dx=:(I)+(II)+(III).
\end{align*}
\textit{Ad $(I)$}: Hölder's inequality yields $|(I)|\leq\|a(s^\pm)|_{(.,t)}\Lc_{\varepsilon,t}^\pm\phi_A^\varepsilon|_{(.,t)}\|_{L^2(\Omega_t^{B\pm})}\|\psi\|_{L^2(\Omega_t^{B\pm})}$, where
\[
\|a(s^\pm)|_{(.,t)}\Lc_{\varepsilon,t}^\pm\phi_A^\varepsilon|_{(.,t)}\|_{L^2(\Omega_t^{B\pm})}^2=\int_0^{2\delta_0}a^2(\sigma)\int_{-\delta}^\delta (\Lc_{\varepsilon,t}^\pm\phi_A^\varepsilon(.,t))^2|_{X^\pm(r,\sigma,t)}\,J_t^\pm(r,\sigma)\,dr\,d\sigma.
\]
To estimate the inner integral we use Lemma \ref{th_SE_phi_A} to obtain
\begin{align*}
|\varepsilon(\Lc_{\varepsilon,t}^\pm\phi_A^\varepsilon|_{(.,t)})^2-\frac{1}{\varepsilon^2}(\Delta r|_{\overline{X}_0(s(.,t),t)}q^\pm|_{(s^\pm(.,t),t)}\theta_0''(\rho_\varepsilon(.,t)))^2|\leq\\ 
\leq Ce^{-c|\rho_\varepsilon(.,t)|}(Ce^{-c|\rho_\varepsilon(.,t)|}+\frac{2}{\varepsilon}|\Delta r||_{\overline{X}_0(s(.,t),t)}q^\pm|_{(s^\pm(.,t),t)}).
\end{align*}
With Lemma \ref{th_remainder} we get that the inner integral is estimated by $\frac{C}{\varepsilon^2}$ and hence $|(I)|\leq \frac{C}{\varepsilon}\|a\|_{L^2(0,2\delta_0)}\|\psi\|_{L^2(\Omega_t^{B\pm})}$ for all $t\in[0,T]$ and $\varepsilon\in(0,\varepsilon_0]$ if $\varepsilon_0>0$ is small. Because of Lemma \ref{th_SE_split_L2}, 1. this fits the desired estimate.\\
\newline
\textit{Ad $(II)$}: It holds $|(II)|\leq\|\tr\,\psi\|_{L^2(\partial\Omega_t^{B\pm})}\|\tr(a(s^\pm(.,t)))\partial_\nu\phi_A^\varepsilon|_{(.,t)}\|_{L^2(\partial\Omega_t^{B\pm})}$ because of Hölder's inequality. Similarly as in the proof of Lemma \ref{th_SE_VxV} we have
\[
\|a(s^\pm)\partial_\nu\phi_A^\varepsilon|_{(.,t)}\|_{L^2(\partial\Omega_t^{B\pm})}\leq C\varepsilon|a(0)|+ Ce^{-c/\varepsilon}(\|a\|_{L^2(0,2\delta_0)}+|a(2\delta_0)|)\leq C\varepsilon\|a\|_{H^1(0,2\delta_0)}.
\]
To estimate $\|\tr\,\psi\|_{L^2(\partial\Omega_t^{B\pm})}$ we need
\begin{Lemma}\label{th_inpol_trace}
Fix $\overline{\varepsilon}>0$. Then there is a $C>0$ such that for all $\varepsilon\in(0,\overline{\varepsilon}]$ and $t\in[0,T]$
\[
\|\tr\,\psi\|_{L^2(\partial\Omega_t^{B\pm})}^2\leq C(\varepsilon\|\nabla \psi\|_{L^2(\Omega_t^{B\pm})}^2+\frac{1}{\varepsilon}\|\psi\|_{L^2(\Omega_t^{B\pm})}^2)\quad\text{ for all }\psi\in H^1(\Omega_t^{B\pm}).
\]
\end{Lemma}
\textit{Proof.} It is enough to prove the assertion for $S:=(-\delta,\delta)\times(0,2\delta_0)$ instead of $\Omega_t^{B\pm}$ since all appearing terms are equivalent to the transformed ones under $X^\pm(.,t)$ uniformly in $t\in[0,T]$. For the $S$-case we use an idea from Evans \cite{Evans}, 5.10, problem 7. Let $w\in C^1(\overline{S})^2$ with $w|_{\partial S}\cdot\nu_{\partial S}\geq 1$. Then
\[
\|\tr\,\psi\|_{L^2(\partial S)}^2\leq\int_{\partial S}\psi^2 w\cdot\nu\,d\Hc^1=\int_S\diverg(\psi^2w)\,dx=\int_S\psi^2\diverg\,w+2\psi w\cdot\nabla\psi\,dx.
\]
Hence by Young's inequality $\|\tr\,\psi\|_{L^2(\partial S)}^2\leq C\left[\|\psi\|_{L^2(S)}^2+\varepsilon\|\nabla\psi\|_{L^2(S)}^2+\frac{1}{\varepsilon}\|\psi\|_{L^2(S)}^2\right]$. Since $1\leq\frac{\overline{\varepsilon}}{\varepsilon}$, the claim is shown.\hfill$\square$\\
\newline
\textit{Proof of Lemma \ref{th_SE_VxVperp} (Continuation).} Young's inequality and Lemma \ref{th_SE_VperpxVperp} yield
\[
|(II)|\leq\frac{\nu_0}{8\varepsilon C}\|\tr\,\psi\|_{L^2(\partial\Omega_t^{B\pm})}^2+\tilde{C}\varepsilon^3\|a\|_{H^1(0,2\delta_0)}^2\leq\frac{1}{8}B_{\varepsilon,t}^\pm(\psi,\psi)+\tilde{C}\varepsilon^3\|a\|_{H^1(0,2\delta_0)}^2.
\]
\textit{Ad $(III)$}: It holds $(III)=\int_0^{2\delta_0}a'(\sigma)g^\pm(\sigma)\,d\sigma$ with
\[
g^\pm(\sigma):=\mp\int_{-\delta}^\delta \nabla s|_{\overline{X}^\pm(r,\sigma,t)}\cdot\left[\phi_A^\varepsilon(.,t)\nabla\psi-\nabla\phi_A^\varepsilon(.,t)\psi\right]|_{X^\pm(r,\sigma,t)}\,J_t^\pm(r,\sigma)\,dr.
\]
Here $\nabla\psi|_{X^\pm(.,t)}=\nabla r|_{\overline{X}^\pm(.,t)}\partial_r\tilde{\psi}_t^\pm\mp\nabla s|_{\overline{X}^\pm(.,t)}\partial_\sigma\tilde{\psi}_t^\pm$ with $\tilde{\psi}_t^\pm:=\psi(X^\pm(.,t))$. To estimate the $\partial_\sigma\tilde{\psi}_t^\pm$-term we use that $||\nabla s|^2(\overline{X}^\pm(r,\sigma,t))-4/L(t)^2|\leq C|r|$ because of a Taylor expansion and Remark \ref{th_coord_bem}, 4. Therefore f.a.a. $\sigma\in(0,2\delta_0)$
\begin{align*}
|g^\pm(\sigma)|\leq \frac{4}{L(t)^2}\left|\int_{-\delta}^\delta\phi_A^\varepsilon|_{\overline{X}^\pm(r,\sigma,t)}(\partial_\sigma\tilde{\psi}_t^\pm J_t^\pm)|_{(r,\sigma)}\,dr\right|+\\
+\int_{-\delta}^{\delta}C|r\,\phi_A^\varepsilon|_{\overline{X}^\pm(r,\sigma,t)}(\partial_\sigma\tilde{\psi}_t^\pm J_t^\pm)|_{(r,\sigma)}|\,dr+\int_{-\delta}^{\delta}|\left[\nabla r\cdot\nabla s\,\phi_A^\varepsilon\right]|_{\overline{X}^\pm(r,\sigma,t)}(\partial_r\tilde{\psi}_t^\pm J_t^\pm)|_{(r,\sigma)}|\,dr+\\
+\int_{-\delta}^{\delta}|(\tilde{\psi}_t^\pm J_t^\pm)|_{(r,\sigma)}\left[\nabla s\cdot\nabla\phi_A^\varepsilon\right]|_{\overline{X}^\pm(r,\sigma,t)}|\,dr,
\end{align*}
where 
\begin{align*}
\nabla\phi_A^\varepsilon=\frac{1}{\sqrt{\varepsilon}}\left[\theta_0''(\rho_\varepsilon)q^\pm(s^\pm,t)+\varepsilon\partial_\rho\hat{v}^{B\pm}|_{(\rho_\varepsilon,H_\varepsilon^\pm,t)}\right]\left[\frac{\nabla r}{\varepsilon}-\nabla s\partial_sh_\varepsilon(s,t)\right]+\\
\mp\frac{1}{\sqrt{\varepsilon}}\nabla s\left[\partial_\sigma q^\pm(s^\pm,t)\theta_0'(\rho_\varepsilon)+\partial_H\hat{v}^{B\pm}|_{(\rho_\varepsilon,H_\varepsilon^\pm,t)}\right]
\end{align*}
We rewrite the first term separately with the aid of $\psi\in(V_{\varepsilon,t}^\pm)^\perp$. Because of Lemma \ref{th_prod_sets_spaces} and since integration gives a bounded linear operator on $L^2(-\delta,\delta)$, we can differentiate the identity in Lemma \ref{th_SE_split_L2}, 2. with respect to $\sigma$ and obtain
\[
\int_{-\delta}^\delta\left[\phi_A^\varepsilon|_{\overline{X}^\pm(.,t)}\partial_\sigma\tilde{\psi}_t^\pm J_t^\pm\right]|_{(r,\sigma)}\,dr=-\int_{-\delta}^\delta\left[(\partial_\sigma(\phi_A^\varepsilon|_{\overline{X}^\pm(.,t)})J_t^\pm+\phi_A^\varepsilon|_{\overline{X}^\pm(.,t)}\partial_\sigma J_t^\pm)\tilde{\psi}_t^\pm\right]|_{(r,\sigma)}\,dr,
\]
where 
\begin{align}\begin{split}\label{eq_SE_dsigma_phiA}
\partial_\sigma(\phi_A^\varepsilon|_{\overline{X}^\pm})=-\frac{1}{\sqrt{\varepsilon}}\partial_\sigma h_\varepsilon^\pm|_{(\sigma,t)}\left[\theta_0''(\rho_\varepsilon)q^\pm(s^\pm,t)+\varepsilon\partial_\rho\hat{v}^{B\pm}(\rho_\varepsilon,H_\varepsilon^\pm,t)\right]|_{\overline{X}^\pm}+\\
+\frac{1}{\sqrt{\varepsilon}}\left[\partial_H\hat{v}^{B\pm}(\rho_\varepsilon,H_\varepsilon^\pm,t)+\partial_\sigma q^\pm(s^\pm,t)\theta_0'(\rho_\varepsilon)\right]|_{\overline{X}^\pm}.
\end{split}
\end{align} 
For all terms we use Hölder's inequality and Lemma \ref{th_remainder} to get
\[
|g^\pm(\sigma)|\leq C\|\tilde{\psi}_t^\pm(.,\sigma)\|_{L^2(-\delta,\delta;J_t^\pm(.,\sigma))}+C\varepsilon\|\nabla_{(r,\sigma)}\tilde{\psi}_t^\pm(.,\sigma)\|_{L^2(-\delta,\delta;J_t^\pm(.,\sigma))}
\]
f.a.a. $\sigma\in(0,2\delta_0)$. Since $|\nabla_{(r,\sigma)}\tilde{\psi}_t^\pm|\leq C|\nabla\psi|_{X^\pm(.,t)}|$, we obtain with Young's inequality
\begin{align*}
|(III)|
&\leq C\|a'\|_{L^2(0,2\delta_0)}(\|\psi\|_{L^2(\Omega_t^{B\pm})}+\varepsilon\|\nabla\psi\|_{L^2(\Omega_t^{B\pm})})\leq\\
&\leq C\varepsilon^2\|a'\|_{L^2(0,2\delta_0)}^2+\frac{\nu_0}{8}\left[\frac{1}{\varepsilon^2}\|\psi\|_{L^2(\Omega_t^{B\pm})}^2+\|\nabla\psi\|_{L^2(\Omega_t^{B\pm})}^2\right].
\end{align*}
The last term is dominated by $\frac{1}{8}B_{\varepsilon,t}^\pm(\psi,\psi)$ because of Lemma \ref{th_SE_VperpxVperp}. Altogether we have shown Lemma \ref{th_SE_VxVperp}.\hfill$\square$\\
\newline
Finally, we put all those estimates together:
\begin{Theorem}
There are $\varepsilon_0, c_0>0$ such that 
\[
B_{\varepsilon,t}^\pm(\psi,\psi)\geq -C\|\psi\|_{L^2(\Omega_t^{B\pm})}^2+c_0\|\nabla_\tau\psi\|_{L^2(\Omega_t^{B\pm})}^2
\] 
for all $\varepsilon\in(0,\varepsilon_0], t\in[0,T]$ and $\psi\in H^1(\Omega_t^{B\pm})$ with $\psi|_{X(.,s,t)}=0$ for $s\in(\frac{3}{2}\delta_0,2\delta_0)$.
\end{Theorem}
\begin{Remark}\upshape
The estimate can be refined, cf. the proof below.
\end{Remark}
\textit{Proof.} Because of Lemma \ref{th_SE_split_L2}, any $\psi\in H^1(\Omega_t^{B\pm})$ can be uniquely written as $\psi=\phi+\phi^\perp$ with $\phi=\left[a(s^\pm)\phi_A^\varepsilon\right]|_{(.,t)}\in V_{\varepsilon,t}^\pm$ and $\phi^\perp\in(V_{\varepsilon,t}^\pm)^\perp$. Lemma \ref{th_SE_VxV} and Lemma \ref{th_SE_VxVperp} imply for $t\in[0,T]$ and $\varepsilon\in(0,\varepsilon_0]$, $\varepsilon_0>0$ small that
\begin{align*}
B_{\varepsilon,t}^\pm(\psi,\psi)=B_{\varepsilon,t}^\pm(\phi,\phi)+2B_{\varepsilon,t}^\pm(\phi,\phi^\perp)+B_{\varepsilon,t}^\pm(\phi^\perp,\phi^\perp)\geq\\
\geq -C\|\phi\|_{L^2(\Omega_t^{B\pm})}^2+(c_0-C\varepsilon^2)\|a\|_{H^1(0,2\delta_0)}^2-\frac{C}{\varepsilon}\|\phi\|_{L^2(\Omega_t^{B\pm})}\|\phi^\perp\|_{L^2(\Omega_t^{B\pm})}+\frac{B_{\varepsilon,t}^\pm(\phi^\perp,\phi^\perp)}{2}.
\end{align*}
For the last part we use Lemma \ref{th_SE_VperpxVperp} and we estimate the third part with Young's inequality to $\frac{\nu_0}{4\varepsilon^2}\|\phi^\perp\|_{L^2(\Omega_t^{B\pm})}^2+C\|\phi\|_{L^2(\Omega_t^{B\pm})}^2$. Hence
\[
B_{\varepsilon,t}^\pm(\psi,\psi)\geq-C\|\phi\|_{L^2(\Omega_t^{B\pm})}^2+\frac{\nu_0}{4\varepsilon^2}\|\phi^\perp\|_{L^2(\Omega_t^{B\pm})}^2+\frac{c_0}{2}\|a\|_{H^1(0,2\delta_0)}^2+\frac{\nu_0}{2}\|\nabla(\phi^\perp)\|_{L^2(\Omega_t^{B\pm})}^2.
\]
It remains to include the $\nabla_\tau\psi$-term in the estimate. The triangle inequality implies $\|\nabla_\tau\psi\|_{L^2(\Omega_t^{B\pm})}\leq\|\nabla_\tau\phi\|_{L^2(\Omega_t^{B\pm})}+\|\nabla_\tau(\phi^\perp)\|_{L^2(\Omega_t^{B\pm})}$. By definition 
\[
\mp\nabla_\tau\phi|_{X^\pm(.,t)}=\nabla s|_{\overline{X}^\pm(.,t)}\partial_\sigma(\phi|_{X^\pm(.,t)})=\nabla s|_{\overline{X}^\pm(.,t)}\left[a'(\sigma)\phi_A^\varepsilon|_{\overline{X}^\pm(.,t)}+a(\sigma)\partial_\sigma(\phi_A^\varepsilon|_{\overline{X}^\pm(.,t)})\right],
\]
where $\partial_\sigma(\phi_A^\varepsilon|_{\overline{X}^\pm})$ was computed in \eqref{eq_SE_dsigma_phiA}. With Fubini's Theorem and Lemma \ref{th_remainder} one shows $\|\nabla_\tau\phi\|_{L^2(\Omega_t^{B\pm})}\leq C\|a\|_{H^1(0,2\delta_0)}$. Since $\|\nabla_\tau(\phi^\perp)\|_{L^2(\Omega_t^{B\pm})}\leq C\|\nabla(\phi^\perp)\|_{L^2(\Omega_t^{B\pm})}$, 
we obtain
\[
\|\nabla_\tau\psi\|_{L^2(\Omega_t^{B\pm})}^2\leq C(\|a\|_{H^1(0,2\delta_0)}^2+\|\nabla(\phi^\perp)\|_{L^2(\Omega_t^{B\pm})}^2).
\]
This yields the desired estimate.\hfill$\square$\\
\newline
The proof of Theorem \ref{th_SE_bdry} is completed.
\section{Difference Estimate}\label{sec_difference}
In this section we estimate the difference of the exact and approximate solution. From this we obtain a convergence result in the next section by showing that our construction fulfills the requirements.
\begin{Theorem}\label{th_difference_estimate}
Let $N=2$, $\Gamma=(\Gamma_t)_{t\in[0,T_0]}$ for some $T_0>0$ be as in Section \ref{sec_coord} and $\delta>0$ such that Theorem \ref{th_coord} holds for $2\delta$ instead of $\delta$. For an $\varepsilon_0>0$ let $u_A^\varepsilon:\overline{\Omega}\times[0,T_0]\rightarrow\R$ and $u_0^\varepsilon:\overline{\Omega}\rightarrow\R$ for $\varepsilon\in(0,\varepsilon_0]$ be smooth with $\partial_\nu u_0^\varepsilon=0$ on $\partial\Omega$ and let $u^\varepsilon:\overline{\Omega}\times[0,T_0]\rightarrow\R$ be the smooth exact solution\footnote{ For analytical properties of the Allen-Cahn equation cf. Bartels \cite{Bartels_book}, Chapter 6.1 and references therein.} to \eqref{eq_AC1}-\eqref{eq_AC3} with $u_0^\varepsilon$ as initial values. For some $R>0$ and $M\in\N, M\geq N$ we impose the following conditions:
\begin{enumerate}
\item Boundedness: $\sup_{\varepsilon\in(0,\varepsilon_0]}\|u_A^\varepsilon\|_{L^\infty(\Omega_{T_0})}+\|u_0^\varepsilon\|_{L^\infty(\Omega)}<\infty$.
\item Spectral Estimate: There are $c_0,C>0$ such that
\[
\int_\Omega|\nabla\psi|^2+\frac{1}{\varepsilon^2}f''(u_A^\varepsilon(.,t))\psi^2\,dx\geq -C\|\psi\|_{L^2(\Omega)}^2+\|\nabla\psi\|_{L^2(\Omega\textbackslash\Gamma_t(\delta))}^2+c_0\|\nabla_\tau\psi\|_{L^2(\Gamma_t(\delta))}^2
\]
for all $\psi\in H^1(\Omega)$ and $\varepsilon\in(0,\varepsilon_0],t\in[0,T_0]$.
\item $(u_A^\varepsilon)_{\varepsilon\in(0,\varepsilon_0]}$ is an approximate solution to \eqref{eq_AC1}-\eqref{eq_AC3} in the following sense: It holds $\partial_\nu u_A^\varepsilon=0$ on $\partial\Omega_{T_0}\textbackslash\Gamma(2\delta)$ for $\varepsilon\in(0,\varepsilon_0]$ and 
\[
\|\partial_\nu u_A^\varepsilon(t)\|_{L^2(\partial\Omega)}\leq C\varepsilon^{M+\frac{1}{2}}\quad\text{ for all }t\in[0,T_0],\quad \|u_A^\varepsilon|_{t=0}-u_0^\varepsilon\|_{L^2(\Omega)}\leq R\varepsilon^{M+\frac{1}{2}}
\]
and for the difference $\overline{u}^\varepsilon:=u^\varepsilon-u_A^\varepsilon$ and the remainder $r_A^\varepsilon:=\partial_t u_A^\varepsilon-\Delta u_A^\varepsilon+\frac{1}{\varepsilon^2}f'(u_A^\varepsilon)$ it holds 
\begin{align}\label{eq_rA_estimate}
\left|\int_\Omega r_A^\varepsilon\overline{u}^\varepsilon(t)\,dx\right|\leq C\varepsilon^{M+\frac{1}{2}}(\|\overline{u}^\varepsilon(t)\|_{L^2(\Omega)}+\|\nabla_\tau\overline{u}^\varepsilon(t)\|_{L^2(\Gamma_t(2\delta))})
\end{align}
for all $\varepsilon\in(0,\varepsilon_0],T\in(0,T_0]$.
\end{enumerate} 
Then we obtain
\begin{enumerate} 
\item Let $M\geq 3$. Then there are $\alpha,\varepsilon_1>0$ such that for $g_\alpha:[0,T_0]\rightarrow\R:t\mapsto e^{-\alpha t}$ it holds
\begin{align}
\begin{split}
\sup_{t\in[0,T]}\|g_\alpha\overline{u}^\varepsilon(t)\|_{L^2(\Omega)}^2+\|g_\alpha\nabla\overline{u}^\varepsilon\|_{L^2(Q_T\textbackslash\Gamma(\delta))}^2&\leq 2R^2\varepsilon^{2M+1},\\
c_0\|g_\alpha\nabla_\tau\overline{u}^\varepsilon\|^2_{L^2(Q_T\cap\Gamma(\delta))}+\varepsilon^2\|g_\alpha\partial_n\overline{u}^\varepsilon\|^2_{L^2(Q_T\cap\Gamma(\delta))}&\leq 2R^2\varepsilon^{2M+1}\label{eq_diff_estimate_eps_R_T}\end{split}
\end{align}
for all $\varepsilon\in(0,\varepsilon_1]$ and $T\in(0,T_0]$.
\item Let $M=N=2$. Then there are $\varepsilon_1,T_1>0$ such that $\eqref{eq_diff_estimate_eps_R_T}$ holds for $\alpha=0$ and all $\varepsilon\in(0,\varepsilon_1], T\in(0,T_1]$.
\end{enumerate}
\end{Theorem}
\textit{Proof.} By continuity and the Dominated Convergence Theorem we obtain that 
\[
T_{\varepsilon,\alpha}:=\sup\,\{T\in(0,T_0]: \eqref{eq_diff_estimate_eps_R_T}\text{ holds for all }t\in(0,T]\}
\]
is well-defined for all $\varepsilon\in(0,\varepsilon_0],\alpha\geq 0$ and it holds $T_{\varepsilon,\alpha}>0$. In the case $M\geq 3$ we have to show that there are $\varepsilon_1,\alpha>0$ such that $T_{\varepsilon,\alpha}=T_0$ for all $\varepsilon\in(0,\varepsilon_1]$. If $M=N$ we prove $T_{\varepsilon,0}\geq T_1$ for all $\varepsilon\in(0,\varepsilon_1]$ and some $T_1>0$.\\
We take the difference of the left hand side of the Allen-Cahn equation for $u^\varepsilon$ and $u_A^\varepsilon$ and obtain
\begin{align}\label{eq_diff_proof1}
\partial_t\overline{u}^\varepsilon-\Delta\overline{u}^\varepsilon+\frac{1}{\varepsilon^2}f''(u_A^\varepsilon)\overline{u}^\varepsilon=-r_A^\varepsilon-r^\varepsilon(u^\varepsilon,u_A^\varepsilon),
\end{align} 
where $r^\varepsilon(u^\varepsilon,u_A^\varepsilon):=\frac{1}{\varepsilon^2}\left[f'(u^\varepsilon)-f'(u_A^\varepsilon)-f''(u_A^\varepsilon)\overline{u}^\varepsilon\right]$. We multiply this equation by $g_\alpha^2\overline{u}^\varepsilon$ and integrate over $Q_T$ for $T\in(0,T_{\varepsilon,\alpha}]$, where $\varepsilon\in(0,\varepsilon_0]$ and $\alpha\geq 0$ are fixed. The first term gives
\[
\int_0^T\int_\Omega g_\alpha^2\partial_t\overline{u}^\varepsilon\overline{u}^\varepsilon\,dx\,dt=\frac{1}{2}\|g_\alpha(T)\overline{u}^\varepsilon(T)\|_{L^2(\Omega)}^2-\frac{1}{2}\|\overline{u}^\varepsilon(0)\|_{L^2(\Omega)}^2+\int_0^T\partial_tg_\alpha g_\alpha\|\overline{u}^\varepsilon\|_{L^2(\Omega)}^2\,dt,
\]
where $\|\overline{u}^\varepsilon(0)\|_{L^2(\Omega)}^2\leq R^2\varepsilon^{2M+1}$. For the other term coming from the left hand side in \eqref{eq_diff_proof1} we use integration by parts:
\begin{align*}
\int_0^Tg_\alpha^2\int_\Omega\left[-\Delta\overline{u}^\varepsilon+\frac{1}{\varepsilon^2}f''(u_A^\varepsilon)\overline{u}^\varepsilon\right]\overline{u}^\varepsilon\,dx\,dt&=\int_0^Tg_\alpha^2\int_\Omega|\nabla\overline{u}^\varepsilon|^2+\frac{1}{\varepsilon^2}f''(u_A^\varepsilon)(\overline{u}^\varepsilon)^2\,dx\,dt+\\
&-\int_0^Tg_\alpha^2\int_{\partial\Omega}\partial_\nu u_A^\varepsilon\,\tr\,\overline{u}^\varepsilon\,d\Hc^1\,dt.
\end{align*}
Because of the spectral estimate and \eqref{eq_diff_estimate_eps_R_T}, the first term is bounded from below by
\[
-C\int_0^T g_\alpha^2\|\overline{u}^\varepsilon\|_{L^2(\Omega)}^2\,dt+\|g_\alpha\nabla\overline{u}^\varepsilon\|_{L^2(Q_T\textbackslash\Gamma(\delta))}^2+c_0\|g_\alpha\nabla_\tau\overline{u}^\varepsilon\|_{L^2(Q_T\cap\Gamma(\delta))}^2.
\]
The boundary term can be estimated as follows: Since $\partial_\nu u_A^\varepsilon=0$ on $\partial Q_T\textbackslash\Gamma(2\delta)$, we infer
\[
\left|\int_{\partial\Omega}\partial_\nu u_A^\varepsilon\,\tr\,\overline{u}^\varepsilon(t)\,d\Hc^1\right|\leq\|\partial_\nu u_A^\varepsilon(t)\|_{L^2(\partial\Omega)}\|\tr\,\overline{u}^\varepsilon(t)\|_{L^2(\partial\Omega\cap\Gamma_t(2\delta))}.
\]
With a similar idea as in the proof of Lemma \ref{th_inpol_trace} we obtain
\[
\|\tr\,\overline{u}^\varepsilon(t)\|_{L^2(\partial\Omega\cap\Gamma_t(2\delta))}^2\leq C\|\overline{u}^\varepsilon(t)\|_{L^2(\Gamma_t(2\delta))}^2+C\|\overline{u}^\varepsilon(t)\|_{L^2(\Gamma_t(2\delta))}\|\nabla_\tau\overline{u}^\varepsilon(t)\|_{L^2(\Gamma_t(2\delta))}.
\]
To this end one uses $S=(-2\delta,2\delta)\times(-1,1)$ and $w\in C^1(\overline{S})^2$ such that $w_1=0$ and $w|_{\partial S}\cdot\nu_{\partial S}\geq 1$ on $[-2\delta,2\delta]\times\{\pm1\}$ and $w|_{\partial S}\cdot\nu_{\partial S}=0$, else. Therefore we obtain
\begin{align*}
\left|\int_0^Tg_\alpha^2\int_{\partial\Omega}\partial_\nu u_A^\varepsilon\,\tr\,\overline{u}^\varepsilon\,d\Hc^1\,dt\right|\leq 
C\varepsilon^{M+\frac{1}{2}}\int_0^T g_\alpha^2(\|\overline{u}^\varepsilon\|_{L^2(\Gamma_t(2\delta))}+\|\nabla_\tau\overline{u}^\varepsilon\|_{L^2(\Gamma_t(2\delta))})\,dt.
\end{align*}
Due to $|\nabla_\tau\overline{u}^\varepsilon|\leq C|\nabla\overline{u}^\varepsilon|$ and \eqref{eq_diff_estimate_eps_R_T}, this is estimated by $CR(\|g_\alpha\|_{L^1(0,T)}+\|g_\alpha\|_{L^2(0,T)})\varepsilon^{2M+1}$.\\
It remains to treat the terms coming from the right hand side in \eqref{eq_diff_proof1}. For the one involving $r_A^\varepsilon$ we use \eqref{eq_rA_estimate}. This term has the same structure as the one we obtained from estimating the boundary term. The other term coming from $r^\varepsilon(u^\varepsilon,u_A^\varepsilon)$ can be treated in an analogous way as in Abels, Liu \cite{ALiu}, Lemma 5.3: As in the proof of the weak maximum principle it follows that
\[
\sup_{\varepsilon\in(0,\varepsilon_0]}\|u^\varepsilon\|_{L^\infty(Q_T)}\leq\max\{1,\sup_{\varepsilon\in(0,\varepsilon_0]}\|u_0^\varepsilon\|_{L^\infty(\Omega)}\}.
\] 
Therefore the boundedness assumptions together with a Taylor expansion yield 
\[
\left|\int_0^Tg_\alpha^2\int_\Omega r^\varepsilon(u^\varepsilon,u_A^\varepsilon)\overline{u}^\varepsilon\,dx\,dt\right|
\leq \frac{C}{\varepsilon^2}\int_0^Tg_\alpha^2\|\overline{u}^\varepsilon\|_{L^3(\Omega)}^3\,dt.
\]
The Gagliardo-Nirenberg inequality in two dimensions applied to $|\overline{u}^\varepsilon(t)|^2$ implies
\[
\|\overline{u}^\varepsilon(t)\|_{L^3(\Omega\textbackslash\Gamma_t(\delta))}^3\leq C\|\overline{u}^\varepsilon\|_{H^1(\Omega\textbackslash\Gamma_t(\delta))}\|\overline{u}^\varepsilon\|_{L^2(\Omega\textbackslash\Gamma_t(\delta))}^2\quad\text{ for all }t\in[0,T].
\]
To control the $L^3$-norm of $\overline{u}^\varepsilon$ on $\Gamma_t(\delta)$, one transforms to $(-\delta,\delta)\times I$, uses Lemma \ref{th_prod_sets_spaces}, 1. and one-dimensional Gagliardo-Nirenberg inequalities in both coordinates\footnote{ This is similar to Lemma 5.4 in Abels, Liu \cite{ALiu}. But there the one-dimensional Gagliardo-Nirenberg inequality was only used in normal direction. Here it is also not needed and just added for optimality.}. This yields
\[
\|\overline{u}^\varepsilon\|_{L^3(\Gamma_t(\delta))}^3
\leq C\|(\overline{u}^\varepsilon,\partial_n\overline{u}^\varepsilon)\|_{L^2(\Gamma_t(\delta))}^{\frac{1}{2}}\|(\overline{u}^\varepsilon,\nabla_\tau\overline{u}^\varepsilon)\|_{L^2(\Gamma_t(\delta))}^{\frac{1}{2}}\|\overline{u}^\varepsilon\|_{L^2(\Gamma_t(\delta))}^2\quad\text{ for all }t\in[0,T].
\]
The appearing terms can be estimated with \eqref{eq_diff_estimate_eps_R_T} and Hölder's inequality. We obtain
\begin{align*}
\left|\int_0^T g_\alpha^2\int_\Omega r^\varepsilon(u^\varepsilon,u_A^\varepsilon)\overline{u}^\varepsilon\,dx\,dt\right|\leq\frac{C}{\varepsilon^2}R^3\left[\|g_\alpha^{-1}\|_{L^1(0,T)}+\|g_\alpha^{-1}\|_{L^2(0,T)}\right]\varepsilon^{3(M+\frac{1}{2})}+\\
+\frac{C}{\varepsilon^2}R^3\left[\varepsilon^\frac{1}{2}\|g_\alpha^{-1}\|_{L^1(0,T)}+\|g_\alpha^{-1}\|_{L^{\frac{4}{3}}(0,T)}+\|g_\alpha^{-1}\|_{L^2(0,T)}\right]\varepsilon^{3M+1}.
\end{align*}
Moreover,
\begin{align*}
\varepsilon^2\|g_\alpha\partial_n\overline{u}^\varepsilon\|_{L^2(Q_T\cap\Gamma(\delta))}^2
&\leq C\varepsilon^2\int_0^T g_\alpha^2\int_\Omega|\nabla\overline{u}^\varepsilon|^2+\frac{1}{\varepsilon^2}f''(u_A^\varepsilon)(\overline{u}^\varepsilon)^2\,dx\,dt+\\
&+C\sup_{\varepsilon\in(0,\varepsilon_0]}\|f''(u_A^\varepsilon)\|_{L^\infty(\Omega_{T_0})}\int_0^T g_\alpha^2\|\overline{u}^\varepsilon(t)\|_{L^2(\Omega)}^2\,dt.
\end{align*}
Altogether, if $\varepsilon_0>0$ is small we get for all $T\in(0,T_{\varepsilon,\alpha}], \varepsilon\in(0,\varepsilon_0]$ that
\begin{align}\notag
&\frac{1}{2}\|g_\alpha(T)\overline{u}^\varepsilon(T)\|_{L^2(\Omega)}^2+\frac{1}{2}\|g_\alpha\nabla\overline{u}^\varepsilon\|_{L^2(Q_T\textbackslash\Gamma(\delta))}^2+\\\notag
&+\frac{c_0}{2}\|g_\alpha\nabla_\tau\overline{u}^\varepsilon\|_{L^2(Q_T\cap\Gamma(\delta))}^2+\frac{\varepsilon^2}{2}\|g_\alpha\partial_n\overline{u}^\varepsilon\|_{L^2(Q_T\cap\Gamma(\delta))}^2\leq\\\notag
&\leq \frac{R^2}{2}\varepsilon^{2M+1}+\int_0^T(-\alpha+\tilde{C}_0)g_\alpha^2\|\overline{u}^\varepsilon(t)\|_{L^2(\Omega)}^2\,dt+C_1R\varepsilon^{2M+1}\|g_\alpha\|_{L^2(0,T)}+\\
&+CR^3\varepsilon^{2M+1}\left[\varepsilon^{\frac{1}{2}}\|g^{-1}\|_{L^2(0,T)}+\varepsilon^{M-2}\|g^{-1}\|_{L^2(0,T)}\right].\label{eq_diff_proof2}
\end{align}
Now let $M\geq 3$. Then we choose $\alpha\geq \tilde{C}_0$ large such that $C_1R\|g_\alpha\|_{L^2(0,T)}\leq\frac{R^2}{8}$. Therefore the left hand side in \eqref{eq_diff_proof2} is estimated by $\frac{3}{4}R^2\varepsilon^{2M+1}$ for all $T\in (0,T_{\varepsilon,\alpha}]$ and $\varepsilon\in(0,\varepsilon_1]$, if $\varepsilon_1>0$ is small enough. By definition and continuity this shows $T_{\varepsilon,\alpha}=T_0$ for all $\varepsilon\in(0,\varepsilon_1]$.\\
Finally, we consider $M=N=2$ and $\alpha=0$. Then the left hand side in \eqref{eq_diff_proof2} is dominated by
\[
\left[\frac{R^2}{2}+CR^2T+CRT^{\frac{1}{2}}+CR^3(\varepsilon^\frac{1}{2}+1)T^\frac{1}{2}\right]\varepsilon^{2N+1}.
\]
There are $\varepsilon_1,T_1>0$ such that this is estimated by $\frac{3}{4}R^2\varepsilon^{2N+1}$ for all $T\in(0,\min(T_{\varepsilon,\alpha},T_1)]$ and $\varepsilon\in(0,\varepsilon_1]$. Therefore by definition and continuity $T_{\varepsilon,0}\geq T_1$ for all $\varepsilon\in(0,\varepsilon_1]$.\hfill$\square$\\
\section{Convergence Result}\label{sec_conv_result}
Let $N=2$, $\Gamma=(\Gamma_t)_{t\in[0,T_0]}$ for some $T_0>0$ be a smooth solution to mean curvature flow with $90$\textup{°}-contact angle condition parametrized as in Section \ref{sec_coord} starting from $\Gamma_0$ and let $\delta>0$ such that Theorem \ref{th_coord} holds for $2\delta$. We obtain the following convergence result:
\begin{Theorem}
There is an $\varepsilon_0>0$ and $u_A^\varepsilon:\overline{\Omega}\times[0,T_0]\rightarrow\R$ for $\varepsilon\in(0,\varepsilon_0]$ smooth such that the following assertion holds: If $u_0^\varepsilon:\overline{\Omega}\rightarrow\R$ for $\varepsilon\in(0,\varepsilon_0]$ are smooth with $\partial_\nu u_0^\varepsilon=0$ and for some $R>0$ it holds
\[
\sup_{\varepsilon\in(0,\varepsilon_0]}\|u_0^\varepsilon\|_{L^\infty(\Omega)}<\infty\quad\text{ and }\quad\|u_0^\varepsilon-u_A^\varepsilon|_{t=0}\|_{L^2(\Omega)}\leq R\varepsilon^{N+\frac{1}{2}}\text{ for }\varepsilon\in(0,\varepsilon_0],
\] 
then there are $C,\varepsilon_1,T_1>0$ such that for the unique smooth solution $u^\varepsilon:\overline{\Omega}\times[0,T_0]\rightarrow\R$ to \eqref{eq_AC1}-\eqref{eq_AC3} with initial values $u_0^\varepsilon$ and all $\varepsilon\in(0,\varepsilon_1], T\in(0,T_1]$ it holds
\begin{align*}
\sup_{t\in[0,T]}\|u^\varepsilon(t)-u_A^\varepsilon(t)\|_{L^2(\Omega)}+\|\nabla(u^\varepsilon-u_A^\varepsilon)\|_{L^2(Q_T\textbackslash\Gamma(\delta))}&\leq C\varepsilon^{N+\frac{1}{2}},\\
\|\nabla_\tau(u^\varepsilon-u_A^\varepsilon)\|_{L^2(Q_T\cap\Gamma(\delta))}+\varepsilon\|\partial_n(u^\varepsilon-u_A^\varepsilon)\|_{L^2(Q_T\cap\Gamma(\delta))}&\leq C\varepsilon^{N+\frac{1}{2}}.
\end{align*}
\end{Theorem}
\begin{Remark}\upshape
\begin{enumerate}
\item See Section \ref{sec_approx_sol} for the explicit definition of $u_A^\varepsilon$. In particular it holds $\lim_{\varepsilon\rightarrow 0}u_A^\varepsilon=\pm 1$ uniformly on compact subsets of $Q_T^\pm$.
\item One should also obtain a convergence result for global time $T\in(0,T_0]$ if $N$ in the estimates is replaced by any $M\in\N,M\geq 3$. Therefore only a better approximate solution $u_A^\varepsilon$ is needed, cf. Theorem \ref{th_difference_estimate}. The latter should be obtained by iterating the construction in Section \ref{sec_expansion}, e.g. similar to Chen, Hilhorst, Logak \cite{CHL}.
\end{enumerate}
\end{Remark}
\textit{Proof.} We apply Theorem \ref{th_difference_estimate}. The $u_A^\varepsilon$ for $\varepsilon>0$ are constructed in Section \ref{sec_approx_sol}. They are smooth, uniformly bounded in $\varepsilon$ and fulfill the spectral estimate in Theorem \ref{th_difference_estimate}, 2. because of Theorem \ref{th_SE} for $\varepsilon\in(0,\varepsilon_0]$ if $\varepsilon_0>0$ is small. As in the proof of Lemma \ref{th_remainder} it follows that for $\varepsilon_0>0$ small and all $\varepsilon\in(0,\varepsilon_0]$ it holds
\[ 
\|\partial_\nu u_A^\varepsilon(t)\|_{L^2(\partial\Omega)}\leq C\varepsilon^{N+\frac{1}{2}}\quad\text{ for all }t\in[0,T_0],
\]
where we used Lemma \ref{th_approx_lsq}. Therefore the claim follows from Theorem \ref{th_difference_estimate} if we show \eqref{eq_rA_estimate} for $M=N$. Due to Lemma \ref{th_approx_lsq} we have
\begin{align}\label{eq_conv_result_proof1}
|r_A^\varepsilon|\leq C\left[\varepsilon\sum_\pm e^{-c(|\rho_\varepsilon|+ H_\pm^\varepsilon)}+\varepsilon^2e^{-c|\rho_\varepsilon|}+\varepsilon^3\right]\quad\text{ in }\Gamma(2\delta),\quad
r_A^\varepsilon=0\quad\text{ in }\Omega\textbackslash\Gamma(2\delta).
\end{align}
Hence
\begin{align*}
\left|\int_\Omega r_A^\varepsilon\overline{u}^\varepsilon(t)\,dx\right|
\leq\|r_A^\varepsilon\overline{u}^\varepsilon\|_{L^1(\Omega)}
\leq\int_{-1}^1\int_{-2\delta}^{2\delta}|r_A^\varepsilon\overline{u}^\varepsilon||_{\overline{X}(r,s,t)}J_t(r,s)\,dr\,ds.
\end{align*}
The inner integral is dominated by $C\|\overline{u}^\varepsilon|_{\overline{X}(.,s,t)}\|_{L^2(-2\delta,2\delta)}\|r_A^\varepsilon|_{\overline{X}(.,s,t)}\|_{L^2(-2\delta,2\delta)}$. With \eqref{eq_conv_result_proof1} it follows as in the proof of Lemma \ref{th_remainder} that
\[
\|r_A^\varepsilon|_{\overline{X}(.,s,t)}\|_{L^2(-2\delta,2\delta)}\leq C\left[\varepsilon^\frac{3}{2}\sum_\pm e^{-\frac{c}{\varepsilon}(|s\pm 1|)}+\varepsilon^{\frac{5}{2}}\right]
\]
for $\varepsilon\in(0,\varepsilon_0]$ if $\varepsilon_0>0$ is small. Therefore because of Lemma \ref{th_prod_sets_spaces}
\begin{align*}
\|r_A^\varepsilon\overline{u}^\varepsilon(t)\|_{L^1(\Omega)}
&\leq\|\overline{u}^\varepsilon|_{\overline{X}(.,t)}\|_{L^\infty(-1,1;L^2(-2\delta,2\delta))}\int_{-1}^1 C\left[\varepsilon^\frac{3}{2}\sum_\pm e^{-\frac{c}{\varepsilon}(|s\pm 1|)}+\varepsilon^{\frac{5}{2}}\right]\,ds\leq\\
&\leq C(\|\overline{u}^\varepsilon\|_{L^2(\Gamma_t(2\delta))}+\|\nabla_\tau\overline{u}^\varepsilon\|_{L^2(\Gamma_t(2\delta))})\varepsilon^{N+\frac{1}{2}},
\end{align*}
where we used a scaling argument for the second part.\\
Theorem \ref{th_difference_estimate} implies the existence of $\varepsilon_1,T_1>0$ such that \eqref{eq_diff_estimate_eps_R_T} holds for $\alpha=0$ and all $\varepsilon\in(0,\varepsilon_1],T\in(0,T_1]$.\hfill$\square$\\
\setcounter{secnumdepth}{0}\setcounter{tocdepth}{0}


\end{document}